\RequirePackage{fix-cm}
\documentclass[smallextended]{svjour3}       % onecolumn (second format)
\usepackage{amsmath,amsfonts,amssymb}
\usepackage{mathtools}  
\usepackage{breqn}
\usepackage[dvipsnames]{xcolor}
\usepackage{multirow}
\usepackage{hyperref}
\usepackage{enumerate}
\usepackage{subfigure}
\usepackage[title]{appendix}
\usepackage{lineno}
\usepackage[normalem]{ulem}
\usepackage{lmodern}
\usepackage{cite}
\usepackage[textsize=tiny]{todonotes}
\usepackage{float}
\usepackage[sort,numbers]{natbib}

\hypersetup{
	colorlinks=true, % make the links colored
	linkcolor=blue, % color TOC links in blue
	urlcolor=red, % color URLs in red
	citecolor=magenta,
	linktoc=all % 'all' will create links for everything in the TOC
}
\usepackage{caption} % Add the caption package
\usepackage{float}
\usepackage{comment}
\newcommand{\al}{\alpha}
\newcommand{\ub}{\mathbf{u}}
\newcommand{\qb}{\mathbf{q}}
\newcommand{\sbb}{\mathbf{s}}
\newcommand{\Ub}{\mathbf{U}}
\newcommand{\Wb}{\mathbf{W}}
\newcommand{\fb}{\mathbf{f}}
\newcommand{\Db}{\mathbf{D}}
\newcommand{\Fb}{\mathbf{F}}
\newcommand{\Eb}{\mathbf{E}}
\newcommand{\Bb}{\mathbf{B}}
\newcommand{\Vb}{\mathbf{V}}
\newcommand{\jb}{\mathbf{j}}
\newcommand{\pa}{\partial}
\newcommand{\na}{\nabla}
\newcommand{\eps}{\epsilon}
\newcommand{\Dx}{\Delta x}
\newcommand{\Dy}{\Delta y}
\newcommand{\Dt}{\Delta t}
\newcommand{\iph}{i+\frac{1}{2}}
\newcommand{\jph}{j+\frac{1}{2}}
\newcommand{\imh}{i-\frac{1}{2}}
\newcommand{\jmh}{j-\frac{1}{2}}

\definecolor{lightyellow}{HTML}{fffcf1}

\newcommand{\dx}{\delta_x}
\newcommand{\dy}{\delta_y}
\newcommand{\ax}{\mu_x}
\newcommand{\ay}{\mu_y}
\providecommand{\keywords}[1]
{
	\small	
	\textbf{{Keywords: }} #1
}
\begin{document}
\title{Second order divergence constraint preserving entropy stable finite difference schemes for ideal two-fluid plasma flow equations}
\author{Jaya Agnihotri \and
	Deepak Bhoriya  \and
	Harish Kumar    \and
	Praveen Chandrashekhar \and
	Dinshaw S. Balsara
}

%\author[a]{Jaya Agnihotri}
%\author[b]{Deepak Bhoriya}
%\author[a]{Harish Kumar}
%\author[c]{Praveen Chandrashekar}
%\author[b,d]{Dinshaw S. Balsara}
%\affil[a]{Department of Mathematics, Indian Institute of Technology Delhi, India}
%\affil[b]{Physics Department, University of Notre Dame, USA}
%\affil[c]{Centre for Applicable Mathematics, TIFR, Bangalore, India}
%\affil[d]{ACMS, University of Notre Dame, USA}

\institute{J. Agnihotri (Corresponding author) \at
	{Department of Mathematics, \\Indian Institute of Technology Delhi, New Delhi, India}
	\\ \email{jayaagnihotri96@gmail.com}
	%-----
	\and
	D. Bhoriya \at{Physics Department, University of Notre Dame, USA}
	\\ \email{dbhoriy2@nd.edu} 
%------
	\and
	% ----- Another author
	H. Kumar \at
	{Department of Mathematics, \\Indian Institute of Technology Delhi, New Delhi, India}
	\\ \email{hkumar@iitd.ac.in}
	% ----- 
	\and
	% ----- Another author
	P. Chandrashekhar \at
	{Centre for Applicable Mathematics, \\Tata Institute of Fundamental Research, Bangalore, India}
	\\ \email{praveen@tifrbng.res.in}
		\and
	% ----- Another author
	Dinshaw S. Balsara \at
	{Physics Department University of Notre Dame, USA\\ ACMS, University of Notre Dame, USA}
	\\ \email{dbalsara@nd.edu}
}

\date{Received: date / Accepted: date}
%\linenumbers

\maketitle

%	\tableofcontents
	\begin{abstract}
		Two-fluid plasma flow equations describe the flow of ions and electrons with different densities, velocities, and pressures. We consider the ideal plasma flow i.e. we ignore viscous, resistive and collision effects. The resulting system of equations has flux consisting of three independent components, one for ions, one for electrons, and a linear Maxwell's equation flux for the electromagnetic fields. The coupling of these components is via source terms. In this article, we present {conservative} second-order finite difference schemes that ensure the consistent evolution of the divergence constraints on the electric and magnetic fields. The key idea is to design a numerical solver for Maxwell's equations using the multidimensional Riemann solver at the vertices, ensuring discrete divergence constraints; for the fluid parts, we use an entropy-stable discretization. The proposed schemes are co-located, second-order accurate, entropy stable, and ensure divergence-free evolution of the magnetic field. We use explicit and IMplicit-EXplicit (IMEX) schemes for time discretizations. To demonstrate the accuracy, stability, and divergence constraint-preserving ability of the proposed schemes, we present several test cases in one and two dimensions. We also compare the numerical results with those obtained from schemes with no divergence cleaning and those employing perfectly hyperbolic Maxwell (PHM) equations-based divergence cleaning methods for Maxwell's equations.
	\end{abstract}
%	\linenumbers
	%===========================================================================
\section{Introduction}
Plasma flows are essential to study several interesting physical phenomena and applications, e.g., space propulsion, controlled nuclear fusion, modeling of the solar atmosphere, etc. In addition, most astrophysical phenomena need to consider plasma flows for accurate modeling. Hence, the study and stable simulation of plasma flows is a highly active area of research. One of the simplest plasma flow models is the equations of Magnetohydrodynamics (MHD)~\cite{goedbloed_poedts_2004}. The MHD model is highly successful in a wide variety of plasma flows; however, the MHD model treats plasma as a single fluid while ignoring the presence of ions and electrons. This makes the MHD model unsuitable for several interesting applications, and a more accurate flow description is needed.

Several extensions of the MHD model are based on the two-fluid description of matter~\cite{Baboolal2001, Shumlak2003, Hakim2006, loverich2011, kumar2012entropy, bond_wheatley_samtaney_pullin_2017, amano2015divergence, Abgrall2014, Meena2019, wang2020, polak_fourth-order_2023}, where ions and electrons are considered as the species. The two-fluid plasma flow model allows different densities, velocities, and pressures for each fluid species, and each fluid is described by Euler equations of compressible flow. These equations are coupled via source terms, which model the Lorentz force terms. In addition, Maxwell's equations are considered to model magnetic and electric field evolution. { This allows the modelling of the plasma without the assumption of quasi-neutrality and zero electron mass. Hence, model is able to capture several non-MHD effects e.g. Hall effects, electron inertia effects, displacement current, magnetic reconnection etc.} The complete system of PDEs is a system of Hyperbolic Balance Laws. In addition, the magnetic field needs to be divergence-free, and the electric field needs to satisfy Gauss's law.

Numerical schemes for the two-fluid model have been developed in several articles~\cite{Shumlak2003, Hakim2006, loverich2011, kumar2012entropy, amano2015divergence,  bond_wheatley_samtaney_pullin_2017, Abgrall2014, wang2020}. In~\cite{Shumlak2003}, an approximate Riemann solver is proposed for the model, while in~\cite{Hakim2006}, a wave-propagation-based scheme from~\cite{LeVeque2002} is designed for the two-fluid model. In~\cite{loverich2011}, Discontinuous Galerkin (DG) based schemes are proposed, while entropy stable schemes for the model were proposed in~\cite{kumar2012entropy}. In~\cite{Abgrall2014}, positivity-preserving schemes were proposed, which also discuss the stability of implicit source term updates. More recently, \cite{bond2016plasma} proposed an adaptive mesh refinement-based scheme. In~\cite{wang2020}, an exact solution for the source term was used to discretize the equations.  

To ensure divergence-free magnetic field evolution and evolution of the electric field consistent with Gauss's Law, several approaches are considered in the literature~\cite{Shumlak2003, kumar2012entropy, Hakim2006, amano2015divergence, balsara_high-order_2016}. In~\cite{Shumlak2003}, only one-dimensional problems were considered; hence, a constant magnetic field in the $x$-direction was present. To ensure consistency with Gauss's law, an additional elliptic equation was solved to update the electric field. In~\cite{Hakim2006, loverich2011, kumar2012entropy, Abgrall2014, bond_wheatley_samtaney_pullin_2017}, PHM equation based approach is considered. In~\cite{amano2015divergence}, a quasi-neutral two-fluid plasma model is considered; electric and magnetic fields are evolved on the faces of the cells using a multidimensional Riemann solver at the cell vertices and then coupled with the fluid variable at the cell centers. A similar approach was also considered in~\cite{amano_second-order_2016} for the two-fluid relativistic plasma flow equations. A higher-order scheme based on the multidimensional Riemann solvers approach is presented in~\cite{balsara_high-order_2016}. The works in~ \cite{amano2015divergence, amano_second-order_2016, balsara_high-order_2016} need the evolution of electric and magnetic fields on the cell faces.  {Similarly, for Maxwell's equations, mimetic schemes have been developed in~\cite{yee1966,HYMAN1999881,hyman2001} that mimic certain identities in vector calculus at he discrete level which then helps to satisfy divergence constraints. One key idea to achieve this is to store the electromagnetic variables on the faces and edges of the mesh, i.e., a staggered storage scheme is required.}

In this article, we design second-order accurate, co-located schemes {for ideal two-fluid plasma flow equations,} that ensure divergence-free evolution of the magnetic field and a consistent treatment of the electric field, which respects Gauss's law. We present the schemes in two dimensions and proceed as follows.
\begin{itemize}
\item We first consider the equations for magnetic and electric field evolution. We then use multidimensional Riemann solver~\cite{balsara_multidimensional_2014, balsara2014multidimensional, balsara_high-order_2016, chandrashekar2020} at each corner of the cells. Using these values, we define the numerical fluxes at each cell face. The resulting scheme is then shown to ensure the divergence-free evolution of the magnetic field and an evolution of the electric field that is consistent with a discrete version of Gauss's law. {Such ideas are sometimes called mimetic schemes, and we achieve this property with co-located storage of all unknowns.}
\item For the fluid parts, we consider the entropy stable numerical scheme~\cite{kumar2012entropy}. We note that this choice is not unique. Instead, W e can use any other stable discretization of the fluid flux, e.g., \cite{Abgrall2014}.
\item To achieve second-order accuracy, we reconstruct electric and magnetic fields at the cell corners using MinMod-based reconstruction. 
\item Finally, we consider explicit and Implicit-Explicit (IMEX) time update schemes. In the IMEX scheme, we treat the source term implicitly~\cite{kumar2012entropy, Abgrall2014}. We prove that the completely discrete schemes are also divergence-free and are consistent with the discretization of Gauss's law.
\end{itemize}

The rest of the article is organized as follows. In the next Section, we present the two-fluid plasma flow equations. In Section~\ref{sec:semi_disc}, we present the semi-discrete scheme based on the multidimensional Riemann solver for Maxwell's equations. We also describe the second-order extension and entropy-stable discretization of the fluid components. In Section~\ref{sec:fully_dis}, we present explicit and IMEX time discretizations. We also present the fully discrete evolutions of the divergence constraints. Numerical test cases are presented in Section~\ref{sec:num_results}, followed by concluding remarks in Section~\ref{sec:con}.

\section {Two-fluid plasma flow equations}
\label{sec:tf_plasma_eqns}
Following~\cite{Shumlak2003, Hakim2006, loverich2011, kumar2012entropy, bond_wheatley_samtaney_pullin_2017, amano2015divergence, Abgrall2014, polak_fourth-order_2023}, the fluid part of the two-fluid plasma flow equations can be written as:

	\begin{subequations}\label{eq:tf_fluid}
		\begin{align}
			%--- FLUID PART ---
			% Continuity (ion)
			\dfrac{\partial\rho_\al}{\partial t} + \nabla \cdot (\rho_\al \mathbf{u}_\al) &= 0 \label{eq:mass}\\
			% Momentum (ion)
			\dfrac{\partial (\rho_\al \mathbf{u}_\al)}{\partial t} + \nabla \cdot \left(\rho_\al \mathbf{u}_\al \mathbf{u}_\al^\top + p_\al \mathbf{I} \right) &= r_\al \rho_\al (\mathbf{E}+\mathbf{u}_\al \times \mathbf{B})\label{eq:momentum}\\
			% Energy (ion)
			\dfrac{\partial \mathcal{E}_\al}{\partial t} + \nabla \cdot 
			\big( (\mathcal{E}_\al + p_\al) \mathbf{u}_\al \big)  &= r_\al  \rho_\al (\mathbf{u}_\al \cdot \mathbf{E})\label{eq:energy}
%			%--- MAXWELL PART ---
%			% Faraday's Law
%			\dfrac{\partial \mathbf{B}}{\partial t} + \nabla \times \mathbf{E} &= 0, \label{induc}\\
%			% Generalised Ampere's Law
%			\dfrac{\partial \mathbf{E}}{\partial t} - c^2\nabla \times \mathbf{B} &= -\dfrac{\mathbf{j}}{\epsilon_0},\\
%			% Gauss's law
%			\nabla \cdot \mathbf{E} &= \dfrac{\rho_c }{\epsilon_0},\\
%			% Divergence-free magnetic field
%			\nabla \cdot \mathbf{B} &= 0.
		\end{align}
	\end{subequations}
Here, the subscript $\al\in\{I,E\}$ represents the ion or electron fluid. Within this framework, $\rho_\al$ represents the density, $\mathbf{u}_\al$ denotes the velocity vector, $p_\al$ is the pressure, and $\mathcal{E}_\al$ is the energy of the fluid. The electric and magnetic fields are denoted by $\Eb$ and $\Bb$, respectively. The charge-to-mass ratios are defined as $r_\alpha=\frac{q_\alpha}{m_\alpha}$, where $q_\alpha$ represents the charge and $m_\alpha$ represents the mass of the fluid for the species $\al\in\{I,E\}$. The Eqn.~\eqref{eq:mass} represents the conservation of mass, Eqn.~\eqref{eq:momentum} represents the conservation of momentum, and Eqn.~\eqref{eq:energy} represents the conservation of energy for the fluid $\al\in\{I,E\}$. The source terms in the system~\eqref{eq:tf_fluid} represent the Lorentz force acting on the fluid due to the presence of electromagnetic fields, $\Eb$ and $\Bb$. The fluid Eqns.~\eqref{eq:tf_fluid} are closed by considering the ideal equation of state,
\begin{equation}
	\label{eq:eqn_of_state}
	\mathcal{E}_\al=\frac{p_\al}{\gamma_\al-1}+\frac{1}{2}\rho_\al|\ub_\al|^2
\end{equation}
where $\gamma_\al$ is the ratio of specific heats. The evolution of the electromagnetic fields $\Eb$ and $\Bb$ is governed by Maxwell's equations,
	\begin{subequations}\label{eq:maxwell_eqn}
	\begin{align}
		%			%--- MAXWELL PART ---
		%			% Faraday's Law
				\dfrac{\partial \mathbf{B}}{\partial t} + \nabla \times \mathbf{E} &= 0 \label{eq:max_B}\\
		%			% Generalised Ampere's Law
					\dfrac{\partial \mathbf{E}}{\partial t} - c^2\nabla \times \mathbf{B} &= -\dfrac{\mathbf{j}}{\epsilon_0}\label{eq:max_E}\\
			%			% Divergence-free magnetic field
	\nabla \cdot \mathbf{B} &= 0,\label{eq:div_B}\\
		%			% Gauss's law
					\nabla \cdot \mathbf{E} &= \dfrac{\rho_c }{\epsilon_0} \label{eq:gauss_E}
	\end{align}
\end{subequations}
Here, $c=1/\sqrt{\mu_0\epsilon_0}$ is the speed of light, $\epsilon_0$ is the permittivity, $\mu_0$ is the permeability of free space, $\rho_c$ is total charge density and $\jb$ is the total current density given by,
\begin{equation}
	\label{eq:charge_current_density}
	\rho_c = r_I\rho_I + r_E\rho_E, \qquad \jb=r_I\rho_I \ub_I+ r_E\rho_E\ub_E
\end{equation}
The Eqn.~\eqref{eq:max_B} is Faraday's law for the magnetic field, and Eqn.~\eqref{eq:max_E} represents Ampere's law for the electric field. The magnetic field $\Bb$ has to satisfy the divergence-free condition given in Eqn.~\eqref{eq:div_B}, and the electric field $\Eb$ has to satisfy the divergence constraint Eqn.~\eqref{eq:gauss_E} given by Gauss's law.

The vector of conservative variables that consists of fluid and electromagnetic parts is given by
$$
\Ub=(\Ub_I^\top,\Ub_E^\top,\Ub_M^\top)^\top
$$
where $\Ub_\al=(\rho_\al,\rho_\al\ub_\al,\mathcal{E}_\al)^\top$, $\al\in\{I,E\}$ are the fluid variables and $\Ub_M=(\Bb,\Eb)^\top$ is the vector of electromagnetic variables. The vector quantities $\ub$, $\jb$, $\Eb$ and $\Bb$ are given by $\ub_\al=(u_\al^x,u_\al^y,u_\al^z)$, $\jb=(j_x,j_y,j_z),$ $\Eb=(E_x,E_y,E_z),$ and $\Bb=(B_x,B_y,B_z)$. Then the system~\eqref{eq:tf_fluid} along with the equations~\eqref{eq:max_B} and~\eqref{eq:max_E} can be written as,
\begin{equation}
	\label{eq:tf_cons_form}
	\frac{\pa \Ub}{\pa t} + \frac{\pa \fb^x}{\pa x} + \frac{\pa \fb^y}{\pa y} =\sbb(\Ub)
\end{equation}
where, the fluxes $\fb^x$ and $\fb^y$ can be divided into three independent components given by,
$$
\fb^x=\begin{pmatrix}
	\fb_I^x\\
	\fb_E^x\\
	\fb_M^x
\end{pmatrix},
\qquad 
\fb^y=\begin{pmatrix}
	\fb_I^y\\
	\fb_E^y\\
	\fb_m^y
\end{pmatrix}
$$ 
where
\begin{equation}	
\fb_\al^x=\begin{pmatrix}
\rho_\al u^x_\al\\
\rho_\al (u_\al^x)^2 +p_\al\\
\rho_\al u_\al^x u_\al^y\\
\rho_\al u_\al^x u_\al^z\\
(\mathcal{E}_\al+p_\al)u^x_\al	
\end{pmatrix},
\quad 
	\fb_\al^y=\begin{pmatrix}
		\rho_\al u^y_\al\\
		\rho_\al u_\al^x u_\al^y\\
		\rho_\al (u_\al^y)^2 +p_\al\\
			\rho_\al u_\al^y u_\al^z\\
		(\mathcal{E}_\al+p_\al)u^y_\al	
	\end{pmatrix},
\quad 
\fb_M^x=\begin{pmatrix}
	0\\
	-E_z\\
	E_y\\
	0\\
	c^2B_z\\
	-c^2B_y
\end{pmatrix},
\quad 
\fb_M^y=\begin{pmatrix}
	E_z\\
	0\\
	-E_x\\
	-c^2B_z\\
	0\\
	c^2B_x
\end{pmatrix}	
\end{equation}
These three components are then coupled via source term,
$$
\sbb(\Ub)=\begin{pmatrix}
	\sbb_I(\Ub_I,\Ub_M)\\
	\sbb_E(\Ub_E,\Ub_M)\\
	\sbb_M(\Ub_I,\Ub_E)
\end{pmatrix}
$$
where,
\begin{equation}
	\sbb_\al(\Ub_\al,\Ub_M)=\begin{pmatrix}
		0\\
		r_\al\rho_\al(E_x+u_\al^yB_z-u_\al^zB_y)\\
		r_\al\rho_\al(E_y+u_\al^zB_x-u_\al^xB_z)\\
		r_\al\rho_\al(E_z+u_\al^xB_y-u_\al^yB_x)\\
		r_\al\rho_\al(u_\al^xE_x +u_\al^yE_y+u_\al^zE_z)
	\end{pmatrix},
	\quad 
	\sbb_M(\Ub_I,\Ub_E)=-\frac{1}{\epsilon_0}
	\begin{pmatrix}
		0\\
	    0\\
	    0\\
		j_x\\
		j_y\\
		j_z\\
	\end{pmatrix}
\end{equation}
We also introduce the vector of primitive variables $\Wb=(\Wb_I^\top,\Wb_E^\top,\Wb_M^\top)^\top,$ where $\Wb_\al=(\rho_\al,\ub_\al,p_\al)^\top$, $\al\in\{I,E\}$ and  $\Wb_M=(\Bb,\Eb)^\top$. Following \cite{Abgrall2014}, the set of admissible solutions is given by,
$$
\Omega=\left\{\Ub\in\mathbb{R}^{16} \ | \ \rho_\al >0, p_\al >0,\al\in\{I,E\}\right\}
$$
The system \eqref{eq:tf_cons_form} is hyperbolic for the solutions in $\Omega$. Furthermore, the eigenvalues of the system in $x$ and $y$-direction are given by,
\begin{equation}
	\mathbf{\Lambda}^d=\left\{u_I^d-a_I,u_I^d,u_I^d,u_I^d,u_I^d+a_I,u_E^d-a_E,u_E^d,u_E^d,u_E^d,u_E^d+a_E,-c,-c,0,0,c,c \right\}, ~\text{for}~ d\in\{x,y\}
\end{equation}
where $a_\al=\sqrt{\gamma_\al p_\al/\rho_\al}$, for $\al\in\{I,E\}$. {We note that as the three fluxes $\fb_I^x,\fb_E^x$ and $\fb_M^x$ are independent of each other, and we have independent eigenvalues and eigenvectors for each of them. Similarly, fluxes $\fb_I^y,\fb_E^y$ and $\fb_M^y$ are also independent of each other.}

Following \cite{kumar2012entropy}, we also define the entropy $e_\al$ and the entropy flux $\qb_\al$ as,
\begin{equation}
	e_\al=\frac{-\rho_\al s_\al}{\gamma_\al-1}, \qquad \qb_\al = e_\al \ub_\al, \qquad s_\al=\log{(p_\al)} - \gamma_\al \log{(\rho_\al)}
\end{equation}
We now recall the entropy inequality from \cite{kumar2012entropy}:
\begin{proposition}[see \cite{kumar2012entropy}]
	Smooth solutions of \eqref{eq:tf_cons_form} satisfy
	$$
	\pa_t s_\al + u_\al^x \pa_x s_\al + u_\al^y \pa_y s_\al = 0
	$$
	which  results in,
	\begin{equation}
		\label{eq:ent_eqality}
		\pa_t e_\al + \pa _x q^x_\al + \pa _y q^y_\al =0 
	\end{equation}
\end{proposition}
Furthermore, for non-smooth weak solutions, we should satisfy the entropy inequality,
\begin{equation}
	\label{eq:ent_inq}
		\pa_t e_\al + \pa _x q^x_\al + \pa _y q^y_\al \le 0 
\end{equation}
 
\subsection{Divergence constraints of Maxwell's equations}
In addition to the system of equations in~\eqref{eq:tf_cons_form}, the electric and magnetic fields also need to satisfy the additional divergence constraints in the form of Gauss's law~ \eqref{eq:gauss_E} and the divergence-free condition~\eqref{eq:div_B}. Let us consider the Faraday's law~\eqref{eq:max_B}; taking the divergence of this equation results in
\begin{equation}
	\label{eq:divB_evo_cont}
	\frac{\pa }{\pa t}(\na \cdot \Bb) =0
\end{equation}
Hence, if $\na\cdot\Bb=0$ at the initial time, Eqn.~\eqref{eq:max_B} will ensure~\eqref{eq:div_B} for all time. Therefore, satisfying the divergence constraint~\eqref{eq:div_B} for all time is an outcome of the initial divergence-free magnetic field and Eqn.~\eqref{eq:max_B}. Similarly, taking the divergence of  Ampere's law~\eqref{eq:max_E}, we get,
\begin{equation}
	\label{eq:divE_evo_cont}
	\frac{\pa }{\pa t}(\na \cdot \Eb) = -\frac{1}{\eps_0} \na\cdot \jb = -\frac{1}{\eps_0} (r_I \nabla\cdot(\rho_I \ub_I) + r_E \nabla\cdot(\rho_E\ub_E)) = \frac{1}{\eps_0} \frac{\pa \rho_c}{\pa t} 
\end{equation}
This is the time derivative of the Eqn.~\eqref{eq:gauss_E}; hence,~\eqref{eq:gauss_E} is a consequence of the fact that it is satisfied at the initial time and $\Eb$ is evolved using~\eqref{eq:max_E}. Hence, the system~\eqref{eq:maxwell_eqn} is not over-determined and the constraints should be considered as the conditions on the initial data.

 While the previous two constraint equations are respected at the continuous level, at the discrete level this is not guaranteed ~\cite{jiang_origin_1996} unless a good mimetic discretization is used. Neglecting these constraints can lead to spurious non-physical oscillations in the solution. This is particularly significant when dealing with oscillatory source terms, as it can lead to non-physical oscillations in the overall solution.

To overcome this, several recent articles ~\cite{Hakim2006,loverich2011,kumar2012entropy,Abgrall2014,bond_wheatley_samtaney_pullin_2017} use {\em perfectly hyperbolic Maxwell's equations} (PHM) from~\cite{munz2000}. The Eqns.~\eqref{eq:maxwell_eqn} are modified by introducing "correction potentials", $\phi$ and $\psi$ to get,
 \begin{subequations}\label{eq:phm_eqn}
 	\begin{align}
 		%			%--- MAXWELL PART ---
 		%			% Faraday's Law
 		\dfrac{\partial \mathbf{B}}{\partial t} + \nabla \times \mathbf{E} +\kappa \nabla \psi &= 0 \label{eq:phm_B}\\
 		%			% Generalised Ampere's Law
 		\dfrac{\partial \mathbf{E}}{\partial t} - c^2\nabla \times \mathbf{B}  + \xi c^2 \nabla \phi&= -\dfrac{\mathbf{j}}{\epsilon_0}\label{eq:phm_E}\\
 		%			% Divergence-free magnetic field
 		\frac{\partial \psi}{\pa t}+ \kappa c^2\nabla \cdot \mathbf{B} &= 0\label{eq:phm_div_B}\\
 		%			% Gauss's law
 		\frac{\partial \phi}{\pa t} +\xi  \nabla \cdot \mathbf{E} &= \xi \dfrac{\rho_c }{\epsilon_0} \label{eq:phm_gauss_E}
 	\end{align}
 \end{subequations}
 Here, $\kappa$ and $\xi$ are error propagation speeds. In the limiting case of $\kappa,\xi \to\infty$, both divergence constraints \eqref{eq:div_B} and \eqref{eq:gauss_E} are satisfied. However, in practice, $\kappa,\xi=1$ or $2$ is considered. One key benefit of the PHM formulation is that the resulting equations are still hyperbolic. However, the resulting system has larger eigenvalues, consequently smaller timesteps, leading to an increase in wallclock time. Furthermore, the divergence constraints are satisfied only in the limit of large $\kappa,\xi$.
 
In this section, we have analyzed the continuous problem and corresponding entropy condition. We have also discussed the divergence constraints of Maxwell's equations. In the next Section, we will design numerical schemes for the system~\eqref{eq:tf_cons_form} consistent with the entropy condition~\eqref{eq:ent_eqality} and the divergence constraints~\eqref{eq:divB_evo_cont} and~\eqref{eq:divE_evo_cont}. 
 
 \section{Semi-discrete schemes}
\label{sec:semi_disc}

Let us consider a two-dimensional rectangular domain $D=(x_{\min},x_{max})\times(y_{\min},y_{max}).$ We discretize the domain $D$ using a uniform mesh with the cell size $\Dx\times\Dy$. We now define $x_i=x_{min}+(i+1/2)\Dx$ and $y_j=y_{min}+(j+1/2)\Dy$, with $0\le i< N_x$ and $0 \le j < N_y$ such that, $x_{max}=x_{min}+N_x\Dx$ and $y_{max}=y_{min}+N_y \Dy$. Let us also define $x_\iph=\frac{x_i+x_{i+1}}{2}$ and $y_\jph=\frac{y_j+y_{j+1}}{2}$. A semi-discrete scheme for the system~\eqref{eq:tf_cons_form} is given by,
\begin{equation}
	\label{eq:semi_discrete}
	\frac{d \Ub_{i,j}}{d t} +\frac{1}{\Dx} \left(\Fb^x_{\iph,j} - \Fb^x_{\imh,j}\right) + \frac{1}{\Dy} \left(\Fb^y_{i,\jph} - \Fb^y_{i,\jmh}\right)  = \sbb(\Ub_{i,j})
\end{equation}  
Here, $\Fb^x_{\iph,j}$ and $\Fb^y_{i,\jph}$ are the numerical fluxes consistent with the continuous fluxes $\fb^x$ and $\fb^y$, respectively. As the fluxes  $\fb^x$ and $\fb^y$ have three independent components, the numerical fluxes $\Fb^x_{\iph,j}$ and $\Fb^y_{i,\jph}$ will be of the form,
\begin{equation}
	\Fb^x_{\iph,j}=\begin{pmatrix}
		\Fb_{I,\iph,j}^x\\
		\Fb_{E,\iph,j}^x\\
		\Fb_{M,\iph,j}^x
	\end{pmatrix},
	\qquad 
	\Fb^y_{i,\jph}=\begin{pmatrix}
		\Fb_{I,i,\jph}^y\\
		\Fb_{E,i,\jph}^y\\
		\Fb_{M,i,\jph}^y
	\end{pmatrix}
\end{equation}
Furthermore, each component of the numerical flux is consistent with the corresponding continuous flux. We will first describe $\Fb_{M,\iph,j}^x$ and $\Fb_{M,i,\jph}^y$, which are numerical fluxes for the Maxwell's equations~\eqref{eq:max_B} and~\eqref{eq:max_E}, consistent with the continuous fluxes $\fb_M^x$ and $\fb^y_M$, respectively.

% \begin{defn}{}
% The semi-discrete scheme~\eqref{eq:semi_discrete} is said to be entropy stable if the following discrete entropy inequality is satisfied:
% 	\begin{equation}
% 		\frac{d}{dt}  e_\al(\mathbf{U}_{\alpha,i,j})  +\frac{1}{\Delta x} \left( \hat{q}_{\alpha,i+\frac{1}{2},j}^x - \hat{q}_{\alpha,i-\frac{1}{2},j}^x\right)+\frac{1}{\Delta y}\left( \hat{q}_{\alpha,i,j+\frac{1}{2}}^y - \hat{q}_{\alpha,i,j-\frac{1}{2}}^y\right) \le 0 \ \ \ \alpha \in \{i,e\}, \label{eq:semi_dis_fluid_entropy_inequality}
% 	\end{equation}
% where $\hat{q}^{x}_{\alpha,i+\frac{1}{2},j}$, $\hat{q}^{y}_{\alpha, i,j+\frac{1}{2}}$ are the numerical entropy flux functions consistent with the continuous entropy fluxes $q^x_\al, q^y_\al$ respectively.
% \end{defn}

\subsection{Multidimensional Local Lax-Friedrich flux for Maxwell's equations}
\label{subsec:multid_riem_solver}
At point $\left(x_\iph,y_\jph\right)$, four states $\Ub_{M,i,j},\Ub_{M,i+1,j},\Ub_{M,i,j+1}$ and $\Ub_{M,i+1,j+1}$ interact. Let us define the following:
\begin{equation}
\label{eq:max_face_averages}
\Ub_{M,\iph,j}=\frac{\Ub_{M,i+1,j}+ \Ub_{M,i,j}}{2},\quad	\Ub_{M,i,\jph} = \frac{\Ub_{M,i,j+1} + \Ub_{M,i,j}}{2}
\end{equation}
\begin{equation}
	\label{eq:max_face_averages2}
	\bar{\Ub}_{M,\iph,\jph}=\frac{\Ub_{M,i,j}+ \Ub_{M,i+1,j} +\Ub_{M,i,j+1} + \Ub_{M,i+1,j+1}}{4}
\end{equation}
Following \cite{balsara_multidimensional_2014,balsara2014multidimensional,balsara_high-order_2016,chandrashekar2020}, we now define,
\begin{equation}
     \tilde{E}_{z,\iph,\jph} = \bar{E}_{z,\iph,\jph} +\frac{c}{2}\left( \left( B_{y,i+1,\jph}-B_{y,i,\jph} \right)  - \left( B_{x,\iph,j+1}-B_{x,\iph,j} \right)\right), 	
\label{eq:multid_Ez}
\end{equation}
and 
%
% \begin{equation}
%   c^2\tilde{B}_{z,\iph,\jph}= c^2\bar{B}_{z,\iph,\jph}-\frac{c}{2}\left( (B_{y,i+1,\jph}-B_{y,i,\jph})  - (B_{x,\iph,j+1}-B_{x,\iph,j})\right). 
%     \label{eq:multid_Bz} 
% \end{equation}
%
\begin{equation}
    c^2\tilde{B}_{z,\iph,\jph}
    =
    c^2\bar{B}_{z,\iph,\jph}
    -
    \frac{c}{2}
    \left(
        \left( E_{x,\iph,j+1}-E_{x,\iph,j} \right)
    -   \left( E_{y,i+1,\jph}-E_{y,i,\jph} \right)
    \right).
    \label{eq:multid_Bz} 
\end{equation}
Here, we note that the multidimensional Riemann solver has two parts, namely, the central part denoted the averages $\bar{E}_{z,\iph,\jph}$ and $\bar{B}_{z,\iph,\jph}$ and a multidimensional dissipation given by the remaining terms. {This approximation can also be motivated starting from an artificial viscosity approach, see Section~\ref{sec:Ez}. We use only light speed $c$ in the dissipative part of the flux since the different PDE models are decoupled in terms of the fluxes which determine the wave speeds.} Now, using the values of $\tilde{E}_{z,\iph,\jph}$, and $c^2\tilde{B}_{z,\iph,\jph}$, we can define the numerical flux $\Fb_{m,\iph,j}^x$ as,
\begin{equation}
\label{eq:multiD_flux_x}
\Fb_{M,\iph,j}^x=\begin{pmatrix}
	0\\
	-\dfrac{1}{2} \left(\tilde{E}_{z,\iph,\jph} +\tilde{E}_{z,\iph,\jmh} \right)  \\
	\left(\tilde{F}_{M,\iph,j}^x \right)_{B_z} \\
	0\\
	\dfrac{1}{2}\left(c^2\tilde{B}_{z,\iph,\jph} +c^2\tilde{B}_{z,\iph,\jmh}\right)\\
	\left(\tilde{F}_{M,\iph,j}^x \right)_{E_z}
\end{pmatrix},
\end{equation}
where, $\left(\tilde{F}_{M,\iph,j}^x \right)_{B_z}$ and $\left(\tilde{F}_{M,\iph,j}^x \right)_{E_z}$ are obtained using one-dimensional Rusanov's solver as follows:
% ------
\begin{align}
   \left(\tilde{F}_{M,\iph,j}^x \right)_{B_z} & =  \dfrac{{E}_{y,i,j} + {E}_{y,i+1,j}}{2}
   - 
   \dfrac{c}{2} \left( {B}_{z,i+1,j} - {B}_{z,i,j} \right),
   \label{eq:oned_rus_x_Bz}
   % ------ 
   \\
   % ------
   \left( \tilde{F}_{M,\iph,j}^x \right)_{E_z}  &=
   -c^2\left(
   \dfrac{{B}_{y,i,j} + {B}_{y,i+1,j}}{2}
   \right)
   - 
   \dfrac{c}{2} \left( {E}_{z,i+1,j} - {E}_{z,i,j} 
   \right).
   \label{eq:oned_rus_x_Ez}
\end{align}
% ------
% -------------------------------
% ------
Similarly, the $y$-directional numerical flux $\Fb_{m,i,\jph}^y$ is given by, 
\begin{equation}
	\label{eq:multiD_flux_y}
\Fb_{M,i,\jph}^y=\begin{pmatrix}
	 \dfrac{1}{2} \left(\tilde{E}_{z,\iph,\jph} +\tilde{E}_{z,\imh,\jph} \right)\\
	0
    \\
	\left(\tilde{F}_{M,i,\jph}^y \right)_{B_z}\\
	-\dfrac{1}{2} \left(c^2\tilde{B}_{z,\iph,\jph} +c^2\tilde{B}_{z,\imh,\jph} \right)\\
	0\\
	\left(\tilde{F}_{M,i,\jph}^y \right)_{E_z}
\end{pmatrix},	
\end{equation}
where, $\left(\tilde{F}_{M,i,\jph}^y \right)_{B_z}$ and $\left(\tilde{F}_{M,i,\jph}^y \right)_{E_z}$ are obtained using one-dimensional Rusanov's solver as follows:
% ------
\begin{align}
   \left(\tilde{F}_{M,i,\jph}^y \right)_{B_z} &= 
   - \left(
   \dfrac{{E}_{x,i,j} + {E}_{x,i,j+1}}{2}
   \right)
   - 
   \dfrac{c}{2} \left( {B}_{z,i,j+1} - {B}_{z,i,j} \right),
   \label{eq:oned_rus_y_Bz}
   \\
% ------
   \left(\tilde{F}_{M,i,\jph}^y \right)_{E_z} &= 
   c^2 \left(
   \dfrac{{B}_{x,i,j} + {B}_{x,i,j+1}}{2}
   \right)
   - 
   \dfrac{c}{2} \left( {E}_{z,i,j+1} - {E}_{z,i,j} \right).
   \label{eq:oned_rus_y_Ez}
\end{align}
We note that the numerical fluxes $\Fb_{M,\iph,j}^x$ and $\Fb_{M,i,\jph}^y$ are consistent with the one dimensional fluxes~\cite{chandrashekar2020, balsara_high-order_2016}. Furthermore, as eigenvalues of Maxwell's fluxes $\fb_M^x$ and $\fb_M^y$ are constants, $c$ and $-c$ {(as the fluid variable do not interact with Maxwell's flux)}, the numerical fluxes \eqref{eq:multiD_flux_x} and \eqref{eq:multiD_flux_y} are equivalent to the multidimensional HLL numerical flux.

\subsection{Second-order reconstruction for Maxwell's equations}
\label{subsec:max_2nd_order}
\begin{figure}[ht]
\begin{center}
	\includegraphics[width=3in, height=3in]{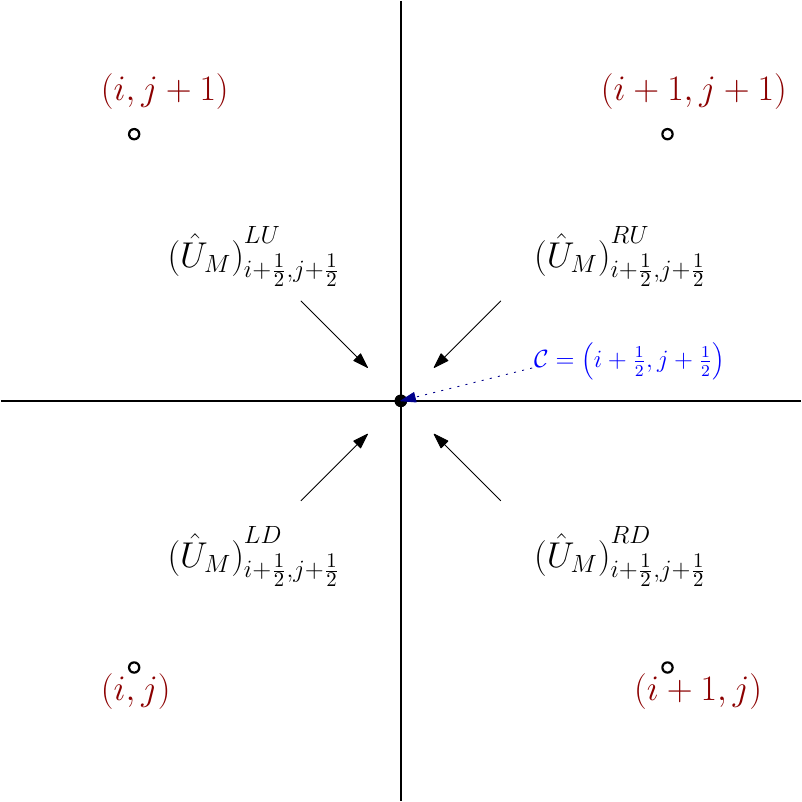}
	\caption{Part of the two-dimensional grid showing the reconstructed values from the four neighboring zones at the vertex point $\mathcal{C}=\left( i+\frac{1}{2}, j+\frac{1}{2} \right)$.}
	\label{fig:grid}
\end{center}
\end{figure}

To achieve second-order accuracy for Maxwell's part, we use the {\em MinMod} slope limiter to obtain reconstructed values at the cell interfaces and cell corners. The {\em MinMod} function of two real numbers $a$ and $b$ is defined as
\begin{align}
   \text{MinMod} (a,b) = \begin{cases}
       \text{sign}(a),  \min \{|a|,|b|\} & \text{if} \quad \text{sign}(a) = \text{sign}(b),
       \\
       0, & \text{otherwise.}
   \end{cases}
   \label{eq:minmod}
\end{align}
Using the slope limiter~\eqref{eq:minmod}, we obtain the traces of $B_z$ and $E_z$ along all the four diagonal by defining (see Figure \eqref{fig:grid}),
\begin{align*}
    (\hat{\Ub}_M)^{LD}_{i+\frac{1}{2},j+\frac{1}{2}} 
    &=  
    \Ub_{M,i,j} + \frac{1}{2} \text{MinMod}
    \Big\{
    \Ub_{M,i,j} - \Ub_{M,i-1,j-1}, \Ub_{M,i+1,j+1} - \Ub_{M,i,j} 
    \Big\},
    % ------
    \\
    % ------
    (\hat{\Ub}_M)^{RD}_{i+\frac{1}{2},j+\frac{1}{2}} 
    &=  
    \Ub_{M,i+1,j} - \frac{1}{2} \text{MinMod}
    \Big\{
    \Ub_{M,i+1,j} - \Ub_{M,i,j+1}, \Ub_{M,i+2,j-1} - \Ub_{M,i+1,j} 
    \Big\},
    % ------
    \\
    % ------
    (\hat{\Ub}_M)^{RU}_{i+\frac{1}{2},j+\frac{1}{2}} 
    &=  
    \Ub_{M,i+1,j+1} - \frac{1}{2} \text{MinMod}
    \Big\{
    \Ub_{M,i+1,j+1} - \Ub_{M,i,j}, \Ub_{M,i+2,j+2} - \Ub_{M,i+1,j+1} 
    \Big\},
    % ------
    \\
    % ------
    (\hat{\Ub}_M)^{LU}_{i+\frac{1}{2},j+\frac{1}{2}} 
    &=  
    \Ub_{M,i,j+1} + \frac{1}{2} \text{MinMod}
    \Big\{
    \Ub_{M,i,j+1} - \Ub_{M,i-1,j+2}, \Ub_{M,i+1,j} - \Ub_{M,i,j+1} 
    \Big\}.
\end{align*}

Now, using these diagonal traces with \eqref{eq:multid_Bz} and \eqref{eq:multid_Ez}, we define the vertex values $\tilde{B}_{z,\iph,\jph}$ and $\tilde{E}_{z,\iph,\jph}$) at the cell vertex $\left(x_\iph,y_\jph\right)$ as,
% ------
\begin{equation}
\begin{aligned}
    \tilde{E}_{z,\iph,\jph}
      & = \dfrac{
        (\hat{E}_z)^{LD}_{i+\frac{1}{2},j+\frac{1}{2}}
      + (\hat{E}_z)^{RD}_{i+\frac{1}{2},j+\frac{1}{2}}
      + (\hat{E}_z)^{RU}_{i+\frac{1}{2},j+\frac{1}{2}}
      + (\hat{E}_z)^{LU}_{i+\frac{1}{2},j+\frac{1}{2}}
        }{4}
      \\ 
      &+ 
      \frac{c}{2}
         \left(
            \dfrac{(\hat{B}_y)^{RD}_{i+\frac{1}{2},j+\frac{1}{2}}
                +  (\hat{B}_y)^{RU}_{i+\frac{1}{2},j+\frac{1}{2}}}{2}
                -
            \dfrac{(\hat{B}_y)^{LD}_{i+\frac{1}{2},j+\frac{1}{2}}
                +  (\hat{B}_y)^{LU}_{i+\frac{1}{2},j+\frac{1}{2}}}{2}  
         \right)
      \\ 
      &- 
      \frac{c}{2}
         \left(
            \dfrac{(\hat{B}_x)^{LU}_{i+\frac{1}{2},j+\frac{1}{2}}
                +  (\hat{B}_x)^{RU}_{i+\frac{1}{2},j+\frac{1}{2}}}{2}
                -
            \dfrac{(\hat{B}_x)^{LD}_{i+\frac{1}{2},j+\frac{1}{2}}
                +  (\hat{B}_x)^{RD}_{i+\frac{1}{2},j+\frac{1}{2}}}{2}  
         \right).	
\label{eq:multid_Ez_2nd}
\end{aligned}
\end{equation}
% ------
and
% ------
\begin{equation}
\begin{aligned}
    c^2 \tilde{B}_{z,\iph,\jph}
      & = \dfrac{c^2}{4} \left( 
        (\hat{B}_z)^{LD}_{i+\frac{1}{2},j+\frac{1}{2}}
      + (\hat{B}_z)^{RD}_{i+\frac{1}{2},j+\frac{1}{2}}
      + (\hat{B}_z)^{RU}_{i+\frac{1}{2},j+\frac{1}{2}}
      + (\hat{B}_z)^{LU}_{i+\frac{1}{2},j+\frac{1}{2}} \right)
      \\ 
      & -
      \frac{c}{2}
         \left(
            \dfrac{(\hat{E}_x)^{LU}_{i+\frac{1}{2},j+\frac{1}{2}}
                +  (\hat{E}_x)^{RU}_{i+\frac{1}{2},j+\frac{1}{2}}}{2}
                -
            \dfrac{(\hat{E}_x)^{LD}_{i+\frac{1}{2},j+\frac{1}{2}}
                +  (\hat{E}_x)^{RD}_{i+\frac{1}{2},j+\frac{1}{2}}}{2}  
         \right)
      \\ 
      & + 
      \frac{c}{2}
         \left(
            \dfrac{(\hat{E}_y)^{RD}_{i+\frac{1}{2},j+\frac{1}{2}}
                +  (\hat{E}_y)^{RU}_{i+\frac{1}{2},j+\frac{1}{2}}}{2}
                -
            \dfrac{(\hat{E}_y)^{LD}_{i+\frac{1}{2},j+\frac{1}{2}}
                +  (\hat{E}_y)^{LU}_{i+\frac{1}{2},j+\frac{1}{2}}}{2}  
         \right).	
\label{eq:multid_Bz_2nd}
\end{aligned}
\end{equation}
To compute second-order accurate $\left(\tilde{F}_{M,\iph,j}^x \right)_{B_z}$ and $\left(\tilde{F}_{M,\iph,j}^x \right)_{E_z}$, we use {\em MinMod} limiter in $x$-direction and one-dimensional Rusanov's solver. We first compute the traces in $x$- direction,
\begin{align}
   {\check{\Ub}}^{-}_{M,\iph,j}
   &=
   \Ub_{M,i,j} + \dfrac{1}{2} 
   \text{MinMod} \Big\{
    \Ub_{M,i,j} - \Ub_{M,i-1,j}, \Ub_{M,i+1,j} - \Ub_{M,i,j} 
    \Big\},
    \\
    {\check{\Ub}}^{+}_{M,\iph,j}
    &=
    \Ub_{M,i+1,j} - \dfrac{1}{2} 
   \text{MinMod} \Big\{
    \Ub_{M,i+1,j} - \Ub_{M,i,j}, \Ub_{M,i+2,j} - \Ub_{M,i+1,j} 
    \Big\}.
\end{align}
and use them in one-dimensional Rusanov's solver \eqref{eq:oned_rus_x_Bz} and \eqref{eq:oned_rus_x_Ez}, to define, 

\begin{align}
   \left(\tilde{F}_{M,\iph,j}^x \right)_{B_z} &=  \dfrac{{\check{E}}^{-}_{y,\iph,j} + {\check{E}}^+_{y,\iph,j}}{2}
   - 
   \dfrac{c}{2} \left( {\check{B}}^+_{z,\iph,j} - {\check{B}}^-_{z,\iph,j} \right),
   \label{eq:oned_rus_x_Bz_2nd}
   \\
% ------
   \left(\tilde{F}_{M,\iph,j}^x \right)_{E_z} &=  \dfrac{ -c^2   }{2}
   \left( {\check{B}}^-_{y,\iph,j} + {\check{B}}^+_{y,\iph,j} \right)
   - 
   \dfrac{c}{2} \left( {\check{E}}^+_{z,\iph,j} - {\check{E}}^-_{z,\iph,j} \right).
   \label{eq:oned_rus_x_Ez_2nd}
\end{align}

Using \eqref{eq:multid_Bz_2nd},\eqref{eq:multid_Ez_2nd},\eqref{eq:oned_rus_x_Bz_2nd} and \eqref{eq:oned_rus_x_Ez_2nd} in \eqref{eq:multiD_flux_x}, we get the second-order accurate approximation of the flux in $x$-direction. We approximate second-order accurate $y$-directional flux, similarly.
% % ------
%
% ------
%

\subsection{Entropy stable discretization of fluid equations}
\label{subsec:ent_stable_fluid}
We now consider the fluid components $\Fb_\alpha$ of the numerical fluxes. To ensure the entropy-stable discretizations of the fluid components, we follow \cite{kumar2012entropy} and give a brief description of the numerical scheme. We introduce the following notations:
$$
[a]_{\iph,j} = a_{i+1,j}-a_{i,j}, \qquad [a]_{i,\jph} = a_{i,j+1}-a_{i,j}
$$
and
$$
\bar{a}_{\iph,j} = \frac{a_{i+1,j}+a_{i,j}}{2},\qquad \bar{a}_{i,\jph} = \frac{a_{i,j+1}+a_{i,j}}{2}
$$
A consistent entropy stable numerical flux  has the following form:
\begin{equation}
	\Fb^x_{\al,\iph,j} = \tilde{\Fb}^x_{\al,\iph,j}-\frac{1}{2}\Db^x_{\al,\iph,j}[\Vb_\alpha]_{\iph,j},  \qquad 
	\Fb^x_{\al,\iph,j} = \tilde{\Fb}^y_{\al,i,\jph}-\frac{1}{2}\Db^y_{\al,i,\jph}[\Vb_\alpha]_{i,\jph} 
\end{equation}
where $\Vb_\alpha=\dfrac{\pa e_\alpha}{\pa \Ub_\alpha}$ is the entropy variable, and $\chi^x_\al=\Vb_\al \cdot \Fb^x_{\al}- q^x_\al$, $\chi^y_\al=\Vb_\al \cdot \Fb^y_{\al}- q^y_\al$ are entropy potentials. Furthermore, $\Db^x_{\al,\iph,j}$ and $\Db^y_{\al,i,\jph}$ are some positive definite symmetric matrices which we will define later in this section. The $\tilde{\Fb}^x_{\al,\iph,j}$ and $\tilde{\Fb}^y_{\al,i,\jph}$  are consistent  entropy conservative numerical fluxes satisfying Tadmor conditions,
$$
[\Vb_\al]_{\iph,j}	\cdot \tilde{\Fb}^x_{\al,\iph,j} = [\chi^x_\al]_{\iph,j},\qquad 
[\Vb_\al]_{i,\jph}	\cdot \tilde{\Fb}^x_{\al,i,\jph} = [\chi^x_\al]_{i,\jph}.
$$
Following~\cite{Chandrashekar_2013} they are given by, 
\begin{align*}
	\label{eq:ent_cons_numflux}
	\tilde{\Fb}^x_{\al,\iph,j} = \begin{pmatrix}
		\hat{\rho}_{\al,\iph,j} \bar{u}^x_{\al,\iph,j}\\
		\frac{\bar{\rho}_{\al,\iph,j}}{2\bar{\beta}_{\al,\iph,j}} + 	\hat{\rho}_{\al,\iph,j} (\bar{u}^x_{\al,\iph,j})^2\\
		\hat{\rho}_{\al,\iph,j}\bar{u}^x_{\al,\iph,j}\bar{u}^y_{\al,\iph,j}\\
		\hat{\rho}_{\al,\iph,j} \bar{u}^x_{\al,\iph,j}\bar{u}^z_{\al,\iph,j}\\
		\frac{1}{2}\bigg(\frac{1}{(\gamma-1)\hat{\beta}_{\al,\iph,j}} - \bar{\mathbf{u}}^2_{\al,\iph,j}\bigg)\hat{\rho}_{\al,\iph,j} \bar{u}^x_{\al,\iph,j} + (\bar{u}^x_{\al,\iph,j})^3\hat{\rho}_{\al,\iph,j}+\bar{u}^x_{\al,\iph,j}\frac{\bar{\rho}_{\al,\iph,j}}{2\bar{\beta}_{\al,\iph,j}} \\
		+ (\bar{u}^y_{\al,\iph,j})^2\bar{u}^x_{\al,\iph,j}\hat{\rho}_{\al,\iph,j} + (\bar{u}^z_{\al,\iph,j})^2\bar{u}^x_{\al,\iph,j}\hat{\rho}_{\al,\iph,j}
	\end{pmatrix} 
\end{align*}
and
\begin{align*}
	\tilde{\Fb}^y_{\al,i,\jph} = \begin{pmatrix}
		\hat{\rho}_{\al,i,\jph} \bar{u}^y_{\al,i,\jph}\\
	\hat{\rho}_{\al,i,\jph}\bar{u}^x_{\al,i,\jph,}\bar{u}^y_{\al,i,\jph}\\
		\hat{\rho}_{\al,i,\jph} (\bar{u}^y_{\al,i,\jph})^2 +\frac{\bar{\rho}_{\al,i,\jph}}{2\bar{\beta}_{\al,i,\jph}}\\
		\hat{\rho}_{\al,i,\jph}\bar{u}^z_{\al,i,\jph}\bar{u}^y_{\al,i,\jph}\\
			\frac{1}{2}\bigg(\frac{1}{(\gamma-1)\hat{\beta}_{\al,i,\jph}} - \bar{\mathbf{u}}^2_{\al,i,\jph}\bigg)\hat{\rho}_{\al,i,\jph} \bar{u}^y_{\al,i,\jph} + (\bar{u}^x_{\al,i,\jph})^2\bar{u}^y_{\al,i,\jph}\hat{\rho}_{\al,i,\jph} \\
		+(\bar{u}^y_{\al,i,\jph})^3\hat{\rho}_{\al,i,\jph}+\bar{u}^y_{\al,i,\jph}\frac{\bar{\rho}_{\al,i,\jph}}{2\bar{\beta}_{\al,i,\jph}}	
 + (\bar{u}^z_{\al,i,\jph})^2\bar{u}^y_{\al,i,\jph}\hat{\rho}_{\al,i,\jph}
		\end{pmatrix} 
\end{align*}

Here, $\beta_\al=\frac{\rho_\al}{p_\al}$ and we define the logarithmic average of a strictly positive number as
$$
	\hat{a}_{\iph,j} =\frac{[a]_{\iph,j}}{[\log a]_{\iph,j}}, \qquad
	\hat{a}_{i,\jph} =\frac{[a]_{i,\jph}}{[\log a]_{i,\jph}}
$$ 
When the left and right states get very near to one another, the logarithmic average definition might not work well numerically. We will apply a reliable approximation technique, as detailed in~\cite{Ismail2009}, to deal with this scenario.

The entropy diffusion matrices are given by, 
\begin{equation}
	\label{eq:diffusion_matrices}
	\Db^x_{\al,\iph,j}=\tilde{\mathbf{R}}^x_{\al,\iph,j}\Lambda^x_{\al,\iph,j} \tilde{\mathbf{R}}^{x\top}_{\al,\iph,j}, \qquad \Db^y_{\al,i,\jph}=\tilde{\mathbf{R}}^y_{\al,i,\jph}\Lambda^y_{\al,i,\jph} \tilde{\mathbf{R}}^{y\top}_{\al,i,\jph}
\end{equation}
where $\tilde{\mathbf{R}}^x_{\al}$ and $\tilde{\mathbf{R}}^y_{\al}$ are the right eigenvector matrices for the $x$ and $y-$directional flux jacobian matrices. Following~\cite{Fjordholm2012}, to achieve higher order entropy stable schemes, we reconstruct the scaled entropy variables using the MinMod reconstruction, which has the {\em sign preserving property}. We will only illustrate the procedure for the $x$-direction, as the $y$-directional case is similar. We first define the change of variable,
$$
\mathcal{V}_{\al,k,j}^{x,\pm} =(\tilde{\mathbf{R}}^x_{\al,i\pm\frac{1}{2},j})^\top \Vb_{\al,k,j},
$$
where $k$ are neighbours of cell $(i,j)$ along the $x-$direction. Applying the MinMod reconstruction process on index $k$, we select a cell stencil and construct a linear polynomial $P_{\al,i,j}$, we get the traces, 
$$
\tilde{\mathcal{V}}_{\al,k,j}^{x,\pm}=P_{\al,i,j}^{x,\pm}(x_{i\pm\frac{1}{2}}).
$$
Finally, we define the reconstructed entropy variables as,
$$
\tilde{\Vb}^{x,\pm}_{\al,\iph,j}= \{(\tilde{\mathbf{R}}^x_{\al,i\pm\frac{1}{2},j})^\top\}^{(-1)} \tilde{\mathcal{V}}_{\al,k,j}^{x,\pm}
$$
The second-order entropy stable numerical flux in the $x-$direction is then given by,
\begin{equation}
	\label{eq:ent_stable_numflux_x}
		\Fb^x_{\al,\iph,j} = \tilde{\Fb}^x_{\al,\iph,j}-\frac{1}{2}\Db^x_{\al,\iph,j}[\![\tilde{\Vb}^x_{\al}]\!]_{\iph,j}
\end{equation}
% Similarly, the second-order entropy stable numerical $y-$directional flux is given by,
% \begin{equation}
% 		\label{eq:ent_stable_numflux_y}
% 	\Fb^y_{\al,i,\jph} = \tilde{\Fb}^y_{\al,i,\jph}-\frac{1}{2}\Db^y_{\al,i,\jph}[\tilde{\Vb}_{\al}]_{i,\jph}
% \end{equation}
where $[\![ \mathbf{\tilde{V}}^x_\alpha]\!]_{i+\frac{1}{2},j}$ defined as, $$[\![ \mathbf{\tilde{V}}^x_\alpha]\!]_{i+\frac{1}{2},j}\,=\,\mathbf{\tilde{V}}^{x,-}_{\alpha,{i+1,j}}\,-\,\mathbf{\tilde{V}}^{x,+}_{\alpha,{i,j}}.$$
Using these numerical fluxes, we have the following result from \cite{kumar2012entropy}:
\begin{theorem}
The semi-discrete scheme \eqref{eq:semi_discrete} with numerical fluxes \eqref{eq:ent_stable_numflux_x} is second-order accurate and entropy stable i.e., it satisfies the following entropy inequality,
\begin{equation}
		\frac{d}{dt}  e_\al(\mathbf{U}_{\alpha,i,j})  +\frac{1}{\Delta x} \left( \hat{q}_{\alpha,i+\frac{1}{2},j}^x - \hat{q}_{\alpha,i-\frac{1}{2},j}^x\right)+\frac{1}{\Delta y}\left( \hat{q}_{\alpha,i,j+\frac{1}{2}}^y - \hat{q}_{\alpha,i,j-\frac{1}{2}}^y\right) \le 0 \label{eq:semi_dis_fluid_entropy_inequality}
	\end{equation}
for $\alpha \in \{I,E\}$, where $\hat{q}^{x}_{\alpha,i+\frac{1}{2},j}$, $\hat{q}^{y}_{\alpha, i,j+\frac{1}{2}}$ are the numerical entropy flux functions consistent with the continuous entropy fluxes $q^x_\al, q^y_\al$ respectively, given by,
$$
\hat{q}_{\alpha,i+\frac{1}{2},j}^x = \bar{\textbf{V}}_{\alpha,\iph,j} \cdot \Fb^x_{\al,\iph,j} - \bar{\chi}_{\iph,j}^x
$$
and
$$
\hat{q}_{\alpha,i,\jph}^y = \bar{\textbf{V}}_{\alpha,i,\jph} \cdot \Fb^y_{\al,i,\jph} - \bar{\chi}_{i,\jph}^y.
$$
This holds true for any choice of Maxwell's equation numerical fluxes $\Fb_{M,\iph,j}$ and $\Fb_{M,i,\jph}$.
\end{theorem}

Note that the three fluxes are independent of one another, and the fluid interacts with Maxwell's equations via source terms only. Furthermore, we note that the source terms do not contribute to the entropy production, i.e. 
$$
\Vb_{\al} \cdot \sbb_\al=0.
$$
Hence, we can take any numerical flux for Maxwell's equations, and we still get the above entropy inequality.

In this Section, we have presented spatial discretization for the second-order schemes which are entropy-stable at the semi-discrete level. We will now present the fully discrete numerical schemes and discuss divergence constraint errors for the electromagnetic fields.

\section{Fully discrete numerical schemes and divergence constraints}
\label{sec:fully_dis}
The semi-discrete scheme \eqref{eq:semi_discrete} with the second-order spatial discretization described in the previous section can be written as, 
\begin{equation}
	\frac{d \Ub_{i,j}}{dt} = \mathcal{L}_{i,j}(\Ub(t)) + \sbb(\Ub_{i,j}(t))
\end{equation}
where
$$
\mathcal{L}_{i,j}(\Ub(t))=-\frac{1}{\Dx} \left(\Fb^x_{\iph,j} - \Fb^x_{\imh,j}\right) - \frac{1}{\Dy} \left(\Fb^y_{i,\jph} - \Fb^y_{i,\jmh}\right)  
$$
We first describe the second-order explicit scheme followed by the second-order IMEX scheme.
\subsection{Explicit scheme}
\label{subsec:exp}
To achieve second-order accuracy, we use the second-order strong stability preserving (SSP) Rung-Kutta scheme from \cite{Gottlieb2001}. Let $\mathbf{U}^n$ denote the solution at time $t^n$ and $\mathbf{U}^{n+1}$ denote the solution at time $t^{n+1} \equiv t^n + \Delta t$ where $\Delta t$ is the step size in the time direction. The RK scheme is given by the following steps.
\begin{enumerate}
	\item Set $\Ub^{(0)}=\Ub^n$.
	\item Compute 
	\begin{equation}
		\label{eq:exp_1st_step}
		\Ub^{(1)}_{i,j}=\Ub^{(0)}_{i,j} + \Dt \ \mathcal{L}_{i,j}(\Ub^{(0)}) +\Dt \ \sbb(\Ub_{i,j}^{(0)})
	\end{equation}
	
	\item 
	Define
\begin{equation}
		\label{eq:exp_u2}
			\Ub^{(2)}_{i,j}=\Ub^{(1)}_{i,j}+\Dt \ \mathcal{L}_{i,j}(\Ub^{(1)}) +\Dt \ \sbb(\Ub_{i,j}^{(1)})
\end{equation}	
	Finally set
		\begin{equation}
		\label{eq:exp_2nd_step}
	\Ub^{n+1}_{i,j}=\frac{1}{2}\Ub^{(0)}_{i,j} + \frac{1}{2}\Ub^{(2)}_{i,j}
	\end{equation}
	
\end{enumerate}
The corresponding scheme is denoted by {\bf O2EXP-MultiD}.
We will now analyze the divergence constraints \eqref{eq:divB_evo_cont} and \eqref{eq:divE_evo_cont} for the fully discrete second-order explicit scheme.

Let us define discrete divergence  at $\left(x_\iph,y_\jph\right)$ for the magnetic field $\Bb$
\begin{equation}
\begin{aligned}
	\label{eq:dis_divB}
	\na \cdot \Bb^n_{\iph,\jph} = & \dfrac{1}{2}\left(\frac{B^n_{x,i+1,j+1}-B^n_{x,i,j+1}}{\Dx} +\frac{B^n_{x,i+1,j}-B^n_{x,i,j}}{\Dx}\right) \\
	+&  \dfrac{1}{2}\left(\frac{B^n_{y,i+1,j+1}-B^n_{y,i+1,j}}{\Dy} +\frac{B^n_{y,i,j+1}-B^n_{y,i,j}}{\Dy}\right) 
\end{aligned}
\end{equation}
with similar expressions for $\na \cdot \Eb^n_{\iph,\jph}$ and $\na \cdot \jb^n_{\iph,\jph}$. 
\begin{comment}
for the electric field $\Eb$
\begin{equation}
\begin{aligned}
	\label{eq:dis_divEE}
	\na \cdot \Eb^n_{\iph,\jph} = &\frac{1}{2}\left(\frac{E^n_{x,i+1,j+1}-E^n_{x,i,j+1}}{\Dx} +\frac{E^n_{x,i+1,j}-E^n_{x,i,j}}{\Dx}\right) \\
	+&  \frac{1}{2}\left(\frac{E^n_{y,i+1,j+1}-E^n_{y,i+1,j}}{\Dy} +\frac{E^n_{y,i,j+1}-E^n_{y,i,j}}{\Dy}\right)
 \end{aligned}
\end{equation}
and for the current density $\jb$
\begin{equation}
\begin{aligned}
	\label{eq:dis_divE}
	\na \cdot \jb^n_{\iph,\jph} = &\frac{1}{2}\left(\frac{j^n_{x,i+1,j+1}-j^n_{x,i,j+1}}{\Dx} +\frac{j^n_{x,i+1,j}-j^n_{x,i,j}}{\Dx}\right) \\
	+&  \frac{1}{2}\left(\frac{j^n_{y,i+1,j+1}-j^n_{y,i+1,j}}{\Dy} +\frac{j^n_{y,i,j+1}-j^n_{y,i,j}}{\Dy}\right)
 \end{aligned}
\end{equation}
\end{comment}
Now, we have the following result on the divergence evolution.
\begin{theorem}[Divergence evolution for the explicit scheme] The explicit scheme update \eqref{eq:exp_2nd_step} satisfies,
	\begin{equation}
		\label{eq:exp_divB_error}
			\na \cdot \Bb^{n+1}_{\iph,\jph} = 	\na \cdot \Bb^n_{\iph,\jph}
	\end{equation}
	for the magnetic field $\Bb$ which is consistent with \eqref{eq:divB_evo_cont}. Similarly, for the electric field $\Eb$, we have the discrete divergence evolution given by, 
	\begin{equation}
		\label{eq:exp_divE_error}
		\na \cdot \Eb^{n+1}_{\iph,\jph} = 	\na \cdot \Eb^n_{\iph,\jph} -\frac{\Dt}{2\eps_0} \left( \na \cdot \jb^{n}_{\iph,\jph} + \na \cdot \jb^{(1)}_{\iph,\jph}\right)
	\end{equation}
which is consistent with \eqref{eq:divE_evo_cont}.
\end{theorem}
\begin{proof}
	Considering~\eqref{eq:exp_1st_step} for the magnetic field, we get
	$$
	B_{x,i,j}^{(1)}=	B_{x,i,j}^{(0)} 
	-  \dfrac{\Delta t}{2\Delta y } 
	\bigg[
	\left(\tilde{E}^{(0)}_{z,i+\frac{1}{2},j+\frac{1}{2}} + \tilde{E}^{(0)}_{z,i-\frac{1}{2},j+\frac{1}{2}}\right)
	-\left(\tilde{E}^{(0)}_{z,i+\frac{1}{2},j-\frac{1}{2}} + \tilde{E}^{(0)}_{z,i-\frac{1}{2},j-\frac{1}{2}}\right)
	\bigg]
	$$
	where electric field components $\tilde{E}^{(0)}_{z,i+\frac{1}{2},j+\frac{1}{2}}$ at the edges are computed using multidimensional Riemann solver \eqref{eq:multid_Ez} from the reconstructed values using MinMod limiter. Similarly, we have,
$$
B_{y,i,j}^{(1)}=	B_{y,i,j}^{(0)} 
	+  \dfrac{\Delta t}{2\Delta x } 
	\bigg[
	\left(\tilde{E}^{(0)} _{z,i+\frac{1}{2},j+\frac{1}{2}} + \tilde{E}^{(0)} _{z,i+\frac{1}{2},j-\frac{1}{2}}\right)
	-\left(\tilde{E}^{(0)} _{z,i-\frac{1}{2},j+\frac{1}{2}} + \tilde{E}^{(0)} _{z,i-\frac{1}{2},j-\frac{1}{2}}\right)
	\bigg]
$$
then, we compute the discrete divergence of $\Bb^{(1)}$ as follows:
		\begin{align*}
			\na \cdot \Bb^{(1)} _{\iph,\jph} 
			& =
			\dfrac{1}{2} 
			\Bigg(
			\frac{{B}^{(1)}_{x,i+1,j+1}-{B}^{(1)}_{x,i,j+1} }{\Delta x}
			+
			\frac{{B}^{(1)}_{x,i+1,j}-{B}^{(1)}_{x,i,j} }{\Delta x}
			\Bigg)
			\\ & \qquad +
			\dfrac{1}{2} 
			\Bigg(
			\frac{{B}^{(1)}_{y,i+1,j+1}-{B}^{(1)}_{y,i+1,j} }{\Delta y}
			+
			\frac{{B}^{(1)}_{y,i,j+1}-{B}^{(1)}_{y,i,j} }{\Delta y}
			\Bigg) \\
			& =
			\dfrac{1}{2\Delta x} 
			\Bigg[
			B_{x,i+1,j+1}^{(0)}
			-  \dfrac{\Delta t}{2\Delta y }
			\bigg[
			\left(\tilde{E}^{(0)}_{z,i+\frac{3}{2},j+\frac{3}{2}} + \tilde{E}^{(0)}_{z,i+\frac{1}{2},j+\frac{3}{2}}\right)
			-\left(\tilde{E}^{(0)}_{z,i+\frac{3}{2},j+\frac{1}{2}} + \tilde{E}^{(0)}_{z,i+\frac{1}{2},j+\frac{1}{2}}\right)
			\bigg]\\
			&- 
			\dfrac{1}{2\Delta x} 
			\Bigg[
			B_{x,i,j+1}^{(0)}
			-  \dfrac{\Delta t}{2\Delta y } 
			\bigg[
			\left(\tilde{E}^{(0)}_{z,i+\frac{1}{2},j+\frac{3}{2}} + \tilde{E}^{(0)}_{z,i-\frac{1}{2},j+\frac{3}{2}}\right)
			-\left(\tilde{E}^{(0)}_{z,i+\frac{1}{2},j+\frac{1}{2}} + \tilde{E}^{(0)}_{z,i-\frac{1}{2},j+\frac{1}{2}}\right)
			\bigg]\\
			&+\dfrac{1}{2\Delta x} 
			\Bigg[
			B_{x,i+1,j}^{(0)}
			-  \dfrac{\Delta t}{2\Delta y } 
			\bigg[
			\left(\tilde{E}^{(0)}_{z,i+\frac{3}{2},j+\frac{1}{2}} + \tilde{E}^{(0)}_{z,i+\frac{1}{2},j+\frac{1}{2}}\right)
			-\left(\tilde{E}^{(0)}_{z,i+\frac{3}{2},j-\frac{1}{2}} + \tilde{E}^{(0)}_{z,i+\frac{1}{2},j-\frac{1}{2}}\right)
			\bigg]\\
			&-\dfrac{1}{2\Delta x} 
			\Bigg[
			B_{x,i,j}^{(0)}
			-  \dfrac{\Delta t}{2\Delta y } 
			\bigg[
			\left(\tilde{E}^{(0)}_{z,i+\frac{1}{2},j+\frac{1}{2}} + \tilde{E}^{(0)}_{z,i-\frac{1}{2},j+\frac{1}{2}}\right)
			-\left(\tilde{E}^{(0)}_{z,i+\frac{1}{2},j-\frac{1}{2}} + \tilde{E}^{(0)}_{z,i-\frac{1}{2},j-\frac{1}{2}}\right)
			\bigg]\\
			&+\dfrac{1}{2\Delta y} 
			\Bigg[
			B_{y,i+1,j+1}^{(0)} 
			+  \dfrac{\Delta t}{2\Delta x } 
			\bigg[
			\left(\tilde{E}^{(0)}_{z,i+\frac{3}{2},j+\frac{3}{2}} + \tilde{E}^{(0)}_{z,i+\frac{3}{2},j+\frac{1}{2}}\right)
			-\left(\tilde{E}^{(0)}_{z,i+\frac{1}{2},j+\frac{3}{2}} + \tilde{E}^{(0)}_{z,i+\frac{1}{2},j+\frac{1}{2}}\right)
			\bigg]\\
			&-\dfrac{1}{2\Delta y} 
			\Bigg[
			B_{y,i+1,j}^{(0)} 
			+  \dfrac{\Delta t}{2\Delta x }
			\bigg[
			\left(\tilde{E}^{(0)}_{z,i+\frac{3}{2},j+\frac{1}{2}} + \tilde{E}^{(0)}_{z,i+\frac{3}{2},j-\frac{1}{2}}\right)
			-\left(\tilde{E}^{(0)}_{z,i+\frac{1}{2},j+\frac{1}{2}} + \tilde{E}^{(0)}_{z,i+\frac{1}{2},j-\frac{1}{2}}\right)
			\bigg]\\
			&+\dfrac{1}{2\Delta y} 
			\Bigg[
			B_{y,i,j+1}^{(0)}
			+  \dfrac{\Delta t}{2\Delta x } 
			\bigg[
			\left(\tilde{E}^{(0)}_{z,i+\frac{1}{2},j+\frac{3}{2}} + \tilde{E}^{(0)}_{z,i+\frac{1}{2},j+\frac{1}{2}}\right)
			-\left(\tilde{E}^{(0)}_{z,i-\frac{1}{2},j+\frac{3}{2}} + \tilde{E}^{(0)}_{z,i-\frac{1}{2},j+\frac{1}{2}}\right)
			\bigg]\\
			&-\dfrac{1}{2\Delta y} 
			\Bigg[
			B_{y,i,j}^{(0)} 
			+  \dfrac{\Delta t}{2\Delta x } 
			\bigg[
			\left(\tilde{E}^{(0)}_{z,i+\frac{1}{2},j+\frac{1}{2}} + \tilde{E}^{(0)}_{z,i+\frac{1}{2},j-\frac{1}{2}}\right)
			-\left(\tilde{E}^{(0)}_{z,i-\frac{1}{2},j+\frac{1}{2}} + \tilde{E}^{(0)}_{z,i-\frac{1}{2},j-\frac{1}{2}}\right)
			\bigg]\\
				\end{align*}
			\begin{align*}
			=\na \cdot \Bb^{(0)}_{\iph,\jph}
			+
			\dfrac{\Delta t}{4 \Delta x \Delta y}
			\Bigg[&
			\textcolor{red}{-\tilde{E}^{(0)}_{z,i+\frac{3}{2},j+\frac{3}{2}}}
			 \textcolor{green}{- \tilde{E}^{(0)}_{z,i+\frac{1}{2},j+\frac{3}{2}}}
			\textcolor{blue}{+\tilde{E}^{(0)}_{z,i+\frac{3}{2},j+\frac{1}{2}}} 
			+ \tilde{E}^{(0)}_{z,i+\frac{1}{2},j+\frac{1}{2}}\\
			&\textcolor{green}{+\tilde{E}^{(0)}_{z,i+\frac{1}{2},j+\frac{3}{2}}}
			   \textcolor{brown}{+\tilde{E}^{(0)}_{z,i-\frac{1}{2},j+\frac{3}{2}}}
			-\tilde{E}^{(0)}_{z,i+\frac{1}{2},j+\frac{1}{2}} 
			\textcolor{RubineRed}{- \tilde{E}^{(0)}_{z,i-\frac{1}{2},j+\frac{1}{2}}}\\
			&\textcolor{blue}{-\tilde{E}^{(0)}_{z,i+\frac{3}{2},j+\frac{1}{2}}} 
			- \tilde{E}^{(0)}_{z,i+\frac{1}{2},j+\frac{1}{2}}
			\textcolor{Gray}{+\tilde{E}^{(0)}_{z,i+\frac{3}{2},j-\frac{1}{2}}} 
			\textcolor{Emerald}{+ \tilde{E}^{(0)}_{z,i+\frac{1}{2},j-\frac{1}{2}}}\\
			&+\tilde{E}^{(0)}_{z,i+\frac{1}{2},j+\frac{1}{2}} 
			\textcolor{RubineRed}{+ \tilde{E}^{(0)}_{z,i-\frac{1}{2},j+\frac{1}{2}}}
			\textcolor{Emerald}{-\tilde{E}^{(0)}_{z,i+\frac{1}{2},j-\frac{1}{2}}}
			\textcolor{Melon}{- \tilde{E}^{(0)}_{z,i-\frac{1}{2},j-\frac{1}{2}}}\\
			&\textcolor{red}{+\tilde{E}^{(0)}_{z,i+\frac{3}{2},j+\frac{3}{2}}} 
			\textcolor{blue}{+ \tilde{E}^{(0)}_{z,i+\frac{3}{2},j+\frac{1}{2}}}
			\textcolor{green}{-\tilde{E}^{(0)}_{z,i+\frac{1}{2},j+\frac{3}{2}}} 
			- \tilde{E}^{(0)}_{z,i+\frac{1}{2},j+\frac{1}{2}}\\
			&\textcolor{blue}{-\tilde{E}^{(0)}_{z,i+\frac{3}{2},j+\frac{1}{2}}} 
    			\textcolor{Gray}{- \tilde{E}^{(0)}_{z,i+\frac{3}{2},j-\frac{1}{2}}}
			+\tilde{E}^{(0)}_{z,i+\frac{1}{2},j+\frac{1}{2}} 
			\textcolor{Emerald}{+ \tilde{E}^{(0)}_{z,i+\frac{1}{2},j-\frac{1}{2}}}\\
			&\textcolor{green}{+\tilde{E}^{(0)}_{z,i+\frac{1}{2},j+\frac{3}{2}}} 
			+ \tilde{E}^{(0)}_{z,i+\frac{1}{2},j+\frac{1}{2}}
			\textcolor{brown}{-\tilde{E}^{(0)}_{z,i-\frac{1}{2},j+\frac{3}{2}}} 
			\textcolor{RubineRed}{- \tilde{E}^{(0)}_{z,i-\frac{1}{2},j+\frac{1}{2}}}\\
			&-\tilde{E}^{(0)}_{z,i+\frac{1}{2},j+\frac{1}{2}} 
			\textcolor{Emerald}{- \tilde{E}^{(0)}_{z,i+\frac{1}{2},j-\frac{1}{2}}}
			\textcolor{RubineRed}{+\tilde{E}^{(0)}_{z,i-\frac{1}{2},j+\frac{1}{2}}} 
			\textcolor{Melon}{+ \tilde{E}^{(0)}_{z,i-\frac{1}{2},j-\frac{1}{2}}}
			\Bigg]
		\end{align*}
		We now observe that the same colored terms cancel one another, and we are left with
		$$
		\na \cdot \Bb^{(1)} _{\iph,\jph} =\na \cdot \Bb^{(0)} _{\iph,\jph} 
		$$
Similarly, after the second stage of the RK scheme, we get,
		$$
			\na \cdot \Bb^{(2)} _{\iph,\jph} = 	\na \cdot \Bb^{(1)} _{\iph,\jph} =  \na \cdot \Bb^{(0)} _{\iph,\jph}  
		$$
	Using \eqref{eq:exp_2nd_step}, we get \eqref{eq:exp_divB_error}. Performing a similar calculation for $\na \cdot \Eb^{(1)} _{\iph,\jph} $, we get,
	$$
	\na \cdot \Eb^{(1)} _{\iph,\jph}=	\na \cdot \Eb^{(0)} _{\iph,\jph} -\frac{\Dt}{\eps_0} \left( \na \cdot \jb^{(0)}_{\iph,\jph}\right)
	$$
	and
	$$
	\na \cdot \Eb^{(2)} _{\iph,\jph}=	\na \cdot \Eb^{(1)} _{\iph,\jph} -\frac{\Dt}{\eps_0} \left( \na \cdot \jb^{(1)}_{\iph,\jph}\right)
	$$
	Combining these with the update equation~\eqref{eq:exp_2nd_step}, we get \eqref{eq:exp_divE_error}.
\end{proof}

\begin{remark}
	The above proof is based on the fact that our discretization mimics the vector identities
	$$
	\na \cdot (\na \times \Bb) = \na \cdot( \na \times \Eb)=0
	$$
	at the discrete level (also see Appendix~\ref{sec:Ez}). Furthermore, this property is preserved by the Runge-Kutta time stepping as they are linear invariants.
\end{remark}

\subsection{IMEX scheme}
\label{subsec:imex}
In many cases, the source term can be stiff~\cite{kumar2012entropy}. To overcome this, we use a L-stable second-order accurate IMEX scheme from~\cite{Pareschi2005}. This involves two implicit stages that are given by
\begin{eqnarray}
			\label{eq:imex_1st_step}
	\Ub_{i,j}^{(1)}&=& \Ub_{i,j}^n + \Dt \left(\beta \sbb \left(\Ub_{i,j}^{(1)} \right) \right) \\
			\label{eq:imex_2nd_step}	
	\Ub_{i,j}^{(2)} &=& \Ub_{i,j}^n + \Dt \left[ \mathcal{L}_{i,j} \left(\Ub^{(1)}\right) + (1-2\beta) \sbb \left(\Ub_{i,j}^{(1)} \right) + \beta \sbb \left(\Ub_{i,j}^{(2)} \right) \right]\\
			\label{eq:imex_3rd_step}
\Ub^{n+1}_{i,j} &=&\Ub_{i,j}^n + \frac{1}{2}\Dt\left[ \mathcal{L}_{i,j} \left(\Ub^{(1)}\right) + \mathcal{L}_{i,j} \left(\Ub^{(2)}\right)  + \sbb \left(\Ub_{i,j}^{(1)}\right) +\sbb \left(\Ub_{i,j}^{(2)} \right) \right]
\end{eqnarray} 
In the above update, $\beta$ is a constant given by $\beta = 1 - \frac{1}{\sqrt{2}}$. We follow \cite{Abgrall2014,kumar2012entropy} to solve for $\Ub_{i,j}^{(1)}$ and $\Ub_{i,j}^{(2)}$ in first two stages. We note that the resulting algebraic equations can be solved exactly, and the computational costs are comparable to the explicit schemes above. The corresponding scheme is denoted by {\bf O2IMEX-MultiD}. 

Similar to the explicit scheme, we have the following result  for the IMEX scheme: 
\begin{theorem}[Divergence errors for the IMEX scheme] The IMEX scheme update~\eqref{eq:imex_1st_step}-\eqref{eq:imex_3rd_step} satisfies,
	\begin{equation}
		\label{eq:imex_divB_error}
		\na \cdot \Bb^{n+1}_{\iph,\jph} = 	\na \cdot \Bb^n_{\iph,\jph}
	\end{equation}
	for the magnetic field $\Bb$ which is consistent with \eqref{eq:divB_evo_cont}. Similarly, for the electric field $\Eb$, we have the discrete divergence evolution given by, 
	\begin{equation}
		\label{eq:imex_divE_error}
		\na \cdot \Eb^{n+1}_{\iph,\jph} = 	\na \cdot \Eb^n_{\iph,\jph} -\frac{\Dt}{2\eps_0} \left( \na \cdot \jb^{(1)}_{\iph,\jph} + \na \cdot \jb^{(2)}_{\iph,\jph}\right)
	\end{equation}
	which is consistent with \eqref{eq:divE_evo_cont}.
\end{theorem}
\begin{proof}
Proof for \eqref{eq:imex_divB_error} is similar to the case of the explicit scheme as \eqref{eq:max_B} does not have a source term. We note that from \eqref{eq:imex_1st_step}, we get,
$$
\na \cdot \Eb^{(1)} _{\iph,\jph}=	\na \cdot \Eb^{(0)} _{\iph,\jph} - \beta \frac{\Dt}{\eps_0} \left( \na \cdot \jb^{(1)}_{\iph,\jph}\right)
$$
and update in Eqn.~\eqref{eq:imex_2nd_step} gives,

$$
\na \cdot \Eb^{(2)} _{\iph,\jph}=	\na \cdot \Eb^{(1)} _{\iph,\jph} +(1- 2\beta) \frac{\Dt}{\eps_0} \left( \na \cdot \jb^{(1)}_{\iph,\jph} \right)+ \beta \frac{\Dt}{\eps_0} \left( \na \cdot \jb^{(2)}_{\iph,\jph}\right)
$$
Combining these with \eqref{eq:imex_3rd_step}, we get \eqref{eq:imex_divE_error}.
\end{proof}

In this Section, we have presented the satisfaction of discrete divergence constraints for the fully discrete schemes. We will now present the numerical results for the proposed schemes and discuss the corresponding divergence constraint errors for several two-dimensional test cases.

\section{Numerical results}
\label{sec:num_results}
We will now present various test cases in one and two dimensions. The time step is chosen using,
\[
{\Delta t} = \text{CFL}  \min_{i,j} \left\{ \dfrac{1}{ \dfrac{\Lambda_{max}^x(\mathbf{U}_{i,j})} {\Delta x} + \dfrac{\Lambda_{max}^y(\mathbf{U}_{i,j})}{\Delta y}}: 1 \le i \le N_x, \ 1 \le j \le N_y \right\}
\]
where $\Lambda_{max}^x(\mathbf{U}_{i,j})$ and $\Lambda_{max}^y(\mathbf{U}_{i,j})$ are estimates of the maximum wave speeds in the two directions given by
$$
\Lambda_{max}^x(\mathbf{U}) = \max\{ |u_I^x| + a_I, |u_E^x| + a_E, c\}, \qquad
\Lambda_{max}^y(\mathbf{U}) = \max\{ |u_I^y| + a_I, |u_E^y| + a_E, c\}
$$
We take CFL as $0.45$ for the \textbf{O2IMEX-MultiD} scheme and $0.2$ for the explicit scheme {\bf O2EXP-MultiD}, unless stated otherwise. Several physically significant parameters appear in test cases and are defined as follows:

\begin{itemize}
\item Following~\cite{kumar2012entropy, Abgrall2014}, the normalized ion Larmor radius $\hat{r}_g$:
\begin{align*}
\hat{r}_g = \frac{{r}_g}{x_0} = \frac{m_I {u_I}_0 }{q_I B_0 x_0}
\end{align*}
where ${u_I}_0$ is the reference thermal velocity of the ion, $B_0$ is the reference magnetic field, $x_0$ is the reference length.
\item The normalized ion Debye length $\hat{\lambda}_d$:
\begin{align*}
    \hat{\lambda}_d = \frac{{\lambda}_d}{{r}_g} = \frac{1}{{r}_g}\sqrt{\dfrac{\epsilon_0 {{u_I}_0}^2}{n_0 q_I}}
\end{align*}

\item For the two-dimensional test cases, the charge-to-mass ratio for ions is computed using,
\begin{align*}
	r_I = \frac{1}{d_I \sqrt{\rho_I\mu_0}}
	\qquad
	\text{where} \quad d_I = \dfrac{\text{Domain size}}{\text{Mass ratio}},\qquad \mu_0= 1.0
\end{align*}
Similarly, (see~\cite{wang2020}) the charge-to-mass ratio for electrons, $r_E$ can be found for each test case.
\end{itemize}

In one dimension, the divergence constraints are automatically satisfied by any consistent scheme, therefore we do not discuss the divergence errors in one dimension. For the two-dimensional test cases, we compute the divergence errors for the magnetic field in $L^1$ and $L^2$ norms using the following definitions.
\begin{itemize}
	\item \textbf{$L^1$ error of the magnetic field divergence constraint}: \qquad
	$$
\| \nabla \cdot \Bb^n\|_1= 	\frac{1}{N_x N_y}\sum_{i=1}^{N_x} \sum_{j=1}^{N_y} \left|  (\nabla\cdot\mathbf{B})^{n}_{i+\frac{1}{2},j+\frac{1}{2}}\right|
	$$
	\item \textbf{$L^2$ error of the magnetic field divergence constraint}: \qquad
	$$\| \nabla \cdot \Bb^n\|_2=\left[
	\frac{1}{N_x N_y}\sum_{i=1}^{N_x} \sum_{j=1}^{N_y} \left|  (\nabla\cdot\mathbf{B})^{n}_{i+\frac{1}{2},j+\frac{1}{2}}  \right|^2
	\right]^{1/2} $$
\end{itemize}
For the divergence errors in the electric field constraint, we use the following expressions.
 \begin{itemize}
 	\item \textbf{$L^1$ error of electric field divergence constraint for explicit scheme}:
        {\footnotesize
 	$$
\|\na\cdot \Eb\|_{1}^{E}	 =  \frac{1}{N_x N_y}\sum_{i=1}^{N_x} \sum_{j=1}^{N_y} 
 	\left|  
 	(\nabla\cdot\mathbf{E} )^{n+1}_{i+\frac{1}{2},j+\frac{1}{2}} 
 	-
 	\left(
 	(\nabla\cdot\mathbf{E} )^{n}_{i+\frac{1}{2},j+\frac{1}{2}}
 	- 
 	\dfrac{\Delta t}{2\eps_0}
 	\left[
 	(\nabla\cdot\mathbf{j})^{(n)}_{i+\frac{1}{2},j+\frac{1}{2}}
 	+
 	(\nabla\cdot\mathbf{j})^{1}_{i+\frac{1}{2},j+\frac{1}{2}}
 	\right]
 	\right)
 	\right|
 	$$}
 	\item \textbf{$L^2$ error of electric field divergence constraint for explicit scheme}:
        {\footnotesize
 	$$\|\na\cdot\Eb\|_{2}^{E}=\left[
 	\frac{1}{N_x N_y}\sum_{i=1}^{N_x} \sum_{j=1}^{N_y} 
 	\left|  
 	(\nabla\cdot\mathbf{E} )^{n+1}_{i+\frac{1}{2},j+\frac{1}{2}} 
 	-
 	\left(
 	(\nabla\cdot\mathbf{E} )^{n}_{i+\frac{1}{2},j+\frac{1}{2}}
 	- 
 	\dfrac{\Delta t}{2\eps_0}
 	\left[
 	(\nabla\cdot\mathbf{j})^{(n)}_{i+\frac{1}{2},j+\frac{1}{2}}
 	+
 	(\nabla\cdot\mathbf{j})^{1}_{i+\frac{1}{2},j+\frac{1}{2}}
 	\right]
 	\right)
 	\right|^2
 	\right]^{1/2} $$}
 		\item \textbf{$L^1$error of electric field divergence constraint for IMEX scheme}: 
        {\footnotesize
 	$$
\|\na\cdot\Eb\|_{1}^{I}= 	\frac{1}{N_x N_y}\sum_{i=1}^{N_x} \sum_{j=1}^{N_y} 
 	\left|  
 	(\nabla\cdot\mathbf{E} )^{n+1}_{i+\frac{1}{2},j+\frac{1}{2}} 
 	-
 	\left(
 	(\nabla\cdot\mathbf{E} )^{n}_{i+\frac{1}{2},j+\frac{1}{2}}
 	- 
 	\dfrac{\Delta t}{2\eps_0}
 	\left[
 	(\nabla\cdot\mathbf{j})^{(1)}_{i+\frac{1}{2},j+\frac{1}{2}}
 	+
 	(\nabla\cdot\mathbf{j})^{(2)}_{i+\frac{1}{2},j+\frac{1}{2}}
 	\right]
 	\right)
 	\right|
 	$$}
 	\item \textbf{$L^2$ error of electric field divergence constraint for IMEX scheme}:
        {\footnotesize
 	$$\|\na\cdot\Eb\|_{2}^{I}=\left[
 	\frac{1}{N_x N_y}\sum_{i=1}^{N_x} \sum_{j=1}^{N_y} 
 	\left|  
 	(\nabla\cdot\mathbf{E} )^{n+1}_{i+\frac{1}{2},j+\frac{1}{2}} 
 	-
 	\left(
 	(\nabla\cdot\mathbf{E} )^{n}_{i+\frac{1}{2},j+\frac{1}{2}}
 	- 
 	\dfrac{\Delta t}{2\eps_0}
 	\left[
 	(\nabla\cdot\mathbf{j})^{(1)}_{i+\frac{1}{2},j+\frac{1}{2}}
 	+
 	(\nabla\cdot\mathbf{j})^{(2)}_{i+\frac{1}{2},j+\frac{1}{2}}
 	\right]
 	\right)
 	\right|^2
 	\right]^{1/2} $$ }
 \end{itemize}
%\todo[inline]{
%Why don't we measure this error?
% 	$$
%\|\na\cdot \Eb\|_{1}^{E}	 =  \frac{1}{N_x N_y}\sum_{i=1}^{N_x} \sum_{j=1}^{N_y} 
 %	\left|  
% 	(\nabla\cdot\mathbf{E} )^{n+1}_{i+\frac{1}{2},j+\frac{1}{2}} 
 %	- \frac{1}{\epsilon_0} (\rho_c)_{\iph,\jph}^{n+1}
% 	\right|
 %	$$
%Does the initial condition satisfy
%$$
%\nabla\cdot\mathbf{E}^{0}_{i+\frac{1}{2},j+\frac{1}{2}} 
 %	= \frac{1}{\epsilon_0} (\rho_c)_{\iph,\jph}^{0}
%$$
%}

%\todo[inline,color=green]{(HK)The estimates in error analysis in section 4 are based on the estimates above. At this point, if we change the above error analysis, we need to rewrite those sections and recomputation all the test cases, which will be some work. So, I suggest we go ahead with it for the time being. Also, $$\nabla \cdot\mathbf{E} =\rho_c$$
%is satisfied for all the test cases (in fact, $\rho_c$ is zero), as the 
 %the number densities are always the same for both ions and electrons}

These errors are consistent with the theoretical expression derived in \eqref{eq:exp_divB_error}, \eqref{eq:exp_divE_error}, \eqref{eq:imex_divB_error} and \eqref{eq:imex_divE_error}.

To demonstrate the superiority of the proposed schemes in controlling the divergence constraint errors, we will also compare our two-dimensional results with the two other discretizations of Maxwell's equations. The first one is based on \cite{kumar2012entropy}, where PHM equations are used, and an entropy stable discretization of \eqref{eq:phm_eqn} is used. We use the same time-stepping as given in Section \ref{sec:fully_dis}. The corresponding explicit and IMEX schemes that use the PHM formulations from~\cite{kumar2012entropy} are denoted by {\bf O2EXP-PHM} and {\bf O2IMEX-PHM}, respectively. In another approach, we ignore the divergence constraints \eqref{eq:div_B} and \eqref{eq:gauss_E} and just use entropy stable schemes for \eqref{eq:max_B} and \eqref{eq:max_E}. The corresponding explicit and IMEX schemes are denoted by {\bf O2EXP} and {\bf O2IMEX}, respectively.

\subsection{Accuracy test} 
\label{test:1d_smooth} 
In this test case, we will present the numerical accuracy and convergence rates for the proposed schemes. We follow \cite{kumar2012entropy,Abgrall2014} and consider the forced solution approach, i.e., we solve 
\begin{equation*}
	\frac{\partial \mathbf{U}}{\partial t}
	+ 
	\frac{\partial \mathbf{f}^x}{\partial x}  
	= \mathbf{s} + 
	\mathcal{S}(x,t), 
\end{equation*}
with 
\begin{equation*}
	\mathcal{S}(x,t)=\left(\mathbf{0}_{13},-(2+\sin(2 \pi (x-t))),0, 0 \right)^\top.
\end{equation*}
We set initial ion and electron densities as $\rho_I = \rho_E = 2+\sin(2 \pi (x))$, with initial velocities $u_I^x=u_E^x = 1.0$ and initial pressures $p_I=p_E=1.0$. The $y-$magnetic component is $B_y= \sin(2 \pi (x))$ and the $z-$electric field component is $E_z = -\sin(2\pi (x)$. All other primitive variables are set to zero. We consider $I=[0,1]$ as the computational domain with periodic boundary conditions. To have a non-zero source term in the evolutionary equation for the electric field, we set charge to mass ratios $r_I = 1$ and $r_E = -2$. We use ion-electron adiabatic index $\gamma_I=\gamma_E=5/3$. Under these initial and boundary conditions, it is easy to verify that the exact solution is $\rho_I= \rho_E = 2+\sin(2 \pi (x-t))$. We compute the numerical solutions till a final time $t=2.0$.
\begin{table}[ht]
	\centering
	\begin{tabular}{|c|c|c|c|c|}
		\hline
		Number of cells & \multicolumn{4}{|c|}{ \textbf{O2EXP-MultiD} } \\
		\hline
		-- & $L^1$ error & $L^1$ Order & $L^2$ error & $L^2$ Order  \\
		\hline
		
		32 & 5.91083e-02 &  -- & 1.17787e-01 & --\\ 
		64 & 2.05750e-02 & 1.522 & 4.31622e-02 & 1.448\\ 
		128 & 6.88012e-03 &1.580 & 1.55058e-02 & 1.476 \\
		256 & 1.91456e-03 & 1.845  & 4.97168e-03 & 1.641\\ 
		512 & 5.25264e-04 & 1.865 &  1.57020e-04 & 1.662\\
		1024 & 1.41491e-04 & 1.892 & 4.94566e-04 & 1.666\\
		2048 & 3.73855e-05 & 1.920 &  1.55928e-04 & 1.665\\
		4096 & 9.74487e-06 & 1.939 &  4.92118e-05 & 1.663\\
		8192 & 2.51444e-06 & 1.954 & 1.55401e-05 & 1.664\\
		\hline
	\end{tabular}
	\caption[ht]{\nameref{test:1d_smooth}: $L^1$ and $L^2$ errors and order of convergence for $\rho_I$ using \textbf{O2EXP-MultiD}.}
	\label{tab:acc1}
\end{table}	
\begin{table}[ht]
	\centering
	\begin{tabular}{|c|c|c|c|c|}
		\hline
		Number of cells &
		\multicolumn{4}{|c|}{{\textbf{O2IMEX-MultiD}}} \\
		\hline
		-- & $L^1$ error & $L^1$ Order & $L^2$ error & $L^2$ Order  \\
		\hline
		32 & 5.91047e-02 &--& 1.17783e-01 & -- \\ 
		64 & 2.05747e-02 & 1.522& 4.31618e-02 & 1.448 \\
		128 &  6.88012e-03 & 1.580 & 1.55057e-02 & 1.476 \\
		256 & 1.91453e-03 & 1.845& 4.97167e-03 & 1.641 \\
		512 & 5.25264e-04 & 1.865& 1.57020e-03 & 1.662 \\
		1024 &  1.41491e-04 & 1.892& 4.94566e-04 & 1.666 \\
		2048 & 3.73854e-05 & 1.920& 1.55928e-04 & 1.665 \\
		4096 & 9.74483e-06 & 1.939& 4.92118e-05 & 1.663\\
		8192 &  2.51443e-06 & 1.954& 1.55401e-05 & 1.663\\
		\hline
	\end{tabular}
	\caption[ht]{\nameref{test:1d_smooth}: $L^1$ and $L^2$ errors and order of convergence for $\rho_I$ using \textbf{O2IMEX-MultiD} schemes.}
	\label{tab:acc2}
\end{table}	

{In Tables~\ref{tab:acc1} and~\ref{tab:acc2}, we have presented the $L^1$ and $L^2$ errors and order of convergence for ion density $\rho_I$ using \textbf{O2EXP-MultiD} and \textbf{O2IMEX-MultiD} schemes, respectively at the final time $t=2.0$. We observe that both schemes achieve second-order convergence for $L^1$ errors. However, for $L^2$ errors, we observe less than second order of accuracy which is due to the presence of TVD-type limiters in  entropy diffusion operator. When we switch off the diffusion operator in the fluid equations (in that case, the schemes are only entropy conservative), we observe second order in $L^2$ norm as well but those results are not shown here; this behaviour is expected and is also observed in \cite{xu_high-order_2024}.}

\subsection{Brio-Wu shock tube problem} 
\label{test:1d_bw}
Following \cite{Brio1988,Hakim2006,loverich2011,kumar2012entropy,Abgrall2014}, in this test case, we consider a generalization of the Brio-Wu shock tube test problem for ideal MHD. We consider the computational domain $[0,1]$ with outflow boundary conditions. The initial condition consists of two constant states $\Ub_L, \Ub_R$ with a discontinuity at $x=0.5$. The states are given by,
\begin{align*}
\label{test:1d_bw_init_cond}
	\mathbf{U}_{L}=\begin{cases}
		\rho_I=1\\
		\rho_E=m_E/m_I\\
		p_I=p_E=5\times 10^{-5}\\
		\mathbf{u}_I=\mathbf{u}_E=0\\
		B^x=0.75, B^y=1,B^z=0\\
		\mathbf{E}=0\\
		\phi=\psi=0
	\end{cases}\quad \mathbf{U}_{R}=\begin{cases}
		\rho_I=0.125\\
		\rho_E=0.125m_E/m_I\\
		p_I=p_E=5\times 10^{-6}\\
		\mathbf{u}_I=\mathbf{u}_E=0\\
		B^x=0.75, B^y=-1,B^z=0\\
		\mathbf{E}=0\\
		\phi=\psi=0
	\end{cases}
\end{align*}
The ion-electron mass ratio is taken to be $1836$ with non-dimensional Debye length $0.01$. The gas constants are $\gamma_I=\gamma_E=5/3$. We consider the cases with Larmor radii, $\hat{r}_g$ of $0.1$ and $0.001$, which correspond to ion charge to mass ratios of $10$ and $1000$, respectively. Following \cite{Hakim2006,loverich2011,kumar2012entropy,Abgrall2014}, we expect the solution for Larmor radius $0.001$ to be closer to the ideal MHD solution, with additional dispersive effects.

\begin{figure}[!htbp]
	\begin{center}
		\subfigure[Plot of ion density $\rho_I$ for Larmor radius $0.1$.]{
			\includegraphics[width=0.46\textwidth,clip=]{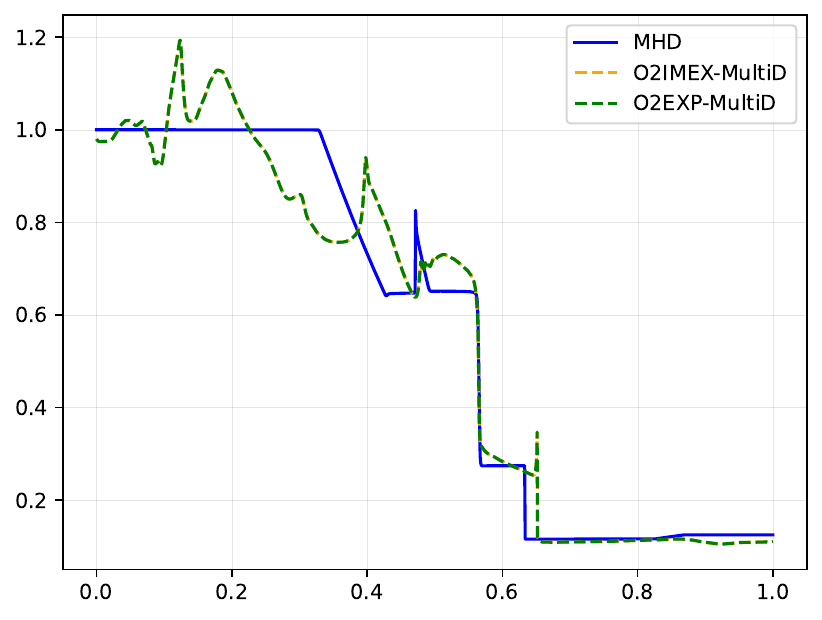}
			\label{fig:bw_re1}}			
		\subfigure[Plot of ion density $\rho_I$ for Larmor radius $0.001$.]{
			\includegraphics[width=0.46\textwidth,clip=]{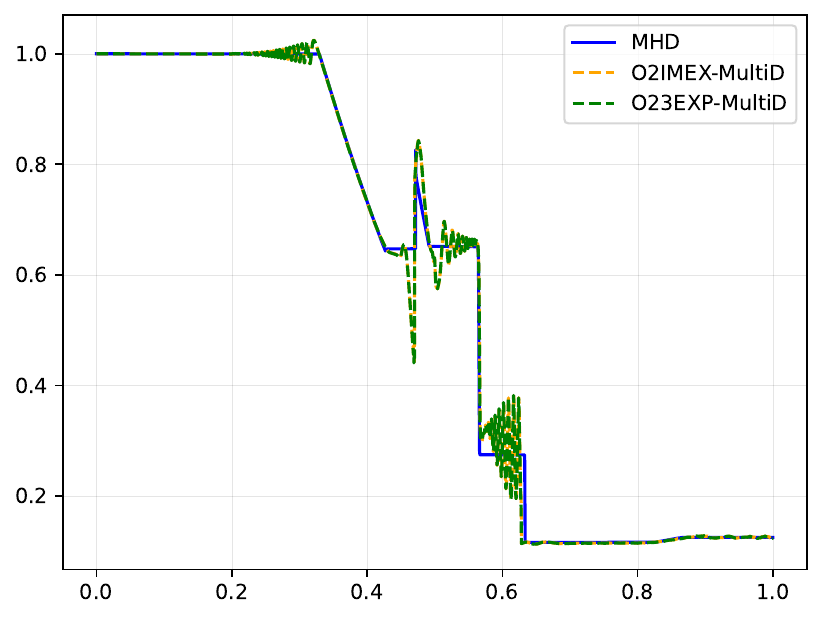}
			\label{fig:bw_re3}}
		\caption{\nameref{test:1d_bw}: Plots of ion density $\rho_I$ for Larmor radius $0.1$ on $10000$ cells and Larmor radius $0.001$ on $50000$ cells.}
		\label{fig:bw}
	\end{center}
\end{figure}

In Figure~\ref{fig:bw}, we present the numerical results at the final time $t=0.1$. We note that both schemes produce similar results. Furthermore, the numerical results are similar to those presented in \cite{Hakim2006,Abgrall2014}. For the case of Larmor radius $0.001$, the second order explicit scheme produces non-physical oscillations, similar to \cite{kumar2012entropy}. Hence, following \cite{kumar2012entropy}, we use SSP third order time stepping (the corresponding scheme is denoted by {\bf O23EXP-MultiD}), for that case. With this change, we note in Figure \ref{fig:bw}b that, using a highly refined mesh of $50000$ cells for the Larmor radius of $0.001$, we are able to resolve MHD waves with additional dispersive effects, as expected.

%In Fig.~\ref{fig:st1_re3}, \ref{fig:st1_re1}, we have shown the numerical results with the $\textbf{O2-ES-MULTID-IMEX}$ and $\textbf{O2-ES-MULTID-EXP}$ schemes using second order IMEX time updates and in both the cases we have compared our results with MHD results.

\subsection{Soliton propagation in one dimension} 
\label{test:1d_sol}
Simulations of soliton propagation in two-fluid plasma are presented in \cite{Baboolal2001,Hakim2006,kumar2012entropy,Abgrall2014}. We consider a one-dimensional computational domain given by the interval $[0,L]$, where $L=12$. We use periodic boundary conditions at both boundaries. The initial ion density profile is taken to be 
\[
\rho_I = 1.0 + \exp\left(-25.0\left|x - \frac{L}{3.0}\right|\right).
\]
We set the ion-electron mass ratio to $\frac{m_I}{m_E} = 25$. The electron pressure is taken to be $p_E = 5.0\rho_I$ with an ion-electron pressure ratio of $1/100$. We take $\gamma_I=\gamma_E=5/3.$ Following \cite{Abgrall2014}, we take normalized Debye length of 1.0 and consider the cases with Larmor radii, $\hat{r}_g$ of $10^{-2}$, $10^{-4}$ and $10^{-6}$. 

\begin{figure}[!htbp]
	\begin{center}
		\subfigure[Plot of ion density $\rho_I$ with Larmor radius $10^{-2}$.]{
			\includegraphics[width=1.9in, height=1.7in]{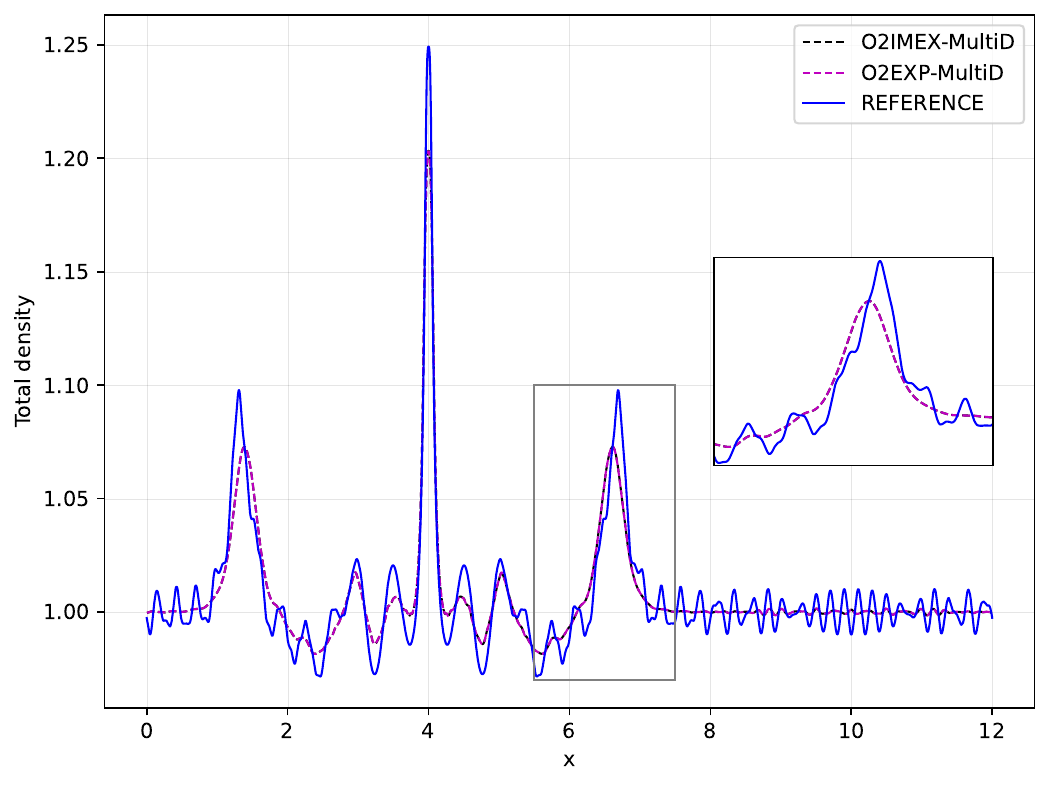}
			\label{fig:sol_1}}
		\subfigure[Plot of ion density $\rho_I$ with Larmor radius $10^{-4}.$]{
			\includegraphics[width=1.9in, height=1.7in]{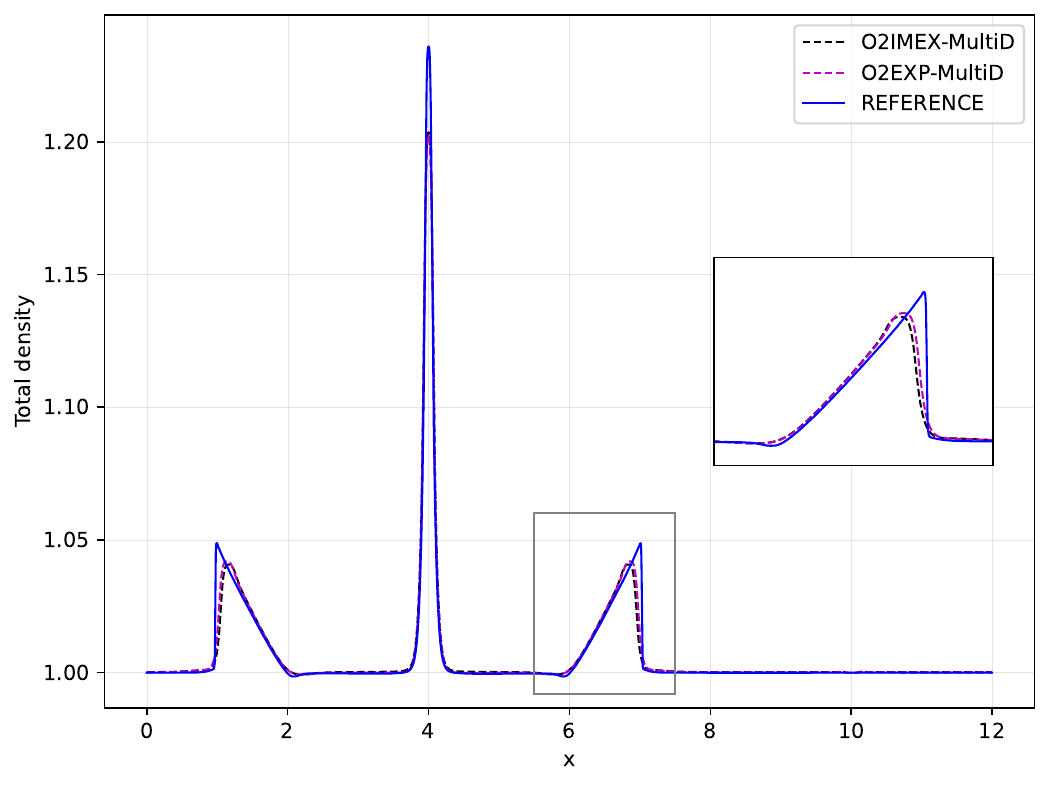}
			\label{fig:sol_2}}
		\subfigure[Plot of ion density $\rho_I$ with Larmor radius $10^{-6}.$]{
			\includegraphics[width=1.9in, height=1.7in]{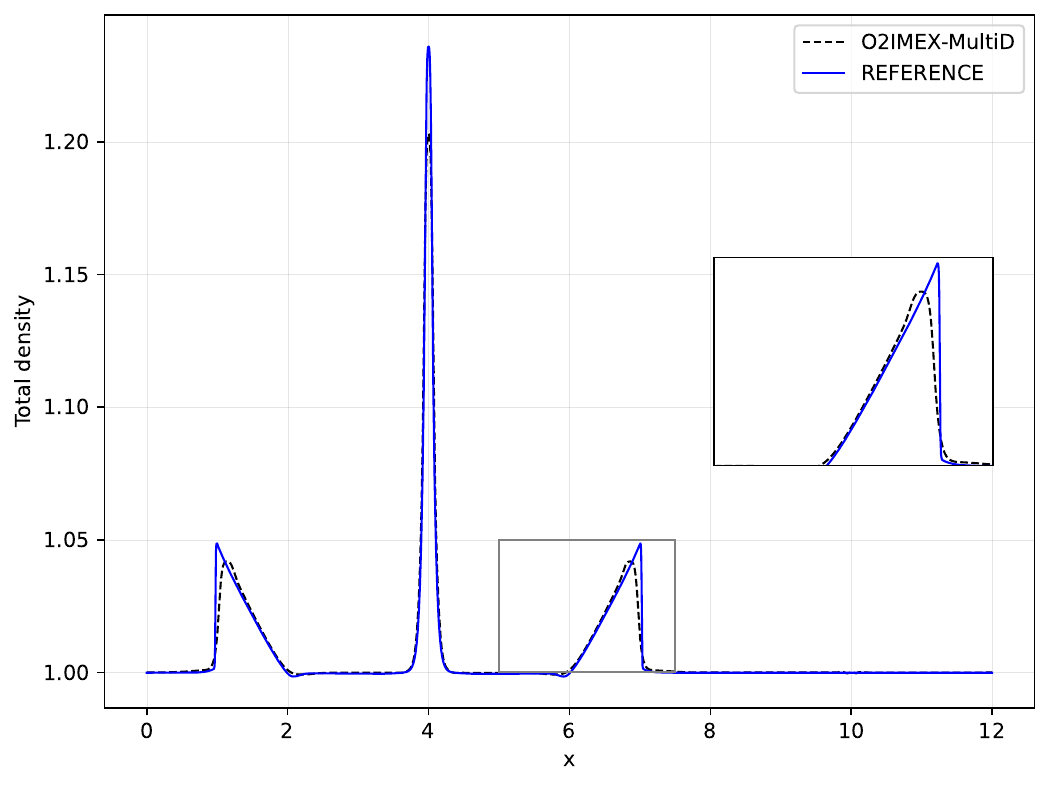}
			\label{fig:sol_3}}
		\caption{\nameref{test:1d_sol}: Plots of ion density using the \textbf{O2IMEX-MultiD } and \textbf{O2EXP-MultiD} schemes with $1500$ cells.}
		\label{fig:sol}
	\end{center}
\end{figure}

In Figure~\ref{fig:sol}, we have presented the solutions from {\bf O2EXP-MultiD}
 and {\bf O2IMEX-MultiD} schemes on $1500$ cells at the final time $t=5.0$. We have also plotted the reference solutions computed using the second-order IMEX scheme (O2-222-imex) in \cite{Abgrall2014} with $20000$ cells. We have not plotted the {\bf O2EXP-MultiD} result for Larmor radius $10^{-6}$ as the source term is highly stiff in this case, and we are not able to compute the results for {\bf O2EXP-MultiD}, similar to the case in \cite{Abgrall2014}. Both schemes produce similar results to those presented in \cite{Abgrall2014}.

\begin{table}[H]
	\renewcommand{\arraystretch}{1.5}
	\centering
	\captionsetup{justification=raggedright,singlelinecheck=false} % Left-align caption	
	\begin{tabular}{cccc}
		\hline
		Schemes & $\hat{r}_g = 10^{-2}$ & $\hat{r}_g = 10^{-4}$  & $\hat{r}_g = 10^{-6}$\\
		\hline
		\textbf{O2IMEX-MultiD} & 176.91 & 174.88 & 173.47\\
		\textbf{O2EXP-MultiD} & 159.12 & 61779.01 & -\\
		\hline
	\end{tabular}
		\caption{Comparison of simulation times of the numerical schemes \textbf{O2EXP-MultiD} and \textbf{O2IMEX-MultiD} for Larmor radii of $10^{-2}, 10^{-4}~\text{and}~10^{-6}$ using 1500 cells.} 
		\label{tab:sol_time}
\end{table}
 
 In Table~\ref{tab:sol_time}, we have presented the computational time for both the schemes using different values of Larmor radii. We have used 4 CPU cores. We note that as the Larmor radius decreases, the stiffness in the source term increases. At the Larmor radius of $10^{-2}$ the source terms are not stiff, and we can see that both {\bf O2IMEX-MultiD} and {\bf O2EXP-MultiD} schemes with  CFL 0.4 are stable. Furthermore, they have similar computational time. This shows that the exact solution approach for the source term discretization in~\cite{Abgrall2014} is highly effective and does not lead to substantial increase in computational time. As we decrease the Larmor radius to $10^{-4}$, the source term becomes more stiff, and we have to decrease the CFL to 0.001 for the {\bf O2EXP-MultiD} scheme. Consequently, we see that the computational time for the {\bf O2EXP-MultiD} has increased substantially, but the {\bf O2IMEX-MultiD} scheme has a similar computational time as for the case of $\hat r_g = 10^{-2}$. For the Larmor radius of $10^{-6}$, the source term is very stiff, and {\bf O2EXP-MultiD} scheme is not able to produce the result. Nevertheless, the {\bf O2IMEX-MultiD} scheme is still stable with CFL $0.4$ and has the same computational time as the other two cases. These results are similar to those presented in \cite{Abgrall2014}.

\subsection{Orszag-Tang vortex} \label{test:2d_ot6}
Orszag-Tang vortex problem was originally considered in~\cite{Orszag1979} for MHD and for the two-fluid equations, the test case was presented in~\cite{wang2020}. We consider a two-dimensional computational domain $[0,2\pi] \times [0,2\pi]$ with periodic boundary conditions. We set adiabatic gas constants as $\gamma_I=\gamma_E=5/3$ and mass ratio $m_I/m_E=25$. Considering the uniform total mass density of $\rho=25/9$, the initial conditions are given by,

\begin{align*}
	\begin{pmatrix*}[c]
		\rho_I \\ u_I^x \\ u_I^y \\ p_I
	\end{pmatrix*} =
	\begin{pmatrix*}[c]
		\frac{25}{26}\rho  \\- {\sin(y)} \\ {\sin(x)} \\ \gamma_I/2
	\end{pmatrix*}, \qquad
	\begin{pmatrix*}[c]
		\rho_E \\u_E^x \\ u_E^y \\ p_E
	\end{pmatrix*} =
	\begin{pmatrix*}[c]
		\frac{1}{26}\rho \\ -\sin(y) \\ \sin(x) \\ \gamma_E/2
	\end{pmatrix*}, \qquad
	\begin{pmatrix*}[c]
		B_x \\ B_y
	\end{pmatrix*} =
	\begin{pmatrix*}[c]
		-\sin(y) \\  \sin (2 x)
	\end{pmatrix*}.	
\end{align*}
The initial electric field is given by Ohm's Law $\mathbf{E} = -\frac{1}{c} \mathbf{u} \times \mathbf{B}$ and all other variables are set to zero. We use the charge-to-mass ratios of $r_I= 2.434602$, $r_E=-60.865062$ and simulate till final time of $t = \pi $ using $400\times 400$ cells.

\begin{figure}[htbp]
	\begin{center}
		\subfigure[Total density $\rho_I +\rho_E$ with \textbf{O2EXP-MultiD} scheme.]{
			\includegraphics[width=0.45\textwidth,clip=]{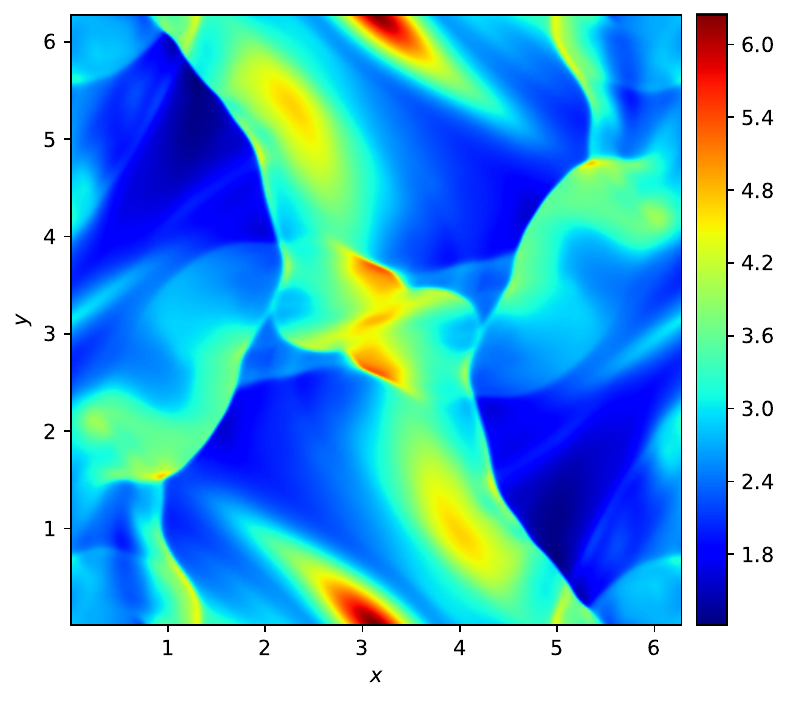}
			\label{fig:ot_o2_exp_multid_rho}}
		\subfigure[Total density $\rho_I +\rho_E$ with \textbf{O2IMEX-MultiD} scheme.]{
			\includegraphics[width=0.45\textwidth,clip=]{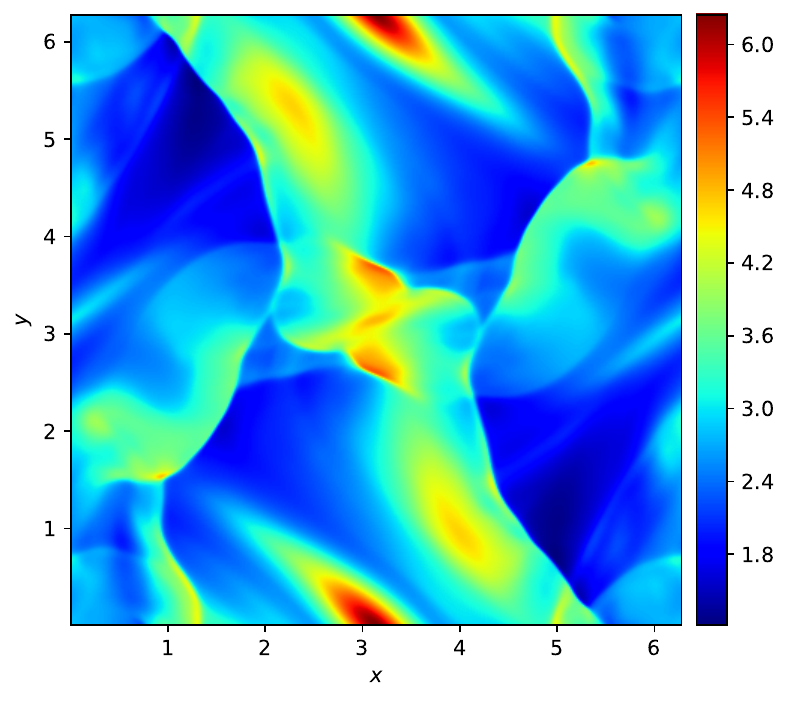}
			\label{fig:ot_o2_imp_multid_rho}}
		\subfigure[Total pressure $p_{i} + p_E$ with \textbf{O2EXP-MultiD} scheme.]{
	\includegraphics[width=0.45\textwidth,clip=]{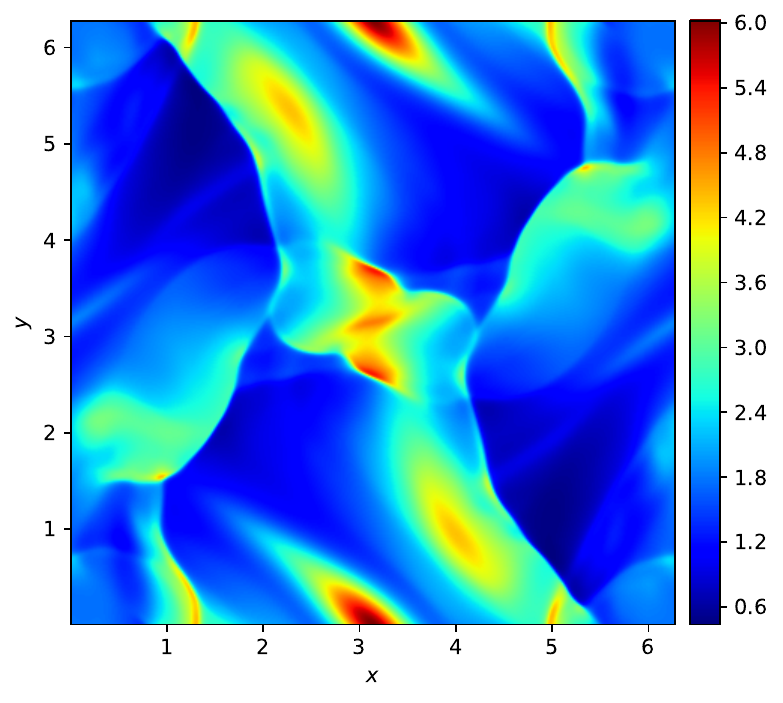}
	\label{fig:ot_o2_exp_multid_p}}
		\subfigure[Total pressure $p_{i} + p_E$ with  \textbf{O2IMEX-MultiD} scheme.]{
			\includegraphics[width=0.45\textwidth,clip=]{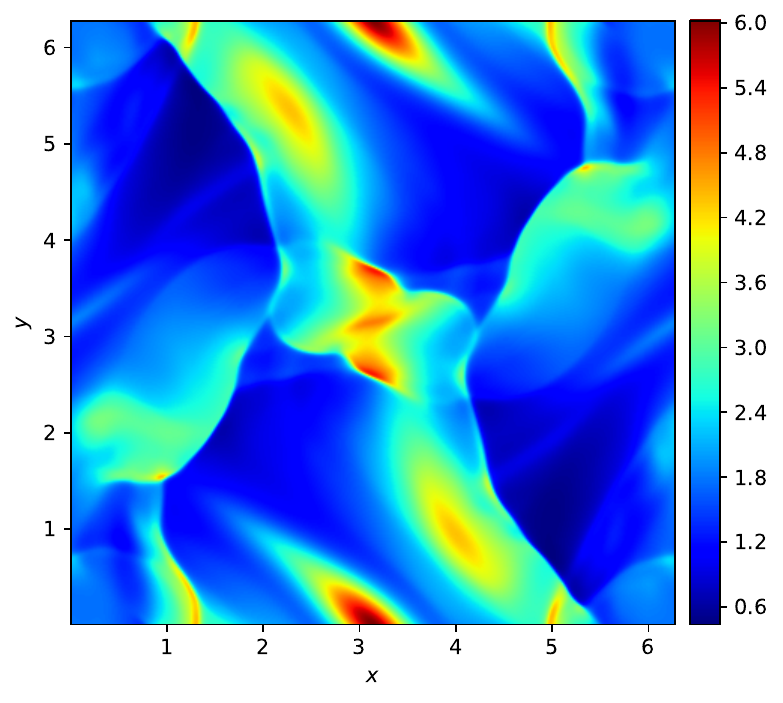}
			\label{fig:ot_o2_imp_multid_p}}
			\subfigure[Magnitude of the Magnetic field, $|\mathbf{B}|$ with \textbf{O2EXP-MultiD} scheme.]{
				\includegraphics[width=0.45\textwidth,clip=]{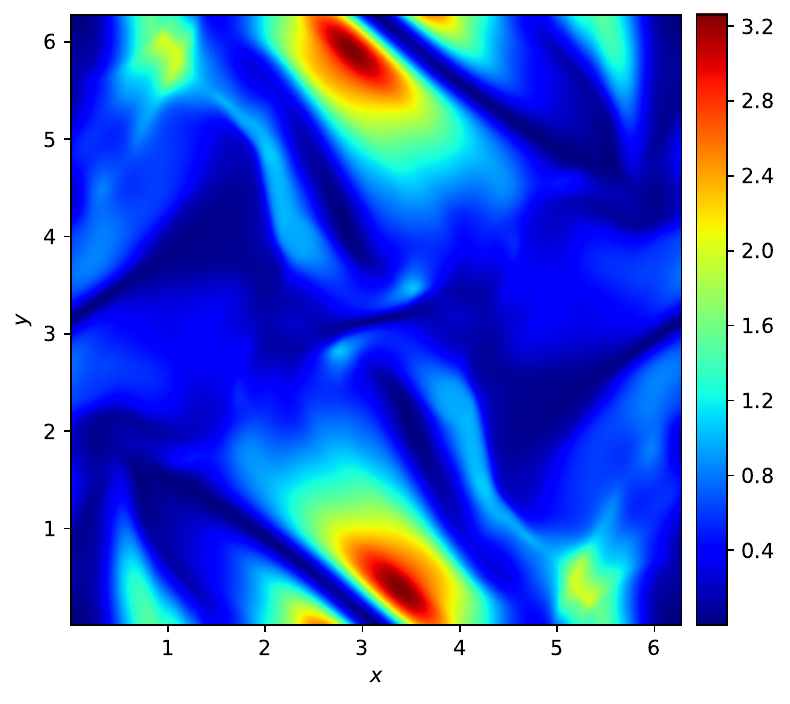}
				\label{fig:ot_o2_exp_multid_MagB}}
		\subfigure[Magnitude of the Magnetic field, $|\mathbf{B}|$ with  \textbf{O2IMEX-MultiD} scheme.]{
			\includegraphics[width=0.45\textwidth,clip=]{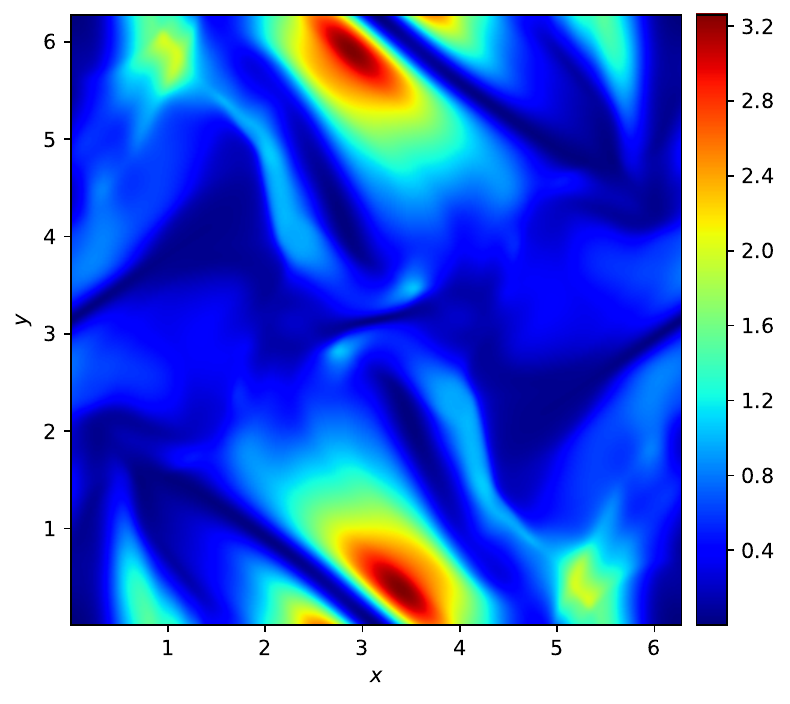}
			\label{fig:ot_o2_imp_multid_MagB}}
		\caption{\nameref{test:2d_ot6}: Plots of total density, total pressure and $|\mathbf{B}|$ with $400\times 400$ cells at time $t=3.14$.}
		\label{fig:ot_multid}
	\end{center}
\end{figure}

%Pressure Cut

\begin{figure}[htbp]
	\begin{center}
\includegraphics[width=0.8\textwidth,clip=]{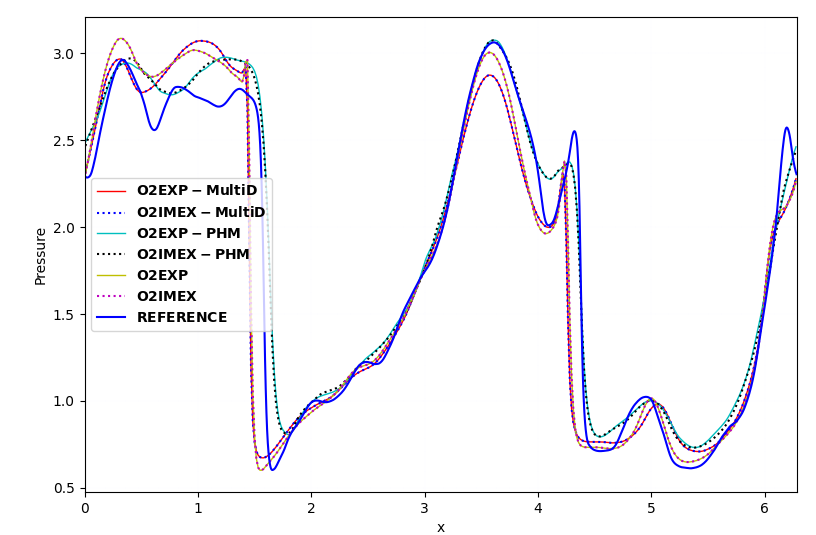}
			\label{fig:ot_o2_pressure}
		\caption{\nameref{test:2d_ot6}: Plot of pressure cut along $y = 1.9634$ for Orszag-Tang test compared with for all schemes.}
		\label{fig:ot_pressure_cut}
	\end{center}
\end{figure}

\begin{figure}[htbp]
	\begin{center}
		\subfigure[$\|\na \cdot \mathbf{B}^n\|_1$,  and $\|\na \cdot \mathbf{B}^n\|_2$  errors for explicit schemes  \textbf{O2EXP-MultiD}, \textbf{O2EXP-PHM} and \textbf{O2EXP}.]{
			\includegraphics[width=0.45\textwidth,clip=]{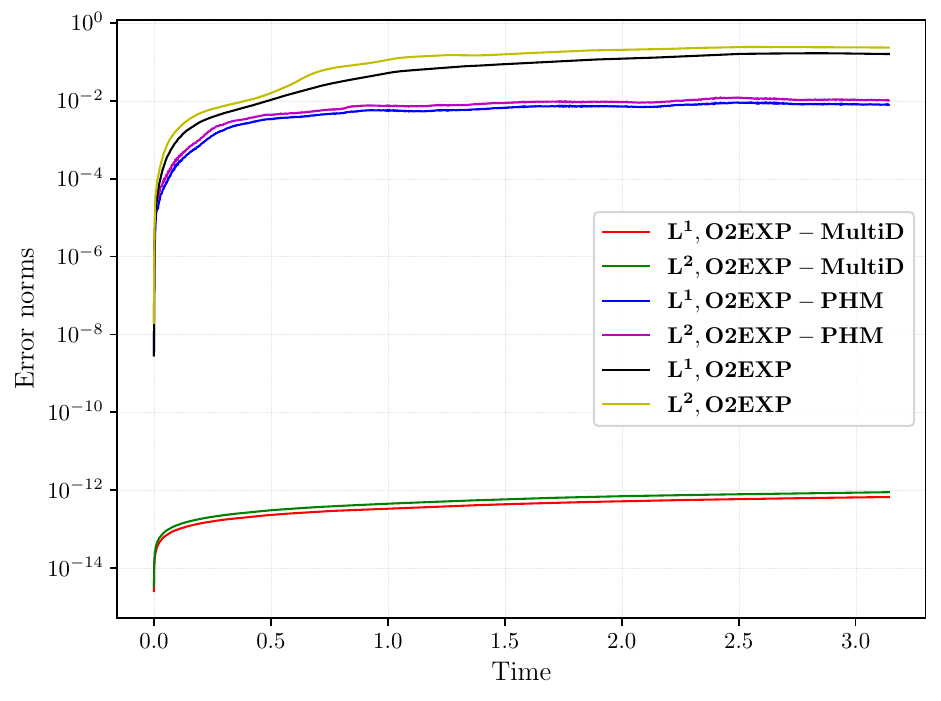}
			\label{fig:ot_exp_divB}}
			\subfigure[$\|\na \cdot \mathbf{B}^n\|_1$,  and $\|\na \cdot \mathbf{B}^n\|_2$  errors for IMEX schemes \textbf{O2IMEX-MultiD}, \textbf{O2IMEX-PHM} and \textbf{O2IMEX}.]{
				\includegraphics[width=0.45\textwidth,clip=]{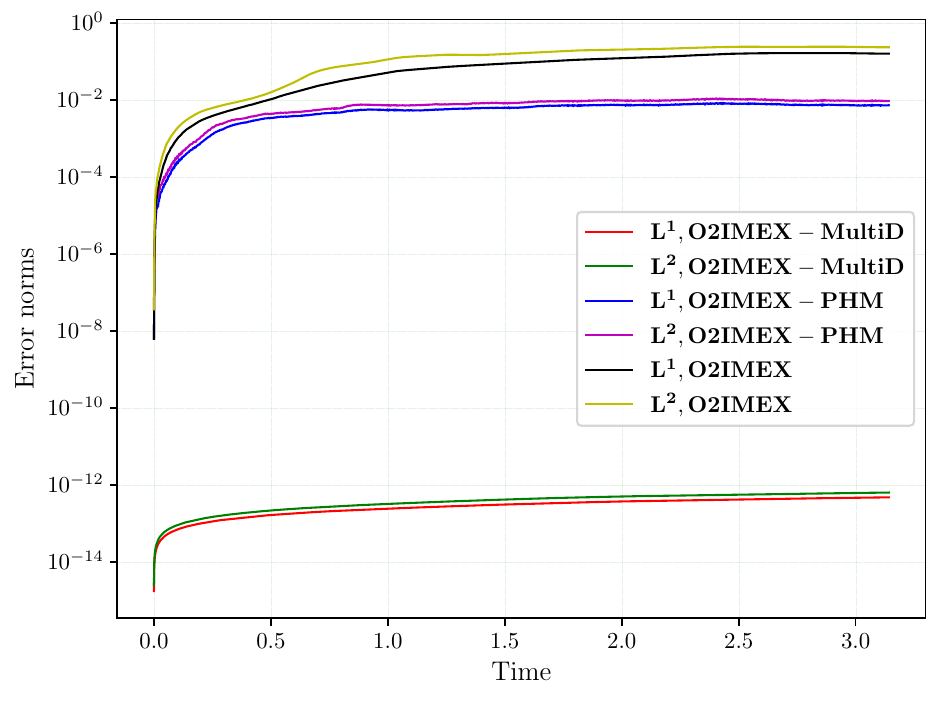}
				\label{fig:ot_imex_divB}}
			\subfigure[$\|\na\cdot\Eb\|_{1}^{E},$ and $\|\na\cdot\Eb\|_{2}^{E}$ errors, for explicit schemes  \textbf{O2EXP-MultiD}, \textbf{O2EXP-PHM} and \textbf{O2EXP}. ]{\includegraphics[width=0.45\textwidth,clip=]{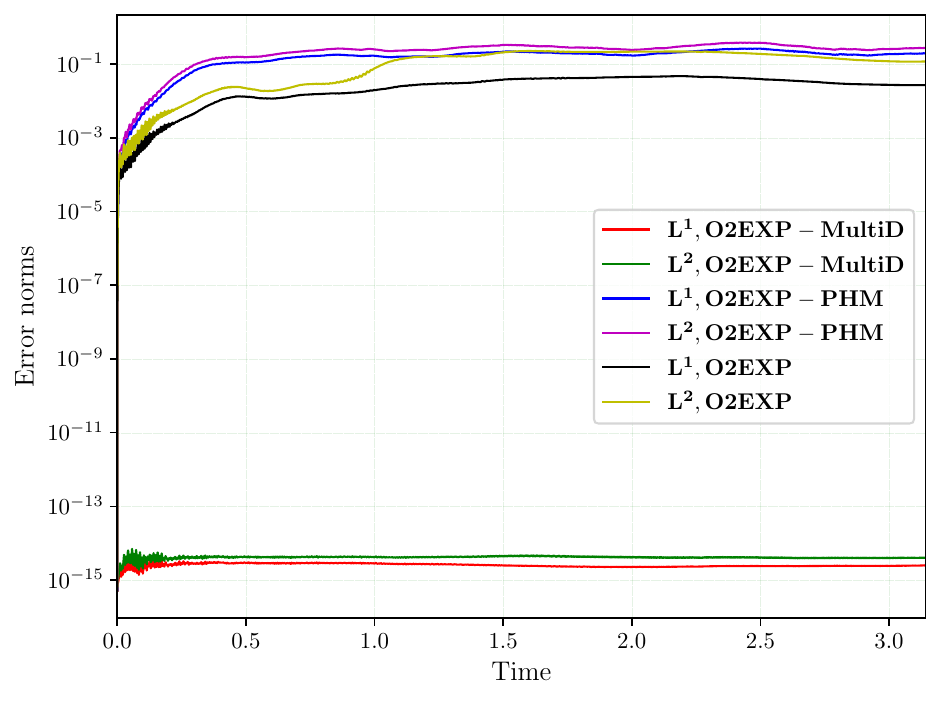}
				\label{fig:ot_exp_divE}}
		\subfigure[$\|\na\cdot\Eb\|_{1}^{I},$ and $\|\na\cdot\Eb\|_{2}^{I}$ errors, for IMEX schemes \textbf{O2IMEX-MultiD}, \textbf{O2IMEX-PHM} and \textbf{O2IMEX}.]{
			\includegraphics[width=0.45\textwidth,clip=]{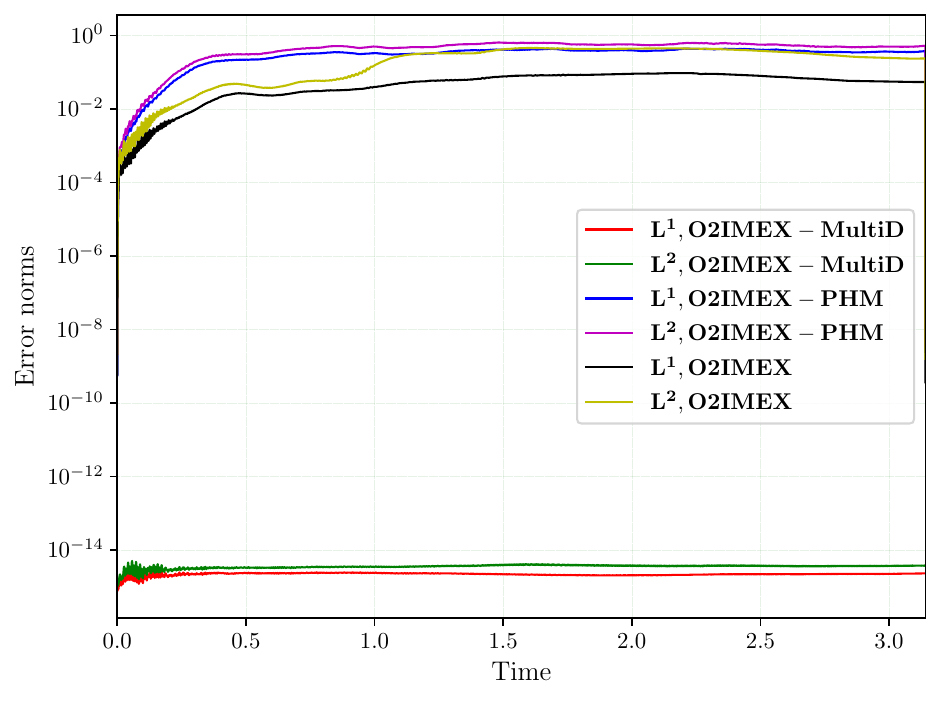}
			\label{fig:ot_imex_divE}}
	\caption{\nameref{test:2d_ot6}: Evolution of the divergence constraint errors for \textbf{O2IMEX-MultiD},  \textbf{O2EXP-MultiD}, \textbf{O2IMEX-PHM}, \textbf{O2EXP-PHM},  \textbf{O2IMEX} and \textbf{O2EXP} schemes using $400\times 400$ cells.}
		\label{fig:ot_div}
	\end{center}
\end{figure}
\begin{figure}[!htbp]
	\begin{center}
		\includegraphics[width=3.0in, height=2.3in]{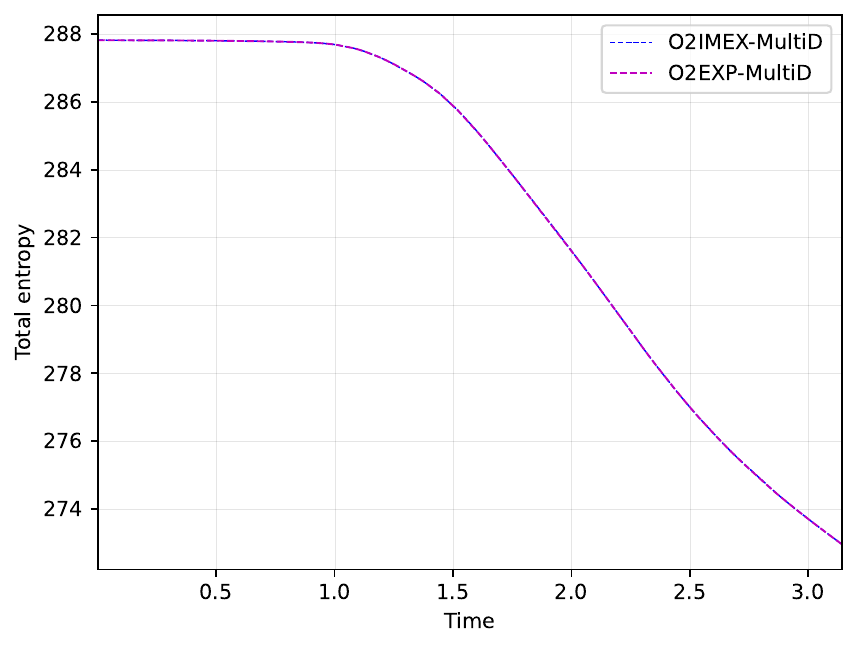}
		\caption{\nameref{test:2d_ot6} Total fluid entropy evolution for the schemes {\bf O2EXP-MultiD} and {\bf O2IMEX-MultiD}.}
		\label{fig:ot_entropy}
	\end{center}
\end{figure}
%\todo{This is stiff or non-stiff case, can we also do stiff case}
In Figure~\ref{fig:ot_multid}, we have plotted the total density, total pressure, and magnitude of the magnetic field  $|\Bb|$ for \textbf{O2EXP-MultiD} and \textbf{O2IMEX-MultiD} schemes. We observe that both the schemes are stable, and we are able to resolve all the waves. The results of both the schemes are comparable to those in \cite{wang2020}. In Figure \ref{fig:ot_pressure_cut}, we have plotted the pressure along the line $y=1.9634$ for different schemes. We have also plotted the reference solution, which is computed using the second order IMEX scheme (O2-222-imex) from \cite{Abgrall2014} with $1024\times 1024$ cells. We note that explicit and the corresponding IMEX scheme have similar accuracy for each discretization. We also note that the \textbf{O2EXP-PHM} and \textbf{OEIMEX-PHM} schemes are much more diffusive than the other schemes. In Figure~\ref{fig:ot_div}, we have plotted the evolution of divergence constraint $L^1$ and $L^2$ errors for, \textbf{O2EXP-PHM}, \textbf{O2IMEX-PHM}, \textbf{O2EXP} and \textbf{O2IMEX} schemes, in addition to the \textbf{O2EXP-MultiD} and \textbf{O2IMEX-MultiD} schemes. We note that the proposed schemes \textbf{O2EXP-MultiD} and \textbf{O2IMEX-MultiD}, ensure divergence constraint errors of machine precision consistent with \eqref{eq:divB_evo_cont} and \eqref{eq:divE_evo_cont}. Other schemes have significantly more divergence constraint errors. For the divergence of magnetic field errors, we observe that the PHM-based schemes \textbf{O2EXP-PHM} and \textbf{O2IMEX-PHM} have slightly lower errors compared to the \textbf{O2EXP} and \textbf{O2IMEX} schemes. However, for the divergence errors related to the electric field, both schemes have similar errors.

We have also plotted the time evolution of total fluid entropy for both schemes in Figure \ref{fig:ot_entropy}. We observe that initially, the flow is smooth, and hence, there is no entropy decay. However, at around time $t=1.0$, shocks start to form, and both the schemes start decaying entropy.

\subsection{Rotor problem} \label{test:2d_rotor}
\begin{figure}[!htbp]
	\begin{center}
		\subfigure[Total density $\rho_I +\rho_E$ with \textbf{O2EXP-MultiD} scheme.]{
			\includegraphics[width=0.45\textwidth,clip=]{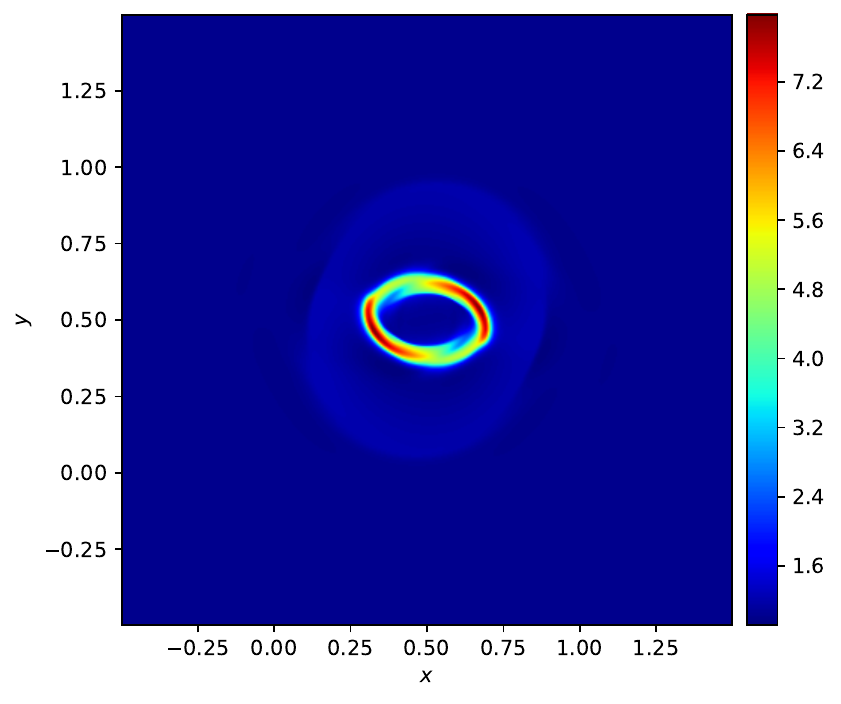}
			\label{fig:density_512_exp}}
		\subfigure[Total density $\rho_I +\rho_E$ with \textbf{O2IMEX-MultiD} scheme.]{
			\includegraphics[width=0.45\textwidth,clip=]{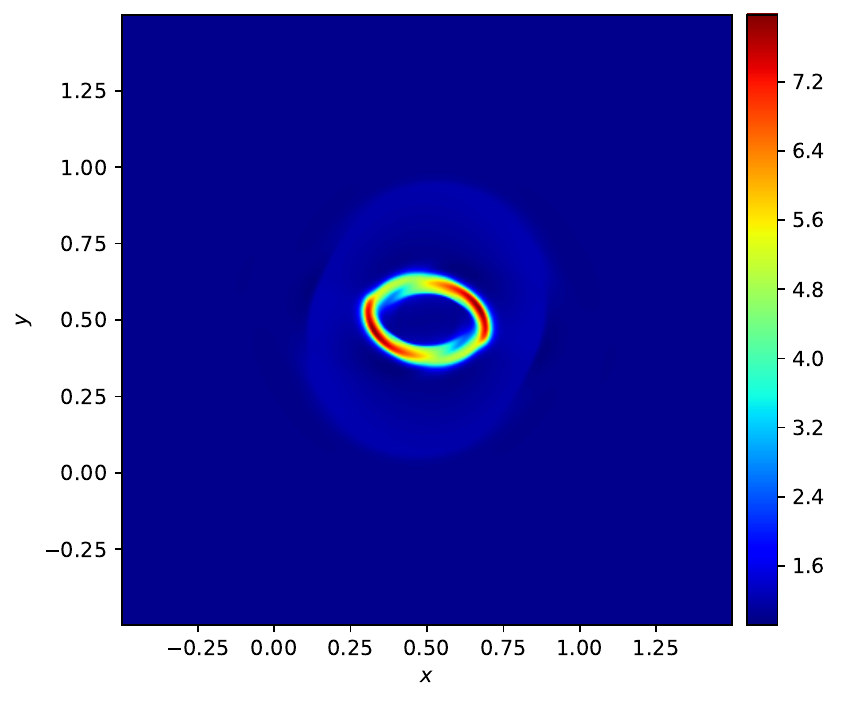}
			\label{fig:density_512_imp}}
		\subfigure[Total pressure $p_I +p_E$ with \textbf{O2EXP-MultiD} scheme.]{
			\includegraphics[width=0.45\textwidth,clip=]{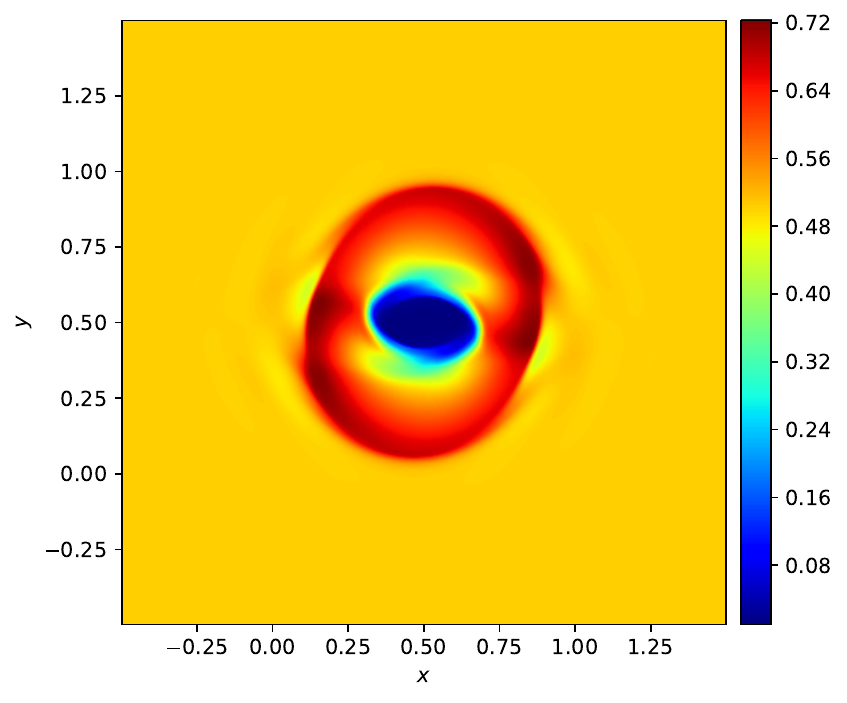}
			\label{fig:p_512_exp}}
			\subfigure[Total pressure $p_I +p_E$ with \textbf{O2IMEX-MultiD} scheme.]{
			\includegraphics[width=0.45\textwidth,clip=]{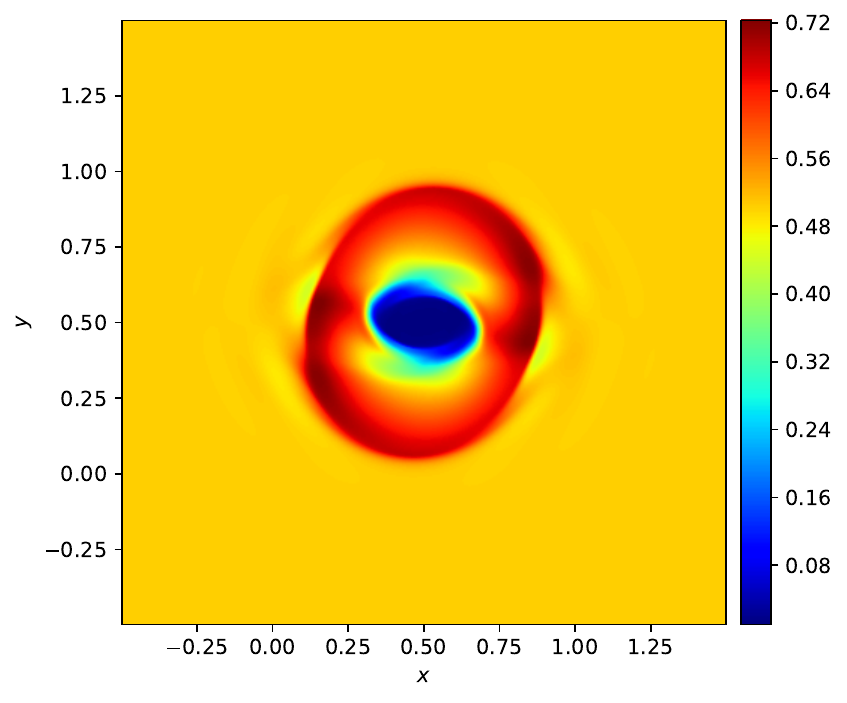}
			\label{fig:p_1024}}
		\caption{\nameref{test:2d_rotor}: Plots of total density and total pressure with $512\times 512$ cells at time $t = 0.295$.}
		\label{fig:rotor_o2}
	\end{center}
\end{figure}

\begin{figure}[!htbp]
	\begin{center}
        \subfigure[Pressure plot using $\kappa=\xi=0.5$ for \textbf{O2EXP-PHM} and \textbf{O2IMEX-PHM} schemes.]{
			\includegraphics[width=0.45\textwidth,clip=]      {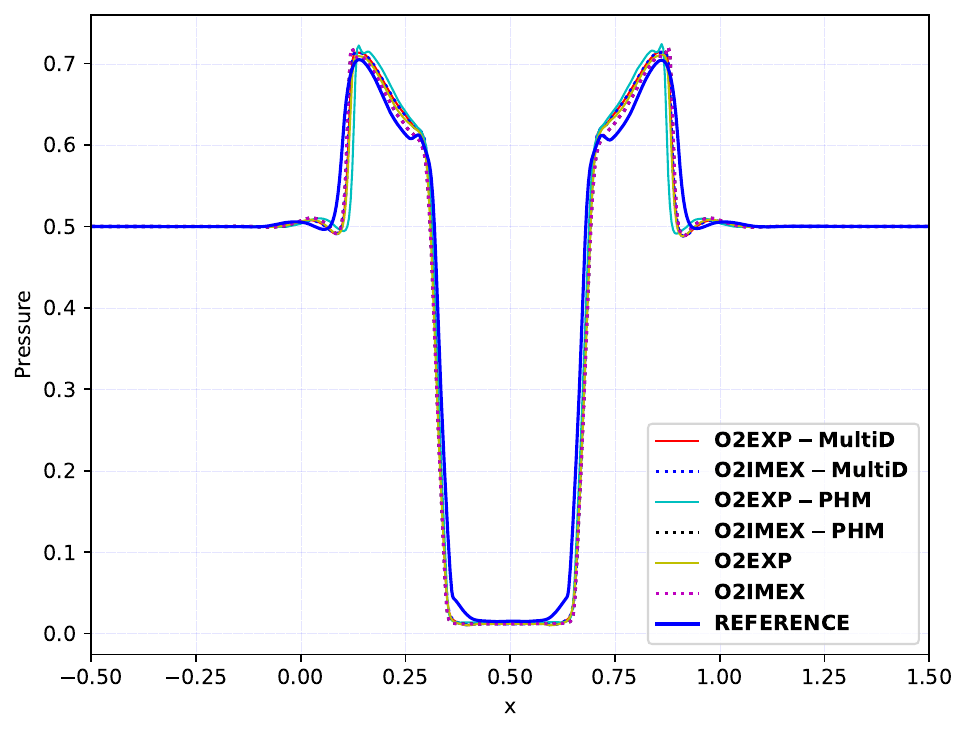}
			\label{fig:rotor_pressure_0p5}}
       \subfigure[Pressure plot using $\kappa=\xi=1$ for \textbf{O2EXP-PHM} and \textbf{O2IMEX-PHM} schemes.]{
			\includegraphics[width=0.45\textwidth,clip=]{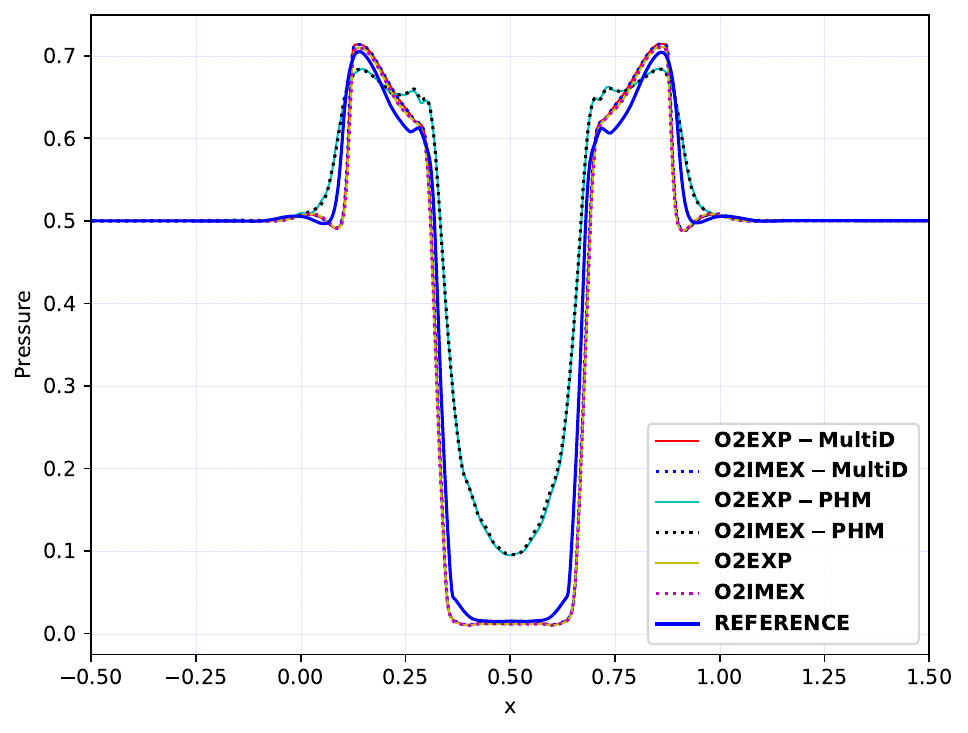}
			\label{fig:rotor_presssure_1p0}}
			%\label{fig:rotor_o2_pressure}
		\caption{\nameref{test:2d_rotor}: Plot of pressure cut along $x = 0.5$ for different schemes and using different values of $\kappa$ and $\xi$.}
		\label{fig:rotor_pressure_cut}
	\end{center}
\end{figure}

To design a test case consisting of strong shocks, we generalize the MHD rotor problem~\cite{balsara1999,fuchs2011,toth2000} to two-fluid plasma flow equations. We take a mass ratio of $m_I/m_E=25$ and consider the computational domain $[-0.5, 1.5]\times[-0.5, 1.5]$ with Neumann boundary conditions. We define  $r(x,y) = |(x,y) - (0.5, 0.5)|$ and  $f(r) = \frac{0.115 - r}{0.015}$. Using these, the initial data is given by,
\begin{itemize}
	\item For $r<0.1$ \[
\begin{cases}
	 \rho_I  = 10\frac{m_I}{(m_I+m_E)}, \\
	 \rho_E  = 10\frac{m_E}{m_I+m_E},\\
	u_I^x = u_E^x =-(y - 0.5)/0.1, \\
    u_I^y = u_E^y  =(x - 0.5)/0.1.
\end{cases}
\]
\item For $r >  0.115$
\[
\begin{cases}
	\rho_I  = \frac{m_I}{(m_I+m_E)}, \\
	\rho_E  = \frac{m_E}{(m_I+m_E)}, \\
	u_I^x = u_E^x =0, \\
	u_I^y = u_E^y =0.	
\end{cases}
\]
\item Otherwise, \[
\begin{cases}
	\rho_I  = \frac{m_I}{(m_I+m_E)}(1 + 9f(r)), \\
	\rho_E  = \frac{m_E}{(m_I+m_E)}(1 + 9f(r)), \\
	u_I^x = u_E^x= -f(r) \dfrac{(y - 0.5)}{r}, \\
	u_I^y = u_E^y =f(r) \dfrac{(x - 0.5)}{r}.
\end{cases}
\]
\end{itemize}
The other variables are set to,
\[
(u_I^z, u_E^z, B_x, B_y, B_z, p_I, p_E) = 
\left(0, 0, \frac{2.5}{\sqrt{4\pi}}, 0, 0, 0.5, 0.5\right),
\]
where $r(x, y)$ and $f(r)$ are same as defined earlier. The charge-to-mass ratios are taken to be $r_I=25.495097$, $r_E = -637.377439$. We compute the results till the final time $t=0.295$. In Figure~\ref{fig:rotor_o2}, we have plotted total density and total pressure for the {\bf O2EXP-MultiD} and {\bf O2IMEX-MultiD} schemes. We observe that both schemes are able to resolve the rotating fluid structures. In Figure \ref{fig:rotor_pressure_cut}, we have plotted the pressure cut along the line $x=0.5$ for all the schemes for $\kappa=\xi=0.5$ and $\kappa=\xi=1.0$. The reference solutions are also plotted which are using the second-order IMEX scheme (O2-222-imex) from \cite{Abgrall2014} with $1024\times 1024$ cells for $\kappa=\xi=0.5$. We observe that \textbf{O2EXP-PHM} and \textbf{O2IMEX-PHM} schemes have similar performance as the other schemes for $\kappa=\xi=0.5$. However, for $\kappa=\xi=1.0$ they are more diffusive when compared to the other schemes.  

\begin{figure}[!htbp]
	\begin{center}
			\subfigure[$\|\na \cdot \mathbf{B}^n\|_1$,  and $\|\na \cdot \mathbf{B}^n\|_2$  errors for explicit schemes  \textbf{O2EXP-MultiD}, \textbf{O2EXP-PHM} and \textbf{O2EXP}.]{
			\includegraphics[width=0.45\textwidth,clip=]{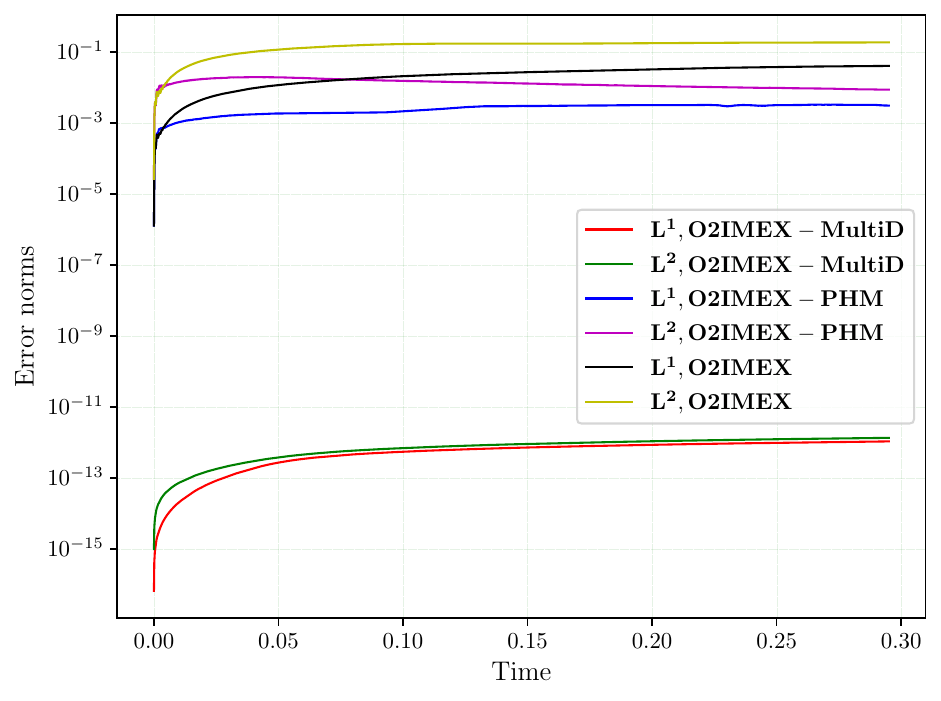}
			\label{fig:rotor_divB_exp}}
		\subfigure[$\|\na \cdot \mathbf{B}^n\|_1$,  and $\|\na \cdot \mathbf{B}^n\|_2$  errors for IMEX schemes \textbf{O2IMEX-MultiD}, \textbf{O2IMEX-PHM} and \textbf{O2IMEX}.]{
			\includegraphics[width=0.45\textwidth,clip=]{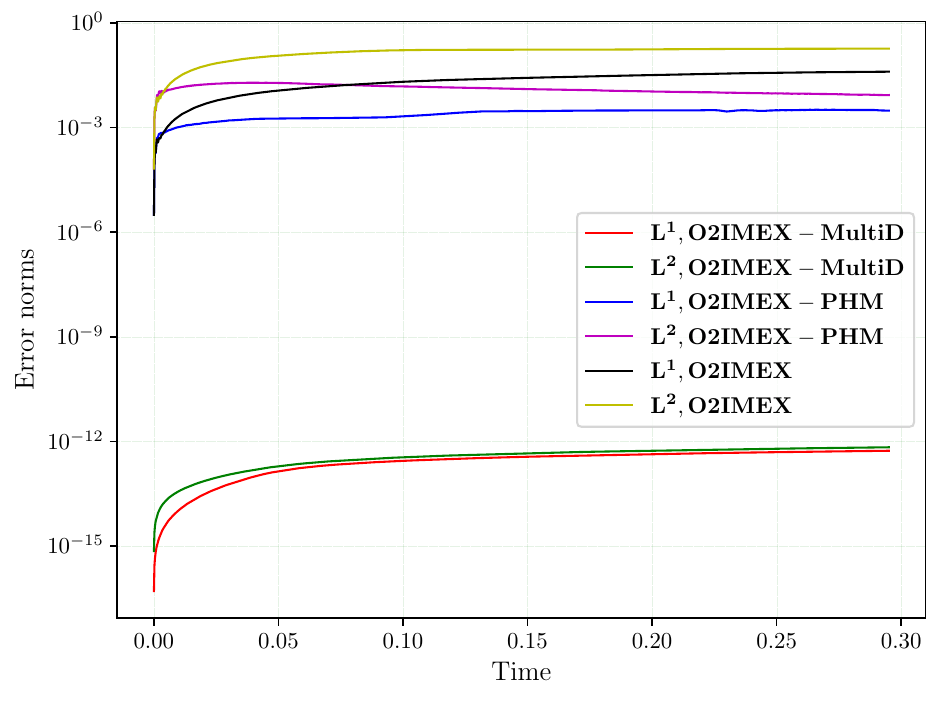}
			\label{fig:rotor_divB}}
		\subfigure[$\|\na\cdot\Eb\|_{1}^{E},$ and $\|\na\cdot\Eb\|_{2}^{E}$ errors, for explicit schemes  \textbf{O2EXP-MultiD}, \textbf{O2EXP-PHM} and \textbf{O2EXP}. ]{
	\includegraphics[width=0.45\textwidth,clip=]{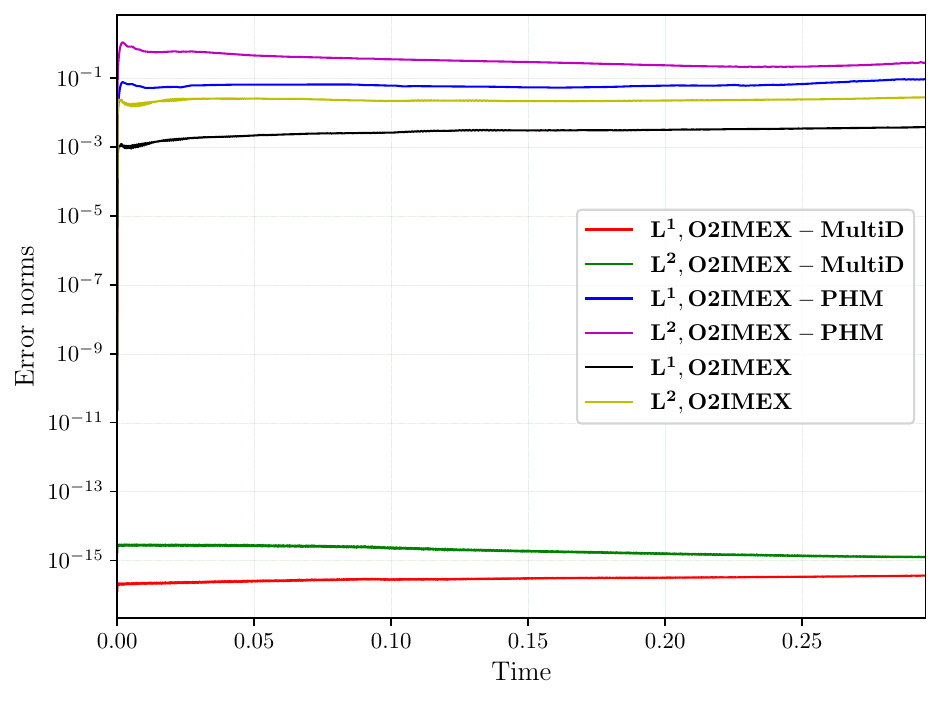}
	\label{fig:rotor_divE_exp}}	
		\subfigure[$\|\na\cdot\Eb\|_{1}^{I},$ and $\|\na\cdot\Eb\|_{2}^{I}$ errors, for IMEX schemes \textbf{O2IMEX-MultiD}, \textbf{O2IMEX-PHM} and \textbf{O2IMEX}.]{
			\includegraphics[width=0.45\textwidth,clip=]{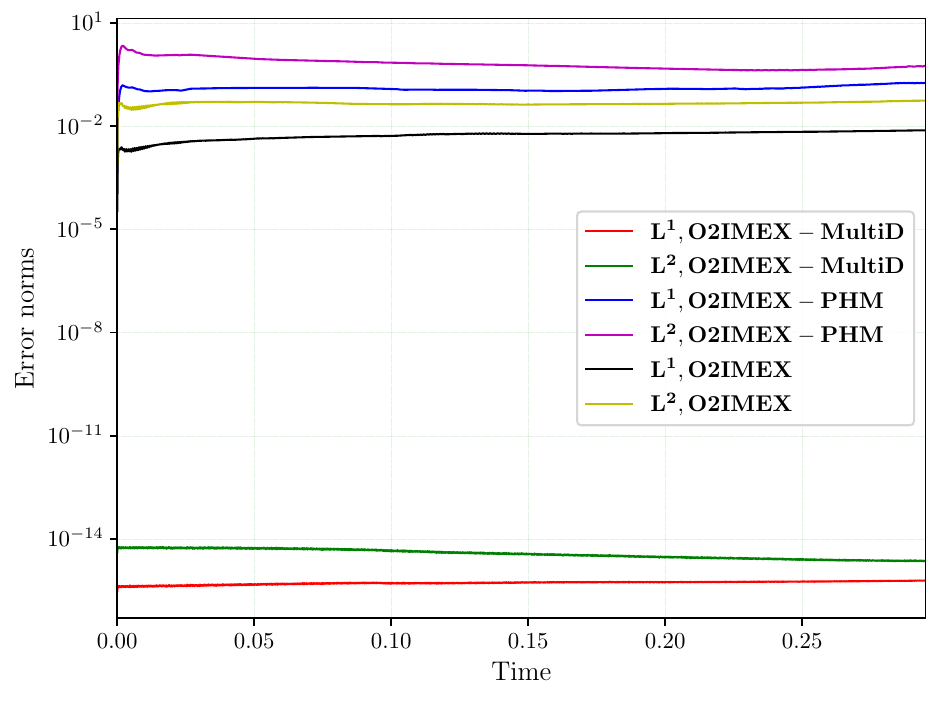}
			\label{fig:rotor_divE}}
		\caption{\nameref{test:2d_rotor}: Evolution of the divergence constraint errors for \textbf{O2IMEX-MultiD},  \textbf{O2EXP-MultiD}, \textbf{O2IMEX-PHM}, \textbf{O2EXP-PHM},  \textbf{O2IMEX} and \textbf{O2EXP} schemes using $512\times 512$ cells.}
		\label{fig:rotor_div_norms}
	\end{center}
\end{figure}

\begin{figure}[!htbp]
	\begin{center}
		\includegraphics[width=3.0in, height=2.3in]{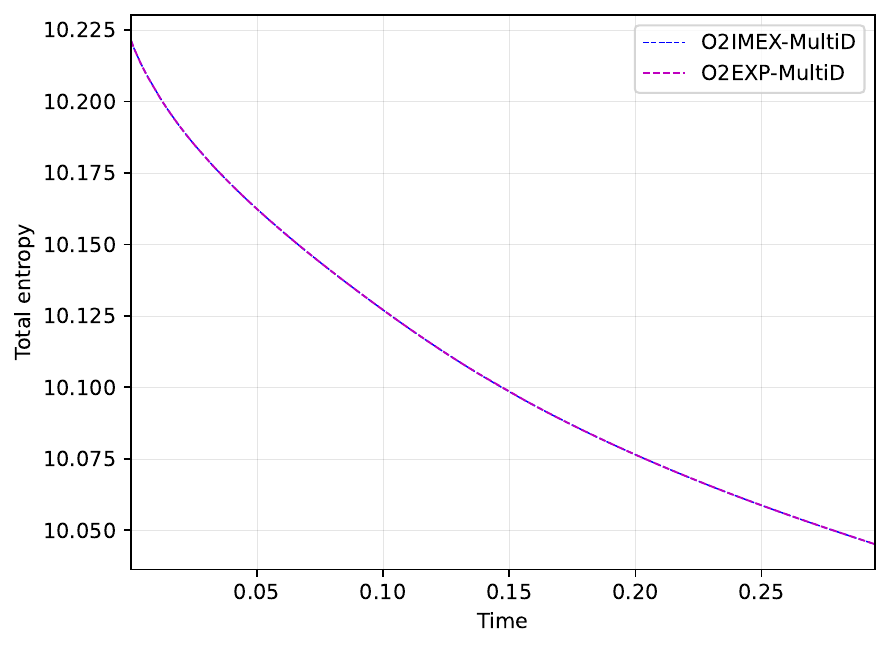}
		\caption{\nameref{test:2d_rotor} Total fluid entropy evolution for the schemes {\bf O2EXP-MultiD} and {\bf O2IMEX-MultiD}.}
		\label{fig:rotor_entropy}
	\end{center}
\end{figure}

In Figure~\ref{fig:rotor_div_norms}, we have plotted the time evolution of the divergence errors. We again observe that {\bf O2EXP-MultiD} and {\bf O2IMEX-MultiD} schemes have significantly low divergence errors for both electric and magnetic field constraints. 
In Figure \ref{fig:rotor_entropy}, we have plotted the time evolution of the total fluid entropy for both schemes. We note that both schemes have similar entropy decays.
% Variables Plot Pseudo

%------------------
%\todo[inline]{Many figures are stacked in a single column, they are going out of page. Rearrange, maybe two columns.}\todo{They look fine when compiling online. We will check this carefully before final submission}
%------------------
\subsection{Two-fluid GEM challenge problem} \label{test:2d_gem}
In this test case, we consider the non-relativistic Geospace Environment Modeling (GEM) magnetic reconnection problem given in~\cite{Birn2001,Hakim2006,loverich2011,wang2020}. The computational domain is given by $[-L_x/2,L_x/2] \times [-L_y/2,L_y/2]$ where $L_x=8\pi$ and $L_y=4\pi$. We consider periodic boundary conditions at $x=\pm L_x/2$ and conducting wall boundary at $y=\pm L_y/2$. The ion-electron mass ratio is taken to be $m_I/m_E = 25$. Accordingly, we take $r_I=1.0$ and $r_E=-25.0.$. The initial conditions are,
\begin{align*}
	\begin{pmatrix}
		\rho_I  \\
		u_I^z\\
		p_I     \\ \\
		\rho_E  \\
		u_E^z \\
		p_E     \\ \\
		B_x     \\
		B_y
	\end{pmatrix} =
	\begin{pmatrix}
		n                                                          \\
		0.0               \\
		5n \frac{B_0}{12}\\ \\
		\frac{m_E}{m_I} n                                          \\
		\frac{J_{e}^z}{r_E\rho_E}                                                   \\
		p_I/5                                                      \\ \\
		B_0 \tanh(y/\lambda) - \psi_0 \frac{\pi}{L_y} \cos(\frac{2\pi x}{L_x}) \sin(\frac{\pi y}{L_y})
		\\
		\psi_0 \frac{2 \pi}{L_x} \sin(\frac{2 \pi x}{L_x}) \cos(\frac{\pi y}{L_y})
	\end{pmatrix},
\end{align*}
where, $n=\mathrm{sech}^2(y/\lambda)+0.2$ and $J_{e}^z = -(B_0/\lambda )\mathrm{sech}^2(y/\lambda),\ \lambda = 0.5,\ B_0 =1.0$. All other variables are set to zero. Adiabatic gas constant are taken to  be $\gamma_I =\gamma_E= 5.0/3.0$ and we compute the solutions till time $t=40.$

\begin{figure}[!htbp]
	\begin{center}
		\subfigure[Total density $\rho_I +\rho_E$ with \textbf{O2EXP-MultiD} scheme.]{
			\includegraphics[width=0.45\textwidth,clip=]{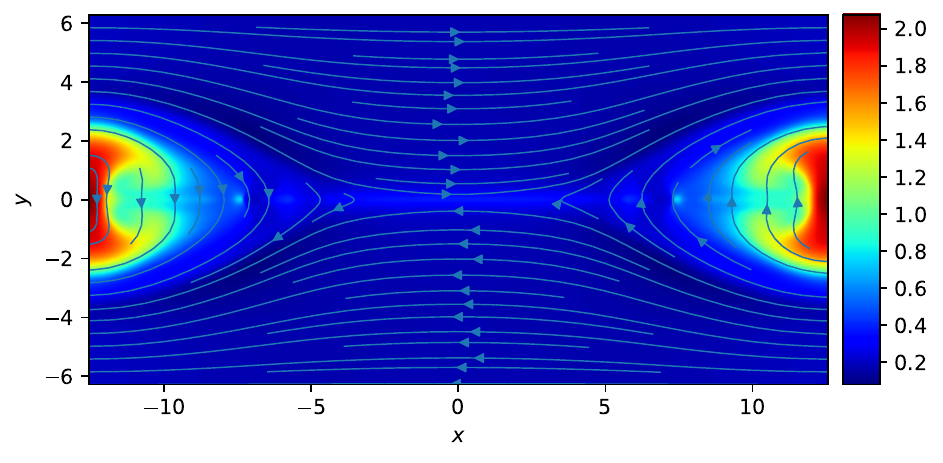}
			\label{fig:gem_exp_t25_512_rho_multid}}
		\subfigure[Total density $\rho_I +\rho_E$ with \textbf{O2IMEX-MultiD} scheme.]{
			\includegraphics[width=0.45\textwidth,clip=]{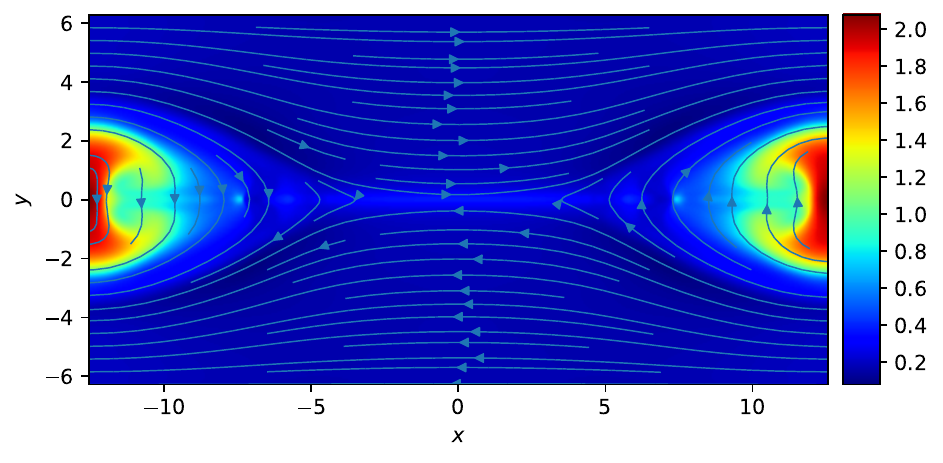}
			\label{fig:gem_imp_t25_512_rho_multid}}
		\subfigure[$B_z$ with \textbf{O2EXP-MultiD} scheme.]{
			\includegraphics[width=0.45\textwidth,clip=]{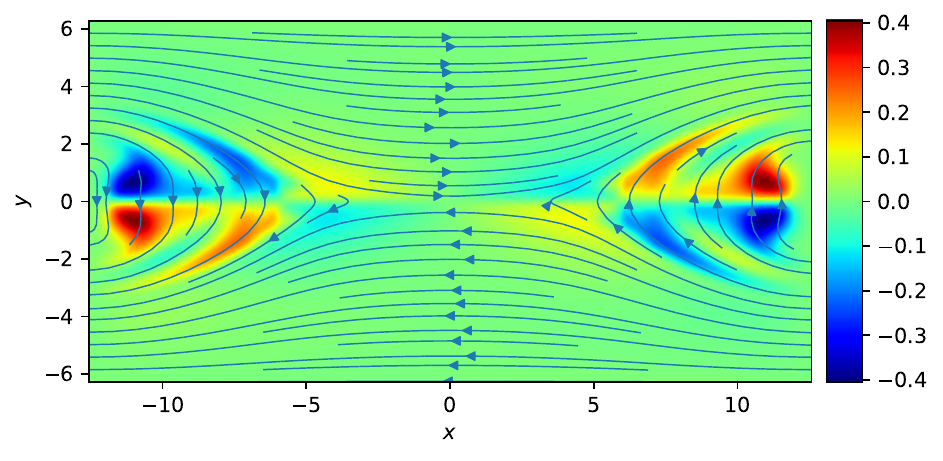}
			\label{fig:gem_exp_t25_512_Bz_multid}}
		\subfigure[$B_z$ with \textbf{O2IMEX-MultiD} scheme.]{
			\includegraphics[width=0.45\textwidth,clip=]{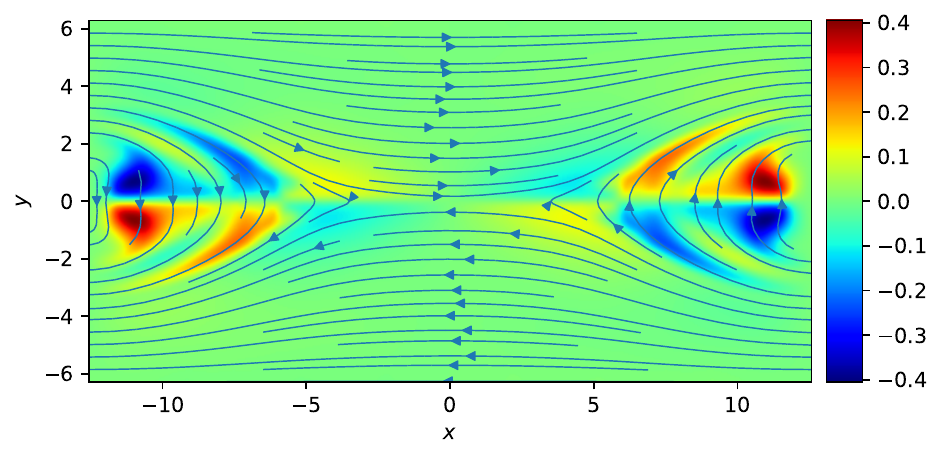}
			\label{fig:gem_imp_t25_512_Bz_multid}}
		\subfigure[$u^x_I$ with \textbf{O2EXP-MultiD} scheme.]{
			\includegraphics[width=0.45\textwidth,clip=]{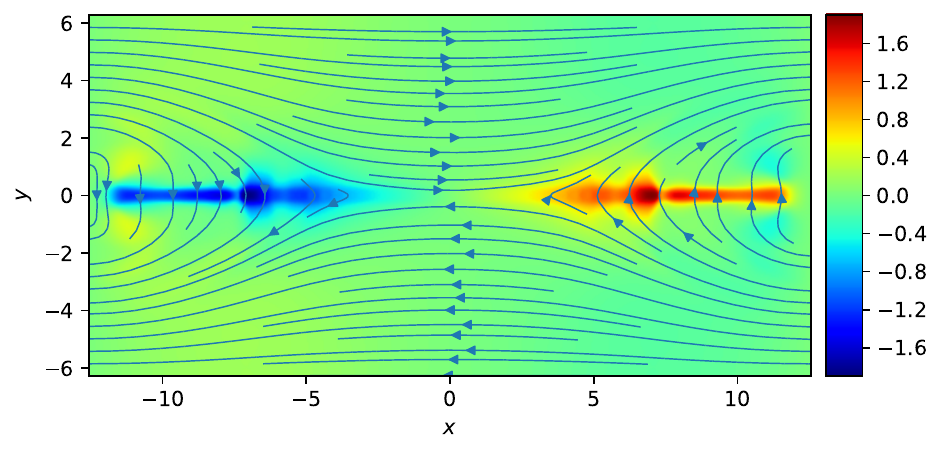}
			\label{fig:gem_exp_t25_512_uxi_multid}}
		\subfigure[$u^x_I$ with \textbf{O2IMEX-MultiD} scheme.]{
			\includegraphics[width=0.45\textwidth,clip=]{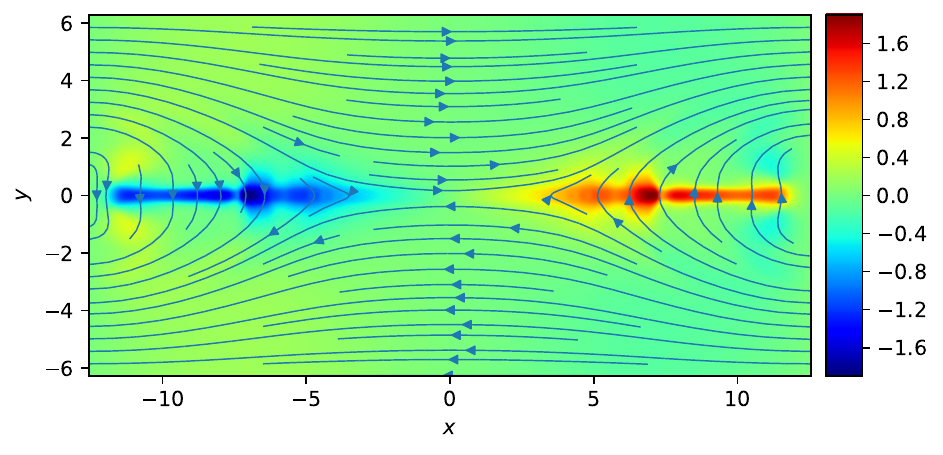}
			\label{fig:gem_imp_t25_512_uxi_multid}}
		\subfigure[$u^x_E$ with \textbf{O2EXP-MultiD} scheme.]{
			\includegraphics[width=0.45\textwidth,clip=]{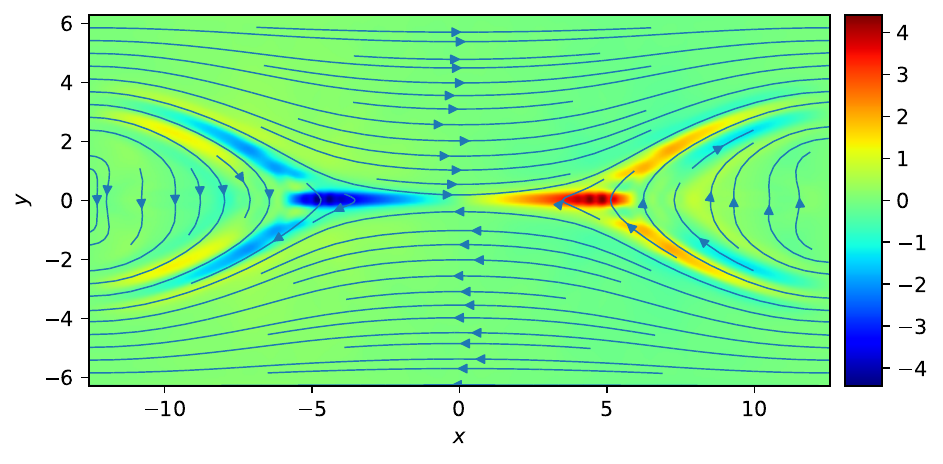}
			\label{fig:gem_exp_t25_512_uxe_multid}}
		\subfigure[$u^x_E$ with \textbf{O2IMEX-MultiD} scheme]{
			\includegraphics[width=0.45\textwidth,clip=]{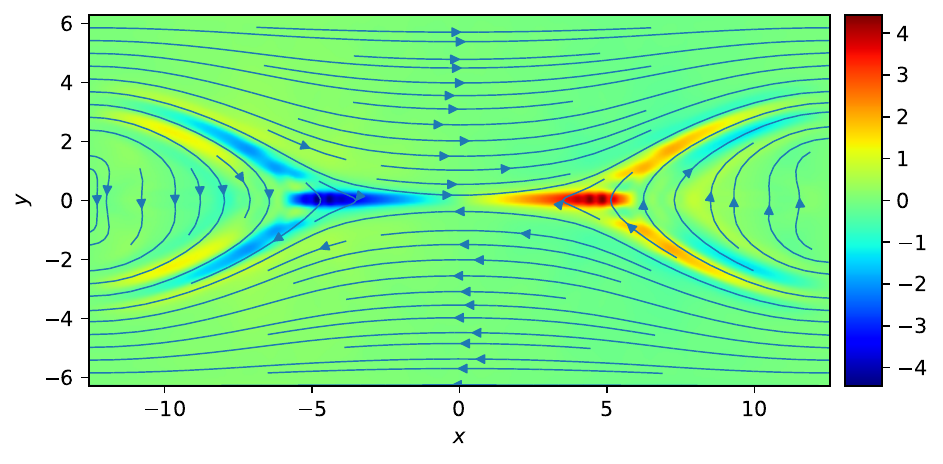}
			\label{fig:gem_imp_t25_512_uxe_multid}}
		\caption{\nameref{test:2d_gem}: Plots for the total density, $B_z$, $u_I^x$, and  $u^x_E$ with $512 \times 256$ cells at time $t=25.0$.}
		\label{fig:gem_o2_25_multiD}
	\end{center}
\end{figure}

\begin{figure}[h]
\begin{center}
	\includegraphics[width=0.55\textwidth,clip=]{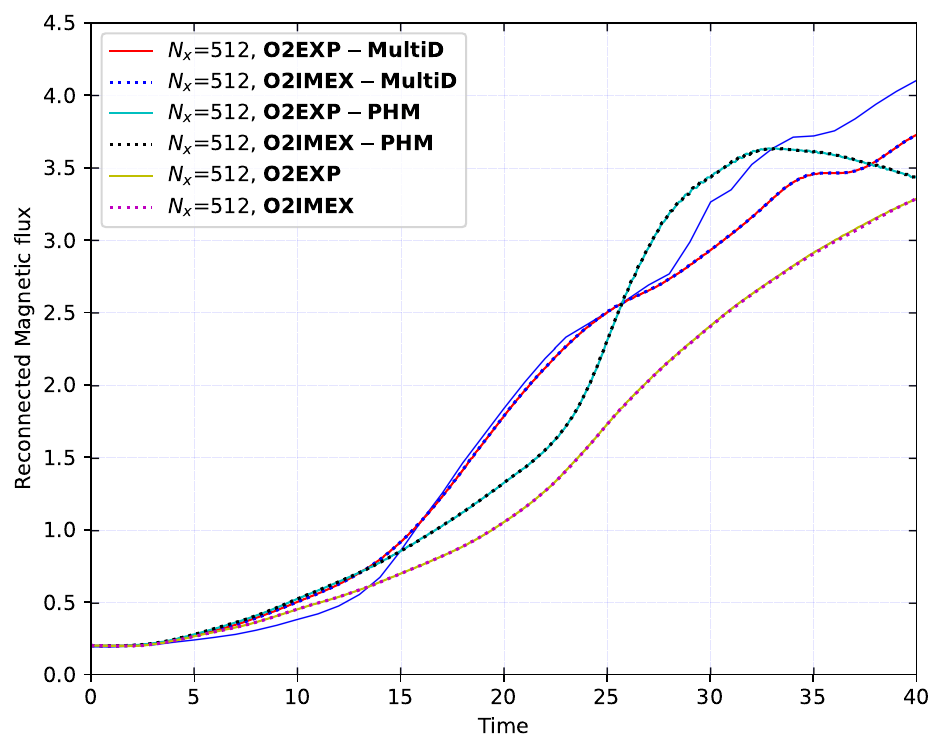}
	\caption{\nameref{test:2d_gem}: Time development of the reconnected magnetic flux using $512\times256$ cells, for different schemes. We overlay the plot on reconnection rates from \cite{Hakim2006} (solid blue lines). }
	\label{fig:reconn_rate}
\end{center}
\end{figure}
\begin{figure}[!htbp]
	\begin{center}
		\includegraphics[width=0.55\textwidth,clip=]{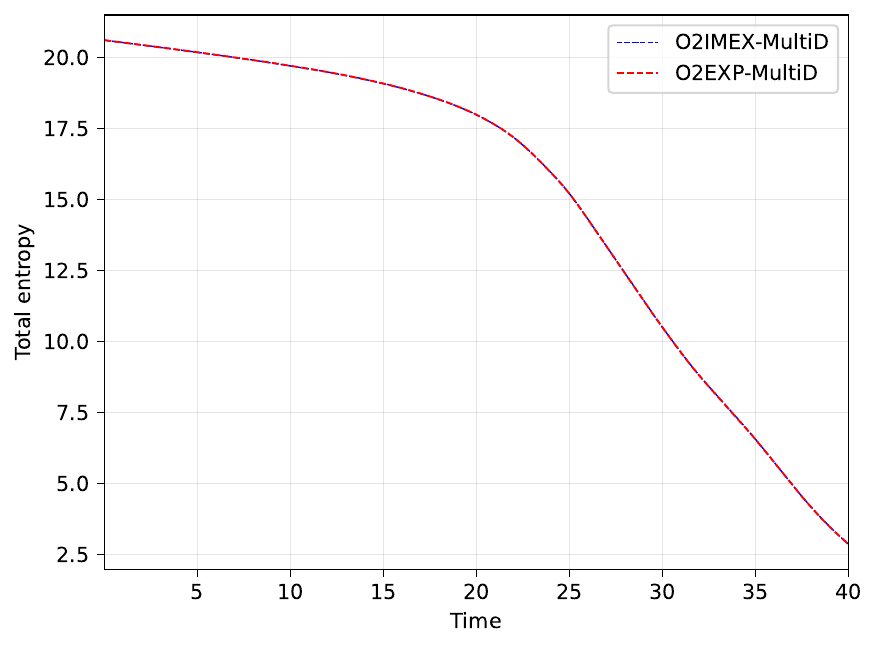}
		\caption{\nameref{test:2d_gem}: Total fluid entropy evolution for the schemes {\bf O2EXP-MultiD} and {\bf O2IMEX-MultiD}.}
		\label{fig:GEM_entropy}
	\end{center}
\end{figure}
\begin{figure}[!htbp]
	\begin{center}
		\subfigure[$\|\na \cdot \mathbf{B}^n\|_1$,  and $\|\na \cdot \mathbf{B}^n\|_2$  errors for explicit schemes  \textbf{O2EXP-MultiD}, \textbf{O2EXP-PHM} and \textbf{O2EXP}.]{
			\includegraphics[width=0.45\textwidth,clip=]{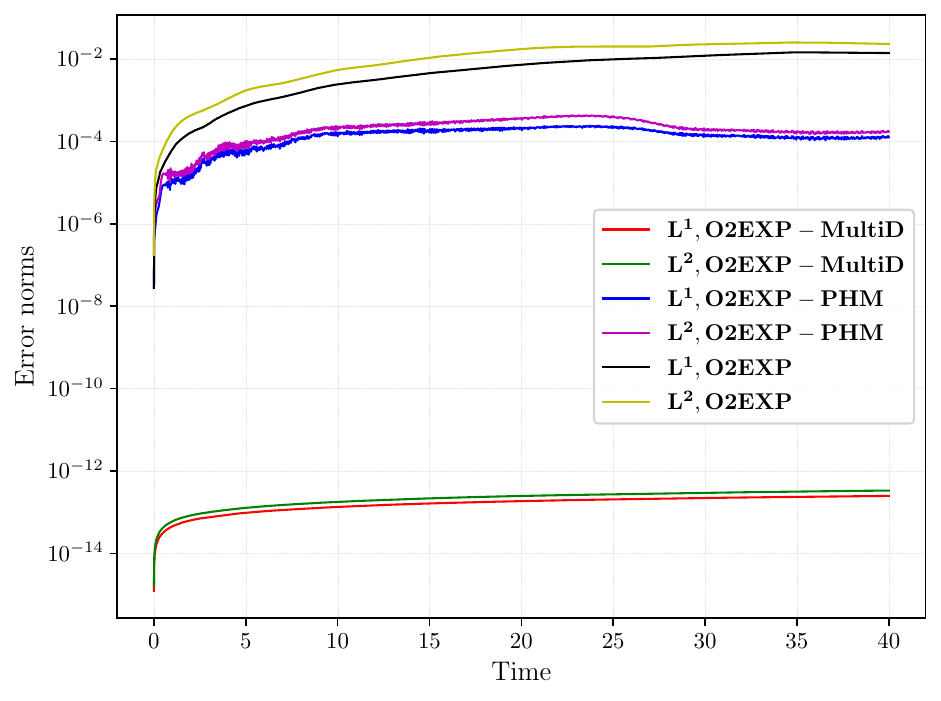}
			\label{fig:gem_divB_exp}}
		\subfigure[$\|\na \cdot \mathbf{B}^n\|_1$,  and $\|\na \cdot \mathbf{B}^n\|_2$  errors for IMEX schemes \textbf{O2IMEX-MultiD}, \textbf{O2IMEX-PHM} and \textbf{O2IMEX}.]{
			\includegraphics[width=0.45\textwidth,clip=]{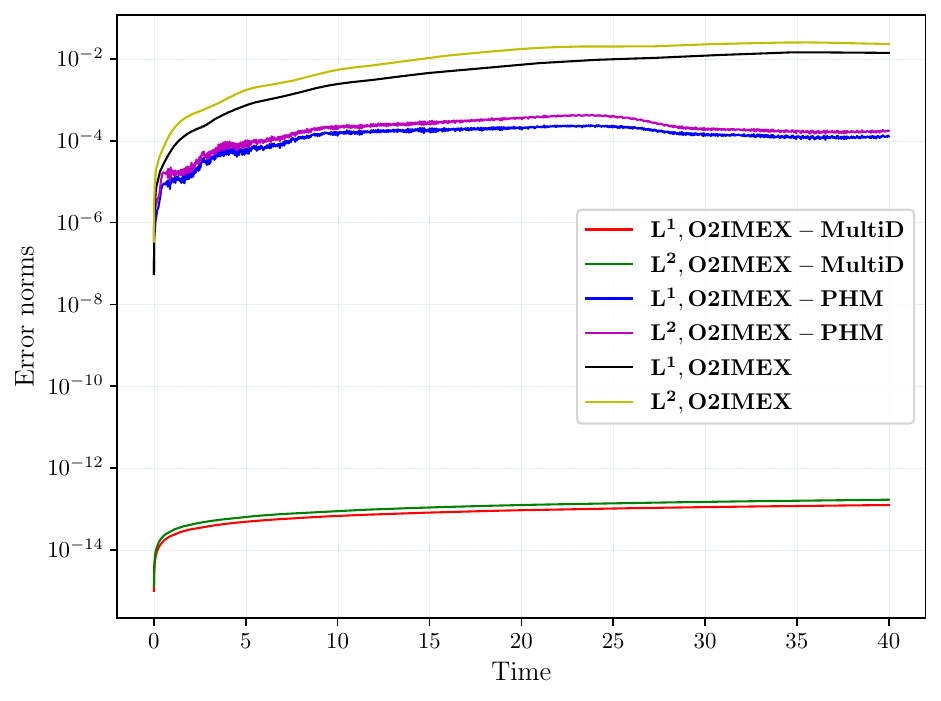}
			\label{fig:gem_divB}}
		\subfigure[$\|\na\cdot\Eb\|_{1}^{E},$ and $\|\na\cdot\Eb\|_{2}^{E}$ errors, for explicit schemes  \textbf{O2EXP-MultiD}, \textbf{O2EXP-PHM} and \textbf{O2EXP}.]{
			\includegraphics[width=0.45\textwidth,clip=]{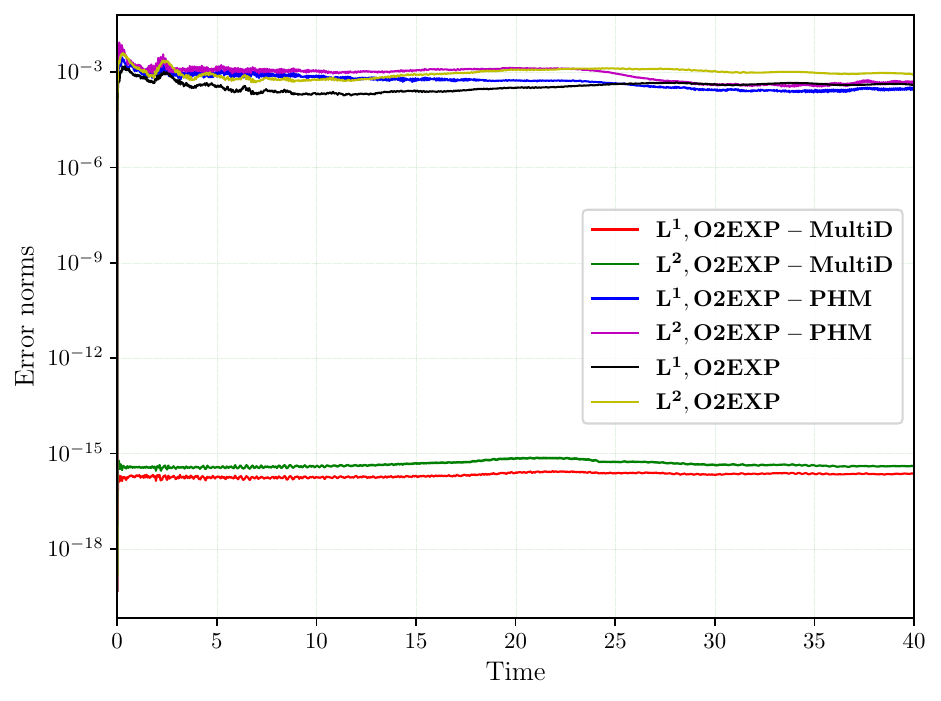}
			\label{fig:gem_divE_exp}}
		\subfigure[$\|\na\cdot\Eb\|_{1}^{I},$ and $\|\na\cdot\Eb\|_{2}^{I}$ errors, for IMEX schemes \textbf{O2IMEX-MultiD}, \textbf{O2IMEX-PHM} and \textbf{O2IMEX}.]{
			\includegraphics[width=0.45\textwidth,clip=]{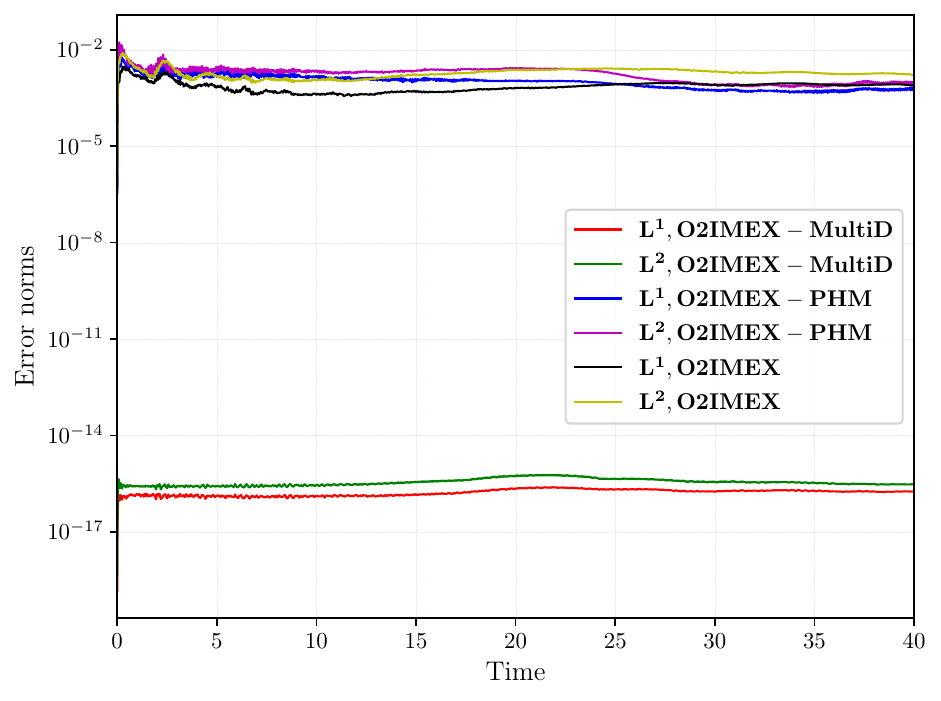}
			\label{fig:gem_divE}}
		\caption{\nameref{test:2d_gem}: Evolution of the divergence constraint errors for \textbf{O2IMEX-MultiD},  \textbf{O2EXP-MultiD}, \textbf{O2IMEX-PHM}, \textbf{O2EXP-PHM},  \textbf{O2IMEX} and \textbf{O2EXP} schemes using $512\times 256$ cells.}
		\label{fig:gem_div_norms}
	\end{center}
\end{figure}	
%Reconnection Rate

In Figure~\ref{fig:gem_o2_25_multiD}, we have plotted total density, $z$-component of the magnetic field, $x$-component of ion velocity, and $x$-component of electron velocity for \textbf{O2EXP-MultiD} and \textbf{O2IMEX-MultiD} schemes using a mesh of $512\times 256$ cells at time $t=25$. The plots also contain field lines for $(B_x,B_y)$. We also observe that the magnetic reconnection is well underway. 
In Figure~\ref{fig:reconn_rate}, we have plotted the magnetic reconnection rates as measured by the quantity
$$\dfrac{1}{2 B_0} \int_{-L_x/2}^{L_x/2}
| B_y(x,y=0,t) | dx,
$$
for all the schemes. We observe that the corresponding explicit and IMEX schemes have similar magnetic reconnection rates for all the different discretizations of Maxwell's equations. We also note that \textbf{O2EXP-MultiD} and \textbf{O2IMEX-MultiD} have magnetic reconnection rates close to the result in \cite{Hakim2006}, which is also shown in the figure as a blue solid line.

In Figure~\ref{fig:gem_div_norms}, we have plotted the evolution of the divergence errors. We again observe that the proposed scheme \textbf{O2EXP-MultiD} and \textbf{O2IMEX-MultiD} have almost machine precision divergence errors throughout the simulations. Errors for the schemes \textbf{O2EXP-PHM} and \textbf{O2IMEX-PHM} are slightly better for the magnetic field divergence errors when compared with  \textbf{O2EXP} and \textbf{O2IMEX} schemes, but the electric field divergence errors are closer to each other.

In Figure~\ref{fig:GEM_entropy}, we have plotted the time evolution of total fluid entropy for both schemes. Initially, the flow is relatively smooth; hence, not much entropy decay is observed. However, around $t=20$, the flow becomes more irregular, and we see both schemes decaying entropy faster. Both schemes have similar entropy decay performance.

\section{Conclusions}
\label{sec:con}
In this article, we have proposed co-located numerical discretizations of the two-fluid plasma flow equations, which are consistent with the divergence constraints of Maxwell's equations. It is based on using a multidimensional Riemann solver at the vertices to define numerical fluxes across the cell faces. The proposed schemes can also be coupled with numerical schemes to preserve other stability properties related to the fluid variables, like entropy stability. We have also presented several test cases in one- and two-dimensions to demonstrate the consistency with the divergence constraints. We have also compared our results with the other strategies commonly employed to preserve the divergence constraints and show that the proposed schemes are superior at preserving the divergence constraints and give more accurate solutions.

% In this article, we propose a novel method of solving the two-fluid plasma flow equations using co-located numerical discretizations. We particularly highlight that our approach guarantees agreement with the divergence limitations present in Maxwell's equations. The fundamental idea of our method is to specify numerical fluxes across cell faces precisely by using a multidimensional Riemann solver at the vertices.

% Our suggested approaches stand out for being compatible with other numerical techniques meant to maintain fluid variable stability features. We demonstrate the effectiveness of our method in upholding consistency with the divergence constraints by means of comprehensive testing in both one- and two-dimensional scenarios.

% We verify our approach by contrasting our findings with standard techniques that are frequently applied to maintain divergence limitations. Our results show that the suggested approaches work better than expected, indicating that they are more effective in preserving the integrity of the divergence constraints. 

%\bibliographystyle{plainnat}
\begin{acknowledgements}
The work of Harish Kumar is supported in parts by VAJRA grant No. VJR/2018/000129 by the Dept. of Science and Technology, Govt. of India. Harish Kumar and Jaya Agnihotri acknowledge the support of FIST Grant Ref No. SR/FST/MS-1/2019/45 by the Dept. of Science and Technology, Govt. of India. The work of Praveen Chandrashekar is supported by the Department of Atomic Energy, Government of India, under project no. 12-R\&D-TFR-5.01-0520. 
\end{acknowledgements}

\section*{Conflict of interest}
The authors declare that they have no conflict of interest.
\section*{Data Availability Declaration}
Data will be made available on reasonable request.

\appendix

\section{Discretization of Maxwell's equations}
\label{sec:Ez}
Define the difference and averaging operators
\[
\dx \phi(x,y) = \frac{\phi(x + \Delta x/2, y) - \phi(x - \Delta x/2, y)}{\Delta x}, \qquad \ax\phi(x,y) = \frac{1}{2}[ \phi(x-\Delta x/2,y) + \phi(x+\Delta x/2,y)]
\]
\[
\dy \phi(x,y) = \frac{\phi(x, y+\Delta y/2) - \phi(x, y - \Delta y/2)}{\Delta y}, \qquad \ay\phi(x,y) = \frac{1}{2}[ \phi(x, y-\Delta y/2) + \phi(x, y+\Delta y/2)]
\]
The discretization of the magnetic field given in Section~\ref{subsec:multid_riem_solver} can be written as
\begin{equation}
\label{eq:BxBy}
\frac{d B_{x,i,j}}{d t} + \dy\ax\tilde E_{z,i,j} = 0, \qquad
\frac{d B_{y,i,j}}{d t} - \dx\ay\tilde E_{z,i,j} = 0
\end{equation}
where $\tilde E_z$ is defined at the vertices. The divergence is measured at the vertices using a central difference approximation
\[
\nabla\cdot \Bb_{\iph,\jph} = (\dx\ay B_x + \dy\ax B_y)_{\iph,\jph}
\]
Then
\[
\frac{d}{dt}  \nabla\cdot \Bb_{\iph,\jph} = \dx\ay \frac{d B_x}{dt} + \dy\ax \frac{d B_y}{dt}  =  -\dx\ay \dy\ax \tilde E_z + \dy\ax \dx\ay \tilde E_z = 0
\]
The last result follows since the difference and averaging operators commute. This is the discrete analogue of the property that $\nabla \cdot \nabla \times \Bb = 0$ for any smooth vector field $\Bb$; such schemes which satisfy vector calculus identifies are usually referred to as being mimetic. Note that this property holds for any consistent definition of the vertex values $\tilde E_z$.

\paragraph{Construction of $\tilde E_z$}. We can motivate the formula \eqref{eq:multid_Ez} for $\tilde E_z$ as a generalization of Lax-Friedrich flux or by starting with the viscous problem
\[
\frac{d \Bb}{d t} + \nabla \times \Eb = -\nu \nabla \times \nabla \times \Bb
\]
where $\nu$ is a numerical viscosity. Here we have used the fact that $\Delta \Bb = -\nabla \times \nabla \times \Bb$ since $\nabla \cdot \Bb = 0$. In two dimensiona, this equation can be written in component form as
\[
\partial_t B_x + \partial_y E_z = -\nu \partial_y ( \partial_x B_y - \partial_y B_x), \qquad
\partial_t B_y - \partial_x E_z = +\nu \partial_x ( \partial_x B_y - \partial_y B_x)
\]
Let us discretize the above two PDEs with a multi-dimensional central difference scheme at each grid point $(i,j)$
\begin{eqnarray*}
\frac{d B_x}{dt} + \dy\ax \ax\ay E_z &=& -\nu \dy\ax ( \dx\ay B_y - \dy\ax B_x) \\
\frac{d B_y}{dt} - \dx\ay \ax\ay E_z &=& +\nu \dx\ay ( \dx\ay B_y - \dy\ax B_x) \\
\end{eqnarray*}
Comparing this with \eqref{eq:BxBy}, we see that the vertex value of $E_z$ is given by
\[
\tilde E_{z,\iph,\jph} = \left[ \ax\ay E_z + \nu ( \dx\ay B_y - \dy\ax B_x) \right]_{\iph,\jph}
\]
which agrees with \eqref{eq:multid_Ez} if the artificial viscosity coefficient is chosen, as in an upwind scheme, as $\nu = \frac{1}{2} c h$ where $h = \Delta x = \Delta y$.
\bibliographystyle{plainnat}

\bibliography{library}

\end{document}